\documentclass[titlepage,11pt]{article}

\usepackage{latexsym}
\usepackage{amsfonts}
\usepackage{amsmath}
\usepackage{float}
\usepackage{lmodern}
\usepackage{tikz}
\usetikzlibrary{decorations.pathmorphing}

\newcommand\blackslug{\hbox{\hskip 1pt \vrule width 4pt height 8pt depth 1.5pt
        \hskip 1pt}}
\newcommand\bbox{\hfill \quad \blackslug \bigbreak}
\def\d{\hbox{-}}
\def\c{\hbox{-}\cdots\hbox{-}}
\def\dom{\hbox{dom}}
\def\ll{,\ldots,}

\title{Induced subgraphs of graphs with large chromatic number.\\
V. Chandeliers and strings}
\author{Maria Chudnovsky\thanks{Supported by NSF grants DMS-1265803 and DMS-1550991.}\\
Princeton University, Princeton, NJ 08544
\\
\\
Alex Scott\\
Oxford University, Oxford, UK
\\
\\
Paul Seymour\thanks{Supported by ONR grant N00014-14-1-0084,
NSF grants DMS-1265563 and DMS-1800053, and AFOSR grant A9550-19-1-0187.}\\
Princeton University, Princeton, NJ 08544}

\date{}

\newtheorem{thm}{}[section]

\newcommand{\Proof}{\noindent{\bf Proof.}\ \ }

\begin{document}
\maketitle
\begin{abstract}
It is known that every graph of sufficiently large chromatic number and bounded clique number contains, as an induced subgraph,
a subdivision of any fixed forest, and a subdivision of any fixed cycle. Equivalently,
every forest is pervasive, and $K_3$ is pervasive, in the class of all graphs, where we say 
a graph $H$ is ``pervasive'' (in some class of graphs) if for all $\ell\ge 1$, every graph 
in the class of bounded clique number and sufficiently large
chromatic number has an induced subgraph that is a subdivision of $H$, in which every edge of $H$ is replaced by
a path of at least $\ell$ edges. 

Which other graphs are pervasive? It was  proved by Chalopin, Esperet, Li and Ossona de Mendez that every such graph  
is a ``forest of lanterns'': roughly, every block is a ``lantern'', a graph obtained from a tree by adding one extra vertex,
and there
are rules about how the blocks fit together. It is not known whether every forest of lanterns is pervasive in the 
class of all graphs; but in another paper two of us prove that all ``banana trees'' are pervasive, that is, multigraphs
obtained from a forest by adding parallel edges, thus generalizing the two results above. This paper contains the first half 
of the proof, which works for
any forest of lanterns, not just for banana trees.

Say a class of graphs is ``$\rho$-controlled'' if for every graph in the class, its chromatic number
is at most some function (determined by the class) of the largest chromatic number of a $\rho$-ball in the graph. In this paper 
we prove that for every $\rho\ge 2$, and for every $\rho$-controlled class, every forest of lanterns is pervasive in this class.

These results turn out particularly nicely when applied to string graphs. A ``chandelier'' is a special lantern,
a graph obtained from a tree
by adding a vertex adjacent to precisely the leaves of the tree.
A ``string graph'' is the intersection graph of a set of curves in the plane. There
are string graphs with clique number two and chromatic number arbitrarily large. We prove that the 
class of string graphs is $2$-controlled, and consequently every forest of lanterns is pervasive in this class;
but in fact something stronger is true, that every string graph of sufficiently large chromatic number and bounded clique
number contains each fixed chandelier as an induced subgraph (not just as a subdivision); and the same for most
forests of chandeliers (there is an extra condition on how the blocks are attached together).
\end{abstract}

\section{Introduction}

All graphs in this paper are finite and simple, and if $G$
is a graph, $\chi(G)$ denotes its chromatic number, and $\omega(G)$ denotes its clique number, that is, the cardinality
of the largest clique of $G$. 
This is the fifth in a series of papers on the induced subgraphs that must be present in graphs that have bounded 
clique number and (sufficiently) large chromatic number. 
The series was originally motivated by three conjectures of Gy\'arf\'as from 
1985~\cite{gyarfas} concerning the lengths of induced cycles in such graphs:
\begin{thm}\label{oddholes}
For every integer $k\ge 0$, every graph $G$ with $\omega(G)\le k$ and $\chi(G)$ sufficiently large contains
an induced cycle of odd length at least $5$.
\end{thm}
\begin{thm}\label{longholes}
For all integers $k,\ell\ge 0$, every graph $G$ with $\omega(G)\le k$ and $\chi(G)$ sufficiently large contains
an induced cycle of length at least $\ell$.
\end{thm}
\begin{thm}\label{longoddholes}
For all integers $k,\ell\ge 0$, every graph $G$ with $\omega(G)\le k$ and $\chi(G)$ sufficiently large contains
an induced odd cycle of length at least $\ell$.
\end{thm}
All three conjectures have now been proved, 
in~\cite{oddholes, longholes,longoddholes} respectively.  Indeed, two of us~\cite{residues} have subsequently proved a much stronger theorem that
contains all these results: 
\begin{thm}\label{residues}
For all integers $k,\ell,m\ge 0$, every graph $G$ with $\omega(G)\le k$ and $\chi(G)$ sufficiently large contains
an induced cycle of length $\ell$ modulo $m$.
\end{thm}

In this paper we we will be interested in proving analogous results for induced subgraphs other than cycles.  In particular, we will be concerned with generalizing \ref{longholes}  (the other results above involve parity constraints and the methods we use
here do not work).

If $G$ has bounded clique number and very large chromatic number, which graphs $H$ must be present in $G$ as  
induced subgraphs? No graph $H$ has this property except for forests, because $G$ can have arbitrarily large girth; and it
is an open conjecture of Gy\'arf\'as~\cite{gyarfastree} and Sumner~\cite{sumner} that forests do have this property. 
This is an interesting question but we have nothing to say about it here (except that we will prove it for string graphs); 
we return to this problem in \cite{distantstars} and \cite{newbrooms}.

We may ask instead for the graphs $H$ with the property that 
every graph $G$ with bounded clique number
and sufficiently large chromatic number must contain an induced subgraph which is a {\em subdivision} of $H$.
This certainly yields a larger class of graphs;
for instance, every cycle has this property, in view of \ref{longholes}, and so does every forest, 
by the following theorem of~\cite{scott}: 

\begin{thm}\label{forests}
For every integer $k$ and every forest $F$, every graph $G$ with $\omega(G)\le k$ and $\chi(G)$ sufficiently large contains
an induced subdivision of $F$.
\end{thm}

This paper is concerned with subdivisions of a graph, so let us clarify some definitions before we go on.
Let $H$ be a graph, and let $H'$ be a graph obtained from $H$ by replacing each edge $uv$ by a path (of length at least one) 
joining $u,v$, such that these paths are vertex-disjoint except for their ends. We say that $H'$ is a {\em subdivision} of $H$; and it is
a {\em proper} subdivision of $H$ if all the paths have length at least two.
If each of the paths has exactly $\ell+1$ edges we call it an {\em $\ell$-subdivision}; if they each have at least $\ell+1$ edges
it is an {\em ($\ge \ell$)-subdivision}; and if they all have at most $\ell+1$ it is an {\em ($\le \ell$)-subdivision}.
If they all have length at least two and at most $\ell+1$ it is a {\em proper} $(\le \ell)$-subdivision.
For $\mu\ge 0$ and $r\ge 1$, we denote the $r$-subdivision of $K_{\mu,\mu}$ by $K_{\mu,\mu}^r$. We will frequently use the fact that 
for every graph $H$, there exists $\mu>0$ such that $K_{\mu,\mu}^1$ contains a subdivision of $H$. (To see this, let $n=|V(H)|$, and 
let $\mu=n(n-1)/2$. There is a subgraph of $K_{n,\mu}$ (not induced) that is isomorphic to the 1-subdivision of the complete 
graph $K_n$; and hence there is an {\em induced} subgraph of $K_{\mu,\mu}^1$ isomorphic to the 3-subdivision of $K_n$, which therefore 
contains a $3$-subdivision of $H$.)

So which graphs $H$ have the property that every graph with large chromatic number contains either a large clique or an induced copy of
a subdivision of $H$?
We have seen in \ref{longholes} and \ref{forests} that this is true for cycles and forests.
 Perhaps many more graphs have the same property? 
 For instance, it is known that $K_4$ has this property (this was proved by Scott; see~\cite{leveque}); but it follows from \ref{esperet} below that there are subdivisions of $K_4$ that do not have the property.
Figuring out which graphs do have the property would be a considerable step forward, 
but unfortunately this still seems out of reach.

Here is what seems to be a more tractable question of the same type, solving which would also extend \ref{longholes} and \ref{forests}.
An {\em ideal} of graphs is a class of graphs $\mathcal{C}$, closed under isomorphism and under induced subgraphs (that is, 
if $G\in \mathcal{C}$ and $H$ is isomorphic to 
an induced subgraph of $G$ then $H\in \mathcal{C}$.)
Let us say a graph $H$
is {\em pervasive} in some ideal of graphs $\mathcal{C}$ if for all $\nu,\ell\ge 0$ there exists $c$ 
such that for every graph $G\in \mathcal{C}$ with $\omega(G)\le \nu$ and $\chi(G)>c$,
there is an induced subgraph of $G$ isomorphic to an ($\ge \ell$)-subdivision of $H$. We say $H$ is {\em pervasive} if 
it is pervasive in the ideal of all graphs. Which graphs are pervasive?

If $H'$ is a subdivision of $H$, then $H'$ is pervasive if and only if $H$ is pervasive; and
\ref{longholes} is equivalent to the statement that all cycles are pervasive (and also equivalent
to the assertion that $K_3$ is pervasive). By \ref{forests}, all forests are pervasive; but what else?

There is a beautiful example of Pawlik, Kozik, Krawczyk, Laso\'{n},
Micek, Trotter and Walczak~\cite{sevenpoles}; they found a sequence of graphs $SP_k$ for $k = 1,2,\ldots$, each with
clique number at most two and with chromatic number at least $k$. Essentially the same graphs were constructed in 
a different way by Burling~\cite{burling}, and they are called {\em Burling graphs}.
These graphs are all string graphs (a {\em string graph} is the 
intersection graph of some set of curves in the plane); and consequently for any non-planar graph $H$, 
no proper subdivision
of $H$ appears in any $SP_k$ as an induced subgraph. 
For every pervasive graph $H$, some $(\ge 2)$-subdivision of $H$ must appear
in some  $SP_k$ as an induced subgraph, and this severely restricts the possibilities for which graphs might be pervasive.
This was analyzed in a paper by Chalopin, Esperet, Li and Ossona de Mendez~\cite{chandeliers}, which we discuss next.

\begin{figure}[ht]
\centering

\begin{tikzpicture}[scale=.5,auto=left]
\node [above] (0,0) {pivot};
\tikzstyle{every node}=[inner sep=1.5pt, fill=black,circle,draw]
\node (z) at (0,0) {};
\node (a) at (-6,-2) {};
\node (b) at ( -4,-2) {};
\node (c) at ( -2,-2) {};
\node (d) at ( 0,-2) {};
\node (e) at (2,-2) {};
\node (f) at (4,-2) {};
\node (g) at (6,-2) {};
\node (a1) at (-3.5, -5) {};
\node (d1) at (1, -3) {};
\node (f1) at (5, -4) {};
\node (d2) at (2, -4.5) {};
\node (a2) at (-2, -6) {};
\node (f2) at (3, -6) {};
\node (x) at (0, -7) {};

\foreach \from/\to in {z/a,z/b, z/c,z/d,z/e,z/f,z/g,a/a1,b/a1,c/a1,d/d1,e/d1,f/f1,g/f1,a1/a2,d1/d2,f1/f2,d2/f2,a2/x,f2/x}
\draw [-] (\from) -- (\to);
\end{tikzpicture}

\caption{A chandelier} \label{fig:1}
\end{figure}
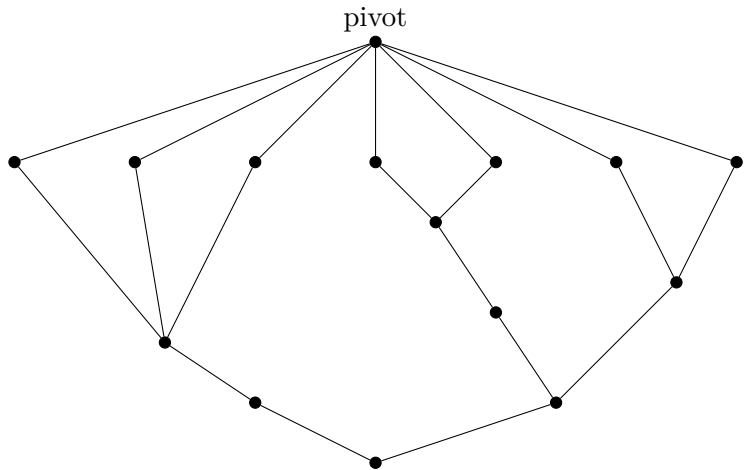

Let $T$ be a tree with $|V(T)|\ge 2$, and let $H$ be obtained from $T$ by adding a new vertex $v$ and 
making $v$ adjacent to every leaf of $T$ (and possibly to some more vertices of $T$); we call $H$ a {\em lantern} with {\em pivot} $v$. 
If $v$ is adjacent only to the leaves of $T$, $H$ is called a {\em chandelier} with {\em pivot} $v$. (In particular, the 1-subdivision of every lantern is a chandelier.) We also count the one- and two-vertex
complete graphs as lanterns and chandeliers, when some vertex is chosen as pivot.
More generally, if we start with a lantern, and repeatedly take a new lantern, and identify its pivot with some vertex 
of what we have already built, what results is called a 
{\em tree of lanterns}, and a {\em tree of chandeliers} is defined similarly. 
If every component of $G$ is a tree of lanterns, $G$ is called a {\em forest of lanterns}, and a {\em forest of chandeliers}
is defined similarly.
It follows from results of Chalopin, Esperet, Li and Ossona de Mendez~\cite{chandeliers} (combine the proof of their
theorem 4.5, their theorem B.4, and the fact that every forest of lanterns is an induced subgraph of some tree of lanterns)
that:
\begin{thm}\label{esperet}
For every graph $H$, there is an $(\ge 2)$-subdivision of $H$ that appears as an induced subgraph in $SP_k$ for some $k$,
if and only if $H$ is a forest of lanterns.\footnote{Incidentally, for a long time we misunderstood the content of the theorem of~\cite{chandeliers},
and thought that \ref{esperet} concerned forests of chandeliers rather than forests of lanterns. In particular we mis-stated
\ref{esperet} in the papers \cite{chibounded, bananatrees}.}

\end{thm}

It follows that every pervasive graph is a forest of lanterns; and perhaps the converse is true, that every forest of
lanterns is pervasive. Whether that is true or not, the goal of this paper is to begin to determine which graphs
are pervasive; and the results we obtain are strong enough that, for pervasiveness in the ideal of string graphs, they tell us the complete answer.
We only have to consider trees of lanterns (since every forest of lanterns is an induced subgraph of a tree of 
lanterns), and they have the convenient property that every subdivision
of a tree of lanterns is another tree of lanterns. Thus, if we could prove that for every tree of lanterns $H$,
every graph with bounded clique number and sufficiently large chromatic number contains a subdivision of $H$
as an induced subgraph, then it would follow that every tree of lanterns is pervasive. We can therefore forget about 
looking for ($\ge \ell$)-subdivisions, and just look for subdivisions. There is also another simplification: every tree
of lanterns has a subdivision that is a tree of chandeliers, and if we can prove that the latter is pervasive, then so is the original tree
of lanterns. So it suffices to prove that every tree of chandeliers is pervasive. The reason for working with
chandeliers instead of lanterns is that nicer things are true for chandeliers, as we shall see.

If $X\subseteq V(G)$, the subgraph of $G$ induced on $X$ is denoted by $G[X]$,
and we often write $\chi(X)$ for $\chi(G[X])$. The {\em distance} between two vertices $u,v$
of $G$ is the length of a shortest path between $u,v$, or $\infty$ if there is no such path.
If $v\in V(G)$ and $\rho\ge 0$ is an integer, $N_G^{\rho}(v)$ (or $N^{\rho}(v)$, when the graph is clear from the context)
 denotes the set of all vertices $u$ with distance
exactly
$\rho$  from $v$, and $N_G^{\rho}[v]$ or $N^{\rho}[v]$ denotes the set of all $v$ with distance at most $\rho$ from $v$.
If $G$ is a nonnull graph  and $\rho\ge 1$,
we define $\chi^{\rho}(G)$ to be the maximum of $\chi(N^{\rho}[v])$ taken over all vertices $v$ of $G$.
(For the graph $G$ with no vertices we define $\chi^{\rho}(G)=0$.)
Let $\mathbb{N}$ denote the set of nonnegative integers, and let $\phi:\mathbb{N}\rightarrow \mathbb{N}$ be a non-decreasing function.
For $\rho\ge 1$, let us say a graph $G$ is {\em $(\rho,\phi)$-controlled} if 
$\chi(H)\le \phi(\chi^{\rho}(H))$ for every induced subgraph $H$ of $G$.
Roughly, this says that in every induced subgraph $H$ of $G$ with
large chromatic number, there is a vertex $v$ such that $N^{\rho}_H[v]$ has large chromatic number. Let us say an ideal of graphs 
$\mathcal{C}$ is {\em $\rho$-controlled} if there is a nondecreasing function $\phi:\mathbb{N}\rightarrow \mathbb{N}$
such that every graph in the ideal is $(\rho,\phi)$-controlled.

Sometimes, it is helpful to know that a statement is true for all $\rho$-controlled ideals, in order to prove that it 
holds for all ideals. For instance, the proof of the main theorem of~\cite{scott} used this approach, as did McGuinness 
in~\cite{mcguinness}, and as we did in~\cite{longholes} and several other papers of this series. 
We hope that the same approach will be helpful for our current
problem of characterizing the pervasive graphs. In this paper we will prove:

\begin{thm}\label{rhopervasive}
For all $\rho\ge 2$, every forest of lanterns is pervasive in every $\rho$-controlled ideal.
\end{thm}
Every $\rho$-controlled ideal is also $(\rho+1)$-controlled, so large values of $\rho$ give more powerful
cases of \ref{rhopervasive}; but we prove \ref{rhopervasive} by induction on $\rho$, and in fact it is the cases
when $\rho$ is small that are most challenging.
The inductive proof of \ref{rhopervasive} is fairly easy for $\rho\ge 4$, slightly more tricky when $\rho=3$, and most difficult 
by far when $\rho=2$.

\begin{figure}[H]
\centering

\begin{tikzpicture}[scale=.8,auto=left, xscale =.5]
\node [above] at (-15.2,3) {plug};
\tikzstyle{every node}=[inner sep=1.5pt, fill=black,circle,draw]
\node (a) at (0,0) {};
\node (b) at (-2,0) {};
\node (c) at (-4,1) {};
\node (d) at (-6,-1) {};
\node (e) at (-8,2) {};
\node (f) at (-10, 1) {};
\node (g) at (-11.8,0) {};
\node (h) at (-14,1) {};
\node (j) at (-16,2) {};
\node (k) at (-16,1) {};
\node (l) at (-16,0) {};
\node (m) at (-16,-1) {};
\tikzstyle{every node}=[inner sep=1.5pt, fill=white,circle,draw]
\node (r) at (-1.2,3) {};
\node (s) at (-3.2,3) {};
\node (t) at (-5.2,3) {};
\node (u) at (-7.2,3) {};
\node (v) at (-9.2,3) {};
\node (w) at (-11.2,3) {};
\node (x) at (-13.2,3) {};
\node (z) at (-15.2,3) {};

\foreach \from/\to in {a/c,a/b,a/d,b/g,c/e,c/f,e/j,e/h,g/l,g/m,h/k}
\draw [very thick] (\from) -- (\to);
\foreach \from/\to in {r/b,r/c,r/d,s/c,s/d,s/g,t/e,t/f,t/g,t/d,u/e,u/f,u/g,v/j,v/h,v/f,v/g,w/j,w/h,w/g,x/j,x/h,x/l,x/m,z/j,z/k,z/l,z/m}
\draw [thin] (\from) -- (\to);

\end{tikzpicture}

\caption{A lamp: each white vertex is adjacent to the left ends of the tree edges below it} \label{fig:0}
\end{figure}
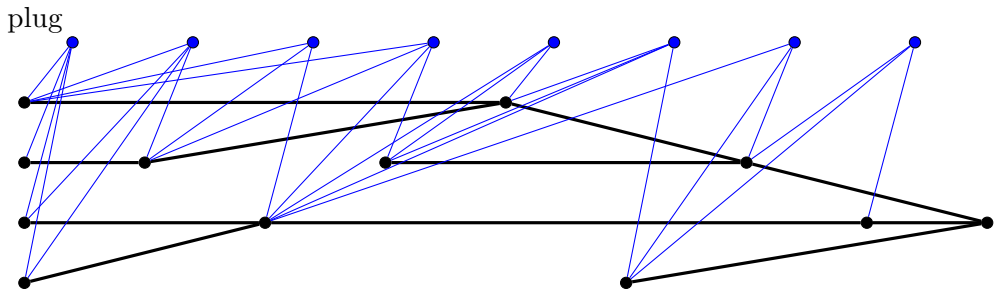

As we saw earlier, to prove \ref{rhopervasive}, it suffices to show that for all $\rho\ge 2$, every tree of chandeliers is pervasive in every $\rho$-controlled ideal.
A ``lamp'' (defined later, see figure 2) is a kind of graph considerably more general than a chandelier, and we will define trees 
of lamps. 
Every chandelier is a lamp, and every chandelier is a lantern, but some lamps are not lanterns (such as the one in 
figure \ref{fig:0}), and some lanterns are not lamps (since lanterns can have triangles and lamps do not, for instance.)
We also think that some trees of chandeliers are also not trees of lamps, because the composition rule for trees of lamps is more restrictive; but for every forest of 
lanterns $H$ there is a tree of lamps that contains a subdivision of $H$ as an induced subgraph.

We will in fact prove something much stronger than \ref{rhopervasive}:
\begin{thm}\label{rhostronger}
For all $\rho\ge 2$, if $\mathcal{C}$ is a $\rho$-controlled ideal of graphs, then either 
\begin{itemize}
\item $\mathcal{C}$ contains every tree of lamps; or
\item $\mathcal{C}$ contains a subdivision of every graph; or
\item for all $\nu\ge 0$, there exists $c$
such that $\chi(G)\le c$ for every graph $G\in \mathcal{C}$ with $\omega(G)\le \nu$.
\end{itemize}
\end{thm}
\noindent{\bf Proof of \ref{rhopervasive}, assuming \ref{rhostronger}\ \ }
Let $\mathcal{C}'$ be a $\rho$-controlled ideal, let $Q'$ be a forest of lanterns, and 
let $\nu,\ell\ge 0$. We must show that there exists $c$ such that  for every graph $G\in \mathcal{C}'$ with $\omega(G)\le \nu$ and 
$\chi(G)>c$,
there is an induced subgraph of $G$ isomorphic to an ($\ge \ell$)-subdivision of $Q'$.
Let $Q$ be a tree of lamps that contains an $(\ge \ell)$-subdivision of $Q'$, 
and let $\mathcal{C}$ be the ideal of all graphs in $\mathcal{C}'$
that contain no subdivision of $Q$. Since $\mathcal{C}$ does not contain every tree of lamps (because it does not contain $Q$),
and $\mathcal{C}$ does not contain a subdivision of every graph (because it does not contain a subdivision of $Q$), it follows from
\ref{rhostronger} that there exists $c$
such that $\chi(G)\le c$ for every graph $G\in \mathcal{C}$ with $\omega(G)\le \nu$. Let $G\in \mathcal{C}'$ with $\omega(G)\le \nu$ and
$\chi(G)>c$. It follows that $G\notin \mathcal{C}$, and so $G$ contains a subdivision of $Q$, and hence contains 
an $(\ge \ell)$-subdivision of $Q'$. This proves that $\mathcal{Q}'$ is pervasive in $\mathcal{C}'$, and so proves \ref{rhopervasive}.~\bbox

Incidentally, the first bullet of \ref{rhostronger} is about trees of lamps; is it also true for trees of lanterns? 
In particular, does \ref{rhostronger} remain true if we replace its first bullet by ``$\mathcal{C}$
contains every lantern''? Let $\mathcal{C}$ be the ideal of all induced subgraphs of Burling graphs; then $\mathcal{C}$ is 
2-controlled (because its members are all string graphs), and the second and third bullets of \ref{rhostronger} are false.
So our question becomes: does every lantern appear as an induced subgraph of some Burling graph? The answer is no.
For instance, lanterns with triangles cannot appear in this way, and nor does 
the lantern consisting of three cycles of length four, with a common edge but otherwise vertex-disjoint. (We thank
Louis Esperet for the latter example.)

To prove the $\rho=2$ case of \ref{rhostronger}, we will show:
\begin{thm}\label{2control}
Let $\nu\ge 0$, let $Q$ be a tree of lamps, and let
$\mu\ge 0$. Let $\mathcal{C}$ be a $2$-controlled ideal of graphs. 
Then there exists $c$ such that every graph $G$ in $\mathcal{C}$ with $\omega(G)\le \nu$ and $\chi(G)>c$
contains one of $K_{\mu,\mu}^1, Q$ as an induced subgraph.
\end{thm}
The general case ($\rho\ge 2$) of the proof of \ref{rhostronger} follows from:
\begin{thm}\label{3control}
Let $\mu\ge 0$, and let
$\rho\ge 2$. Let $\mathcal{C}$ be a $\rho$-controlled ideal of graphs. The ideal of all graphs in $\mathcal{C}$ 
that do not contain
any of $K_{\mu,\mu}^1\ll K_{\mu,\mu}^{\rho+2}$ as an induced subgraph is $2$-controlled.
\end{thm}
We will prove \ref{3control} first, in sections 3--6; and then the sections 7--11 are devoted to proving \ref{2control}.

Why work with lamps rather than chandeliers? For the application
to pervasiveness we could do the whole proof using trees of chandeliers instead of trees of lamps, but there is not much gain; and
\ref{rhostronger} is sufficiently striking that we wanted to prove it for the most general type of graph that we could.

The ideal of all string graphs fits particularly well with \ref{2control}, because:
\begin{itemize}
\item The graph $SP_k$ is a string graph, so only forests of lanterns are pervasive in the ideal of all string graphs.
\item We will prove that the ideal of string graphs is $2$-controlled.
\item Consequently a graph is pervasive in the ideal of all string graphs if and only if it is a forest of lanterns.
\item Since $K_{3,3}^1$ is not a string graph, and hence not an induced subgraph of a string graph, taking $\mu=3$
in \ref{2control} tells us:
if $\nu\ge 0$, and $Q$ is a tree of lamps,
then there exists $c$ such that every string graph $G$ with $\omega(G)\le \nu$ and $\chi(G)>c$
contains $Q$ as an induced subgraph.
\item Consequently we have inadvertently proved the Gy\'arf\'as-Sumner conjecture~\cite{gyarfastree,sumner} for string graphs, 
since every tree is a tree of lamps (and in fact proved much more).
\end{itemize}
We handle string graphs in the final section.

What about ideals that are not $\rho$-controlled? So far, we have not been able to prove that every tree of lanterns is pervasive
in the ideal of all graphs, but two of us prove in~\cite{bananatrees}, using \ref{rhopervasive},
 that all ``banana trees'' are pervasive in 
this ideal (a {\em banana tree} is a multigraph obtained from a tree by adding parallel edges). 


\section{Defining $SP_k$}

Before we go on, let us digress to define $SP_k$. We will not need it in what follows, but 
our work was greatly influenced by
the paper~\cite{chandeliers}, which is based on this construction.

First, here is a composition operation.
We start with a graph $A$, and a stable subset $S$ of $A$. Let $S = \{a_1\ll a_s\}$ say, and for $1\le i\le s$
let $N_i$ be the set of neighbours of $a_i$ in $A$. 

Now take a graph consisting of $s+1$ isomorphic copies of $A\setminus S$, say $A_0\ll A_s$, 
pairwise disjoint and with no edges between them.
For $0\le i,j\le s$, let the isomorphism from $A\setminus S$ to $A_i$ map $N_j$ to $N_{ij}$.
Now add to this $3s^2$ new vertices, namely $x_{ij},y_{ij},z_{ij}$ for all $i,j$ with $1\le i,j\le s$.
Also add edges
so that $x_{ij},y_{ij}$ are both adjacent to every vertex in $N_{0,i}$,
and $x_{ij},z_{ij}$ are both adjacent to every vertex in $N_{ij}$, and $y_{ij}z_{ij}$ an edge, for $1\le i,j\le s$.
Let $G$ be the resulting graph, and let $T$ be the set
$$\{x_{ij},y_{ij}\::1\le i,j\le s\}.$$
We say that $(G,T)$ is obtained by {\em composing $(A,S)$ with itself}. 

To define $SP_k$ let $SP_1$ be the complete graph $K_2$, and let $T_1\subseteq V(SP_1)$ with $|T_1|=1$. Inductively let
$(SP_{k+1},T_{k+1})$ be obtained by composing $(SP_k,T_k)$ with itself. 
It is easy to check that $SP_k$ has no triangles, and for every colouring of $SP_k$ with any number of colours, 
some vertex in $T_k$
has neighbours of $k$ different colours, and in particular $\chi(SP_k)\ge k+1$. 
Moreover, there are graphs $H$ such that no subdivision of $H$ appears as an induced subgraph of any $SP_k$, as discussed in
the previous section. $SP_k$ is the only construction known to the authors with this property. 
Indeed, the following very wild statement might be true as far as we know:
\begin{thm}\label{wildconj}
{\bf Conjecture:} For all $m,i,\nu\ge 0$ there exists $n$ such that if $G$ has $\omega(G)\le \nu$ and $\chi(G)>n$, then
either some $(\ge 1)$-subdivision of $K_{m}$ appears in $G$ as an induced subgraph, or $SP_i$ appears in $G$ as an 
induced subgraph.
\end{thm}
We have little faith in this conjecture; indeed we cannot prove it even for graphs $G$ that are themselves induced subgraphs
of some $SP_k$. We could make it more plausible by weakening it to:
``For all $i,\nu\ge 0$ there exists $n$ such that if $G$ has $\omega(G)\le \nu$ and $\chi(G)>n$, then
some subdivision of $SP_i$ appears in $G$ as an
induced subgraph'', and indeed then we think it might well be true; but first we should disprove the stronger form.

\section{Two routing lemmas}

If $X,Y$ are subsets of the vertex set of a graph $G$, we say
\begin{itemize}
\item $X$ is {\em complete} to $Y$ if $X\cap Y=\emptyset$ and every vertex in $X$
is adjacent to every vertex in $Y$;
\item $X$ is {\em anticomplete} to $Y$ if $X\cap Y=\emptyset$ and every vertex in $X$
is nonadjacent to every vertex in $Y$; and
\item $X$ {\em covers} $Y$ if $X\cap Y=\emptyset$ and every vertex in $Y$ has a neighbour in $X$.
\end{itemize}
(If $X=\{v\}$ we say $v$ is complete to $Y$ instead of $\{v\}$, and so on.)

Throughout the paper, we will be applying various forms of Ramsey's theorem. Here is one that contains all that we need 
(see theorem 5 on page 113 of~\cite{GRS}).
\begin{thm}\label{ramsey}
For all integers $k, n,\alpha, \beta\ge 0 $ there exists $R(k,n,\alpha,\beta)\ge n$ with the following property. Let $A,B$ be disjoint sets, both of cardinality
at least $R(k,n,\alpha,\beta)$. Let $E$ be the set of all sets $X\subseteq A\cup B$ with $|X\cap A| = \alpha$ and $|X\cap B|=\beta$.
If we partition $E$ into $k$ subsets, then there
exist $A'\subseteq A$ and $B'\subseteq B$ with $|A'| = |B'| = n$ such that all the sets $X\in E$ with $X\subseteq A'\cup B'$
belong to the same subset.
\end{thm}

Before we begin the main proofs, we prove two lemmas which will be applied later. We are trying to prove that certain graphs $G$
with bounded clique number
contain a subdivision of some fixed graph $H$ as an induced subgraph. This is true if $G$ has an
induced subgraph which is a proper subdivision of $K_{\mu,\mu}$ for appropriate $\mu$; and so we might as well confine ourselves
to graphs $G$ that do not contain (as an induced subgraph) any proper subdivision of $K_{\mu,\mu}$, for some fixed $\mu$.
This is a little more than we actually need; we only need to exclude subdivisions in which each edge is subdivided  
a small number of times.
For integers $\lambda\ge 2$ and $\mu,\nu\ge 0$, let us say that $G$ is {\em $(\lambda,\mu,\nu)$-restricted} if
$\omega(G)\le \nu$, and no induced subgraph of $G$ is a proper $(\le \lambda)$-subdivision of $K_{\mu,\mu}$.

Let $G,H$ be graphs. An {\em impression} of $H$ in $G$
is a map $\eta$ with domain $V(H)\cup E(H)$, such that:
\begin{itemize}
\item $\eta(v)\in V(G)$ for each $v\in V(H)$;
\item for all distinct $u,v\in V(H)$, $\eta(u)\ne \eta(v)$ and  $\eta(u), \eta(v)$ are nonadjacent in $G$;
\item for every edge $e=uv$ of $H$, $\eta(e)$ is a path of $G$ with ends $\eta(u), \eta(v)$;
\item if $e,f\in E(H)$ have no common end then $V(\eta(e))$ is anticomplete to $V(\eta(f))$.
\end{itemize}
The {\em order} of an impression $\eta$ is the maximum length of the paths $\eta(e)\;(e\in E(H))$.
Our first lemma is:

\begin{thm}\label{impression}
For all $\lambda\ge 1$ and $\mu,\nu\ge 0$, there exists $n$ such that if 
$\omega(G)\le \nu$, and $G$ does not contain
any of $K_{\mu,\mu}^1\ll K_{\mu,\mu}^{\lambda}$ as an induced subgraph (and in particular if $G$ is $(\lambda,\mu,\nu)$-restricted)
then there is 
no impression of $K_{n,n}$ in $G$ of order at most $\lambda+1$.
\end{thm}
\Proof We proceed by induction on $\lambda$. If $\lambda>1$ choose $m_4$ so that the theorem is satisfied with $\lambda$ replaced
by $\lambda-1$ and $n$ by $m_4$, and if $\lambda=1$ let $m_4=0$.
Let 
\begin{eqnarray*}
m_3&=&\max(m_4+1,\mu,\nu+2)\\
m_2&=&R(3^{\lambda^2}, m_3,2,1)\\
 m_1&=&R(3^{\lambda^2}, m_2,1,2)\\
n&=&R(\lambda, m_1,1,1).
\end{eqnarray*}
We claim that $m$ satisfies the theorem. For let $H=K_{n,n}$, and suppose that $\eta$ is an
impression of $H$ in $G$ of order at most $\lambda+1$. 
\\
\\
(1) {\em $\{\eta(v):v\in V(H)\}$ is a stable set of $G$, and if $e\in E(H)$ and $v\in V(H)$ is not incident with~$e$, then 
$\eta(v)$ does not belong
to $\eta(e)$, and has no neighbours in $V(\eta(e))$.}
\\
\\
The first is immediate from the definition of impression. For the second,
if $e\in E(H)$ and $v\in V(H)$ not incident with $e$, then 
there is an edge $f$ of $H$ incident with $v$ and with no common end with $e$, and since 
$V(\eta(e))$ is anticomplete to $V(\eta(f))$, it follows in particular that $\eta(v)$ does not belong
to $\eta(e)$, and has no neighbours in $V(\eta(e))$. This proves (1).

\bigskip

Also we might as well assume that each path $\eta(e)$ is an induced path in $G$.
Let $(A,B)$ be a bipartition of $H=K_{n,n}$. 
There are only $\lambda$ possibilities for the length of each path $\eta(e)\;(e\in E(H))$; and so by \ref{ramsey}, there
exist $A_1\subseteq A$ and $B_1\subseteq B$  with $|A_1| = |B_1| = m_1$ such that the paths $\eta(ab)$ all have the same length,
for all $a\in A_1$ and $b\in B_1$. Let this common length be $\ell$; thus $2\le \ell \le \lambda+1$.
Let us number the vertices of each path $\eta(ab)\;(a\in A_1, b\in B_1)$ as $p_{ab}^0, p_{ab}^1\ll p_{ab}^{\ell}$ in order,
where $p_{ab}^0 = \eta(a)$ and $p_{ab}^{\ell} = \eta(b)$. 

Take an ordering of $B_1$, denoted by $<$. For each $a\in A_1$ and 
all $b,b'\in B_1$ with $b<b'$, let us say the {\em first pattern} of $(a,b,b')$ is the set of all pairs $(i,j)$ with
$1\le i,j\le \ell-1$ such that $p_{ab}^i = p_{ab'}^j$; and the {\em second pattern} of $(a,b,b')$ is the set of all pairs $(i,j)$ with
$1\le i,j\le \ell-1$ such that $p_{ab}^i,p_{ab'}^j$ are distinct and adjacent in $G$. There are only $3^{\lambda^2}$ possibilities
for the first and second patterns; so by \ref{ramsey} there exist $A_2\in A_1$ and $B_2\subseteq B_1$ with $|A_2| = |B_2| = m_2$,
such that all the triples $(a,b,b')$ (for $a\in A_2$ and $b,b'\in B_2$ with $b<b'$) have the same first patterns and they all
have the same second patterns. Let these patterns be $\Pi_1,\Pi_2$ say. 

Similarly, by exchanging $A,B$, choosing an ordering $<$ of $A_2$ and repeating the argument, 
we deduce that there exist $A_3\subseteq A_2$ and $B_3\subseteq B_2$
with $|A_3|=|B_3| = m_3$, and sets $\Pi_3,\Pi_4\subseteq \{1\ll \ell-1\}^2$ such that for all $a,a'\in A_3$ with $a<a'$ and $b\in B_3$,
$p_{ab}^i=p_{a'b}^j$ if and only if $(i,j)\in \Pi_3$, and $p_{ab}^i,p_{a'b}^j$ are different and adjacent if and only if $(i,j)\in \Pi_4$.
\\
\\
(2) {\em $\Pi_1, \Pi_2=\emptyset$.}
\\
\\
For suppose that there exists $(i,j)\in \Pi_1\cup \Pi_2$. By reversing the order on $B$ if necessary, we may assume that $i\le j$.
Choose $b_0\in B_3$, minimal under the ordering of $B_1$.
For each $a\in A_3$ and $b\in B_3\setminus \{b_0\}$, let 
$$Q(ab)=\{ p_{ab}^j, p_{ab}^{j+1},\ll p_{ab}^{\ell}\}.$$
Since $(i,j)\in \Pi_1\cup \Pi_2$, it follows that for each $a\in A_3$ and $b\in B_3\setminus \{b_0\}$,
there is a path $P_{ab}$ of $G$ with ends $p_{ab_0}^i, b$ and with vertex set a subset of $\{p_{ab_0}^i\}\cup Q(ab)$.
For each $b\in B_3\setminus \{b_0\}$ let $\eta'(b)=\eta(b)$; for each $a\in A_3$, let $\eta'(a)=p_{ab_0}^i$;
and for every edge $ab$ of $H=K_{n,n}$ with $a\in A_3$ and $b\in B_3\setminus \{b_0\}$, let $\eta'(ab)=P_{ab}$.
We claim that $\eta'$ is an impression of $K_{m_3,m_3-1}$ in $G$.
To see this, note first that the vertices $\eta'(a)\;(a\in A_3)$ are all distinct; for choose $b\in B_3\setminus \{b_0\}$,
and let $a,a'\in A_3$ be distinct. Then $p_{ab_0}^i$ is equal or adjacent to $p_{ab}^j$, but $p_{a'b_0}^i$ is different from and 
nonadjacent to $p_{ab}^j$
since $V(\eta(a'b_0)),V(\eta(ab))$ are anticomplete, from the definition of an impression. Consequently
$p_{ab_0}^i$ is different from $p_{a'b_0}^i$.
If $(i,i)\in \Pi_4$, then all the vertices $p_{ab_0}^i\;(a\in A_3)$ are pairwise adjacent, contradicting that $\omega(G)\le\nu$;
so $(i,i)\notin \Pi_4$, and the vertices
$\eta'(a)\;(a\in A_3)$ are pairwise nonadjacent. Also for each $a\in A_3$ and $b\in B_3\setminus \{b_0\}$,
$\eta'(a)$ is different from and nonadjacent to $\eta'(b)$ by (1). Thus the first three conditions for an impression are
satisfied. For the final condition, we must check that if $a,a'\in A_3$ are distinct and $b,b'\in B_3\setminus \{b_0\}$ are distinct,
then $V(P_{ab})$ is anticomplete to $V(P_{a'b'})$. We recall that $V(P_{ab})\subseteq \{p_{ab_0}^i\}\cup Q(ab)$, where $Q(ab)$ is 
a subset of the vertex set of $\eta(ab)$, and $V(P_{a'b'})\subseteq \{p_{a'b_0}^i\}\cup Q(a'b')$.
We have seen that $p_{ab_0}^i, p_{a'b_0}^i$ are distinct and nonadjacent, so, exchanging $a,a'$ and $b,b'$ if necessary,
it suffices to show that $V(P_{ab})$ is anticomplete to $Q(a'b')$.
But $V(P_{ab})$ is a subset of $V(\eta(ab_0))\cup V(\eta(ab))$, and both the latter sets are 
anticomplete to $V(\eta(a'b'))\supseteq Q(a'b')$. This proves that $\eta'$ is an impression as claimed. 

Since $m_3-1\ge m_4$, the inductive hypothesis on $\lambda$ implies that the order of $\eta'$ is at least $\lambda+1$. But its order
is at most $\ell-j+1$ if  $(i,j)\in \Pi_2$, and at most $\ell-j$ if $(i,j)\in \Pi_1$. Since $\ell\le \lambda+1$ and $j\ge 1$,
we deduce that $j=1$, and $\ell=\lambda+1$; and so $i=1$, since $i\le j$, and $(1,1)\in \Pi_2$. Choose $a\in A_3$; then all the vertices
$p_{ab}^1\;(b\in B_3\setminus \{b_0\})$ are distinct and pairwise adjacent, contradicting that $\omega(G)\le\nu$. This proves (2).

\bigskip

Similarly $\Pi_3,\Pi_4=\emptyset$. But then $G$ contains an $\ell$-subdivision of $K_{m_3,m_3}$,
contrary to the hypothesis. This proves \ref{impression}.~\bbox

The second lemma is:

\begin{thm}\label{stephanproof}
For all $\mu,\nu\ge 0$, there exists $m$ with the following property. Let $G$ be $(1,\mu,\nu)$-restricted,
and let $X\subseteq V(G)$ with $|X|\ge m$. Then there exist distinct nonadjacent $x,x'\in X$ such that every vertex of $G$ adjacent to both $x,x'$ 
has at least one more neighbour in $X$.
\end{thm}
\Proof Choose $m_4$ so that \ref{impression} holds with $n$ replaced by $m_4$.
Let 
\begin{eqnarray*}
m_3&=&\max(m_4,\nu+1);\\
m_2 &=& R(4,m_3, 2,2);\\
m_1 &=& 2m_2;\\
m&=&R(2,m_1,2,0).
\end{eqnarray*}
We claim that $m$ satisfies the theorem. For suppose that $G, X$ are as in the theorem, and for all distinct nonadjacent
$x,x'\in X$ there exists $w(x,x')$ adjacent to both $x,x'$ and nonadjacent to all other vertices in $X$. Since
$\omega(G)\le \nu< m_1$, there is a stable subset $X_1$ of $X$ with $|X_1| = m_1$, by \ref{ramsey}. 
It follows that all the vertices $w(x,x')\;(x,x'\in M_1, x\ne x')$ are distinct from one another
and distinct from the vertices in $M_1$. Choose two disjoint subsets
$A_2,B_2$ of $X_1$, both of cardinality $m_2$. Take an ordering of $A_2$ and of $B_2$, both denoted by $<$.
Let $E$ be the set of all quadruples $(a,a',b,b')$ such that $a,a'\in A$, $a<a'$, and $b,b'\in B$ and $b<b'$.
For all $(a,a',b,b')\in E$,
we say the {\em first pattern} of $(a,a',b,b')$
is $1$ or $0$ depending whether $w(a,b), w(a',b')$ are adjacent or not; and the {\em second pattern} is $1$ or $0$
depending whether $w(a,b'), w(a', b)$ are adjacent or not. There are four possible choices of first and second pattern;
so by 
\ref{ramsey} there exist $A_3\subseteq A_2$ and $B_3\subseteq B_2$ with $|A_3|=|B_3| = m_3$, 
such that, if $E_3$ denotes the set of $(a,a',b,b')\in E$ with $a,a'\in A_3$ and $b,b'\in B_3$, then
\begin{itemize}
\item either $w(a,b), w(a',b')$ are adjacent for all $(a,a',b,b')\in E_3$, or
$w(a,b), w(a',b')$ are nonadjacent for all $(a,a',b,b')\in E_3$; and
\item either $w(a,b'), w(a',b)$ are adjacent for all $(a,a',b,b')\in E_3$, or
$w(a,b'), w(a',b)$ are nonadjacent for all $(a,a',b,b')\in E_3$.
\end{itemize}
Suppose that $w(a,b), w(a',b')$ are adjacent for all $(a,a',b,b')\in E_3$.
Choose 
\begin{eqnarray*}
a_1<a_2<\cdots<a_{\nu+1}\in A_3\\ 
b_1<b_2<\cdots<b_{\nu+1}\in B_3
\end{eqnarray*}
(this is possible since $m_3\ge \nu+1$); then the
vertices $w(a_1,b_1), w(a_2,b_2)\ll w(a_{\nu+1},b_{\nu+1})$ are pairwise adjacent, contradicting that $\omega(G)\le\nu$.
So the nonadjacency alternative holds in the first bullet above, and similarly nonadjacency holds in the second bullet.
Let $(A',B')$ be a bipartition of $K_{m_3,m_3}$, and choose $\eta$ mapping $A'$ onto $A$ and $B'$ onto $B$; and for
all $a'\in A'$ and $b'\in B'$, let $\eta(a'b')$ be the path of $G$ with vertex set $\{a,w(a,b),b\}$ where $a=\eta(a')$ and $b=\eta(b')$.
Then $\eta$ is an impression of $K_{m_3,m_3}$ in $G$, of order $2$, and the result follows from \ref{impression}. This
proves \ref{stephanproof}.~\bbox

\section{Reducing control}

A {\em levelling} in a graph $G$ is a sequence of pairwise disjoint subsets $(L_0, L_1\ll L_k)$ of $V(G)$
such that
\begin{itemize}
\item $|L_0|=1$;
\item for $1\le i\le k$, $L_{i-1}$ covers $L_i$; and
\item for $0\le i<j\le k$, if $j>i+1$ then $L_i$ is anticomplete to $L_j$.
\end{itemize}
If $\mathcal{L}= (L_0, L_1\ll L_k)$ is a levelling, 
$L_k$ is called the {\em base} of $\mathcal{L}$, and the vertex in $L_0$ is the {\em apex} of $\mathcal{L}$,
and $L_0\cup\cdots\cup L_k$ is the {\em union} of $\mathcal{L}$, denoted by $V(\mathcal{L})$.
If $\mathcal{L}=(L_0,L_1\ll L_k)$ and $\mathcal{L}'=(L_0',L_1'\ll L_k')$ are levellings, we say that 
$\mathcal{L}'$ is {\em contained in} $\mathcal{L}$
if $L_i'\subseteq L_i$ for $0\le i\le k$. For instance, one can obtain a levelling (in a connected graph) 
by classifying all vertices by their distance
from some fixed vertex.

Let $\mathcal{L}=(L_0,L_1\ll L_{\rho-1})$ be a levelling in $G$ with $\rho\ge 2$, and let $C\subseteq V(G)\setminus V(\mathcal{L})$.
We say that $\mathcal{L}$ is a {\em $\rho$-cover for $C$} if $L_{\rho-1}$ covers $C$, and $L_0\ll L_{\rho-2}$ are anticomplete to $C$, 
that is, if $(L_1\ll L_{\rho-1}, C)$ is a levelling.
Let $\mathcal{L} = (L_0\ll L_{\rho-1})$ be a $\rho$-cover for $C$, with apex $x$ say. If $z\in C$, then
$z$ has a neighbour in $L_{\rho-1}$, and that vertex has a neighbour in $L_{\rho-2}$, and so on; and hence there is a path
between $z$ and $x$ of length $\rho$, with exactly one vertex in each of $L_0\ll L_{\rho-1}$. Moreover, this path is induced; 
we call such 
a path an {\em $\mathcal{L}$-radius} for $z$.

If we have a $\rho$-controlled ideal that is not $(\rho-1)$-controlled, there are graphs $G$ in the ideal with $\chi^{\rho-1}(G)$
bounded and $\chi^{\rho}(G)$ arbitrarily large. Choose such a graph $G$, with $\chi^{\rho}(G)$ very large; then there is a vertex
$z_1$ with $\chi(N^\rho[z_1])$ very large (not quite so large). For $0\le j\le \rho$, let $L_{1,j}$ be the set of vertices 
with distance $j$ from $z_1$. Since $\chi^{\rho-1}(G)$ is bounded, it follows that
$\chi(N^{\rho}(z_1))=\chi(L_{1,\rho})$ is very large. The subgraph $G_2$ induced on $L_{1,\rho}$ belongs to the same $\rho$-controlled ideal,
and so there is a vertex $z_2$ in it with $\chi(N^{\rho}_{G_2}[z_2])$; let $L_{2,j}$ be the set of vertices in $G_2$ with distance 
 $j$ in $G_2$ from $z_2$, and then as before $\chi(L_{2,\rho})$ is very large. By continuing this process we obtain a sequence of $\rho$-covers,
and that motivates the following definition.

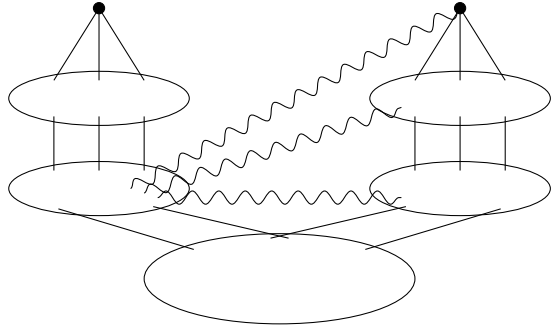
\begin{figure}[H]
\centering
\tikzset{snake it/.style={decorate, decoration=snake}}
\begin{tikzpicture}[scale=.6,auto=left]
\tikzstyle{every node}=[inner sep=1.5pt, fill=black,circle,draw]
\node (x1) at (-4,0) {};
\node (x2) at (4,0) {};
\draw (-4,-2) ellipse (2cm and .6cm);
\draw (-4,-4) ellipse (2cm and .6cm);
\draw (4,-2) ellipse (2cm and .6cm);
\draw (4,-4) ellipse (2cm and .6cm);
\draw (0,-6) ellipse (3cm and 1cm);
\draw (-4,0) -- (-4,-1.6);
\draw (-4,0) -- (-3,-1.6);
\draw (-4,0) -- (-5,-1.6);
\draw (4,0) -- (4,-1.6);
\draw (4,0) -- (3,-1.6);
\draw (4,0) -- (5,-1.6);
\draw (-4,-2.4) -- (-4, -3.6);
\draw (-3,-2.4) -- (-3, -3.6);
\draw (-5,-2.4) -- (-5, -3.6);
\draw (4,-2.4) -- (4, -3.6);
\draw (3,-2.4) -- (3, -3.6);
\draw (5,-2.4) -- (5, -3.6);
\draw (-2.8, -4.4) -- (0.2,-5.1);
\draw (-4.9, -4.45) -- (-1.9,-5.35);
\draw (2.8, -4.4) -- (-0.2,-5.1);
\draw (4.9, -4.45) -- (1.9,-5.35);
\path[draw=black, snake it] (-2.7, -4.2) -- (2.7, -4.2);
\path[draw=black, snake it] (-3.0, -4.1) -- (2.7, -2.2);
\path[draw=black, snake it] (-3.3, -4.0) -- (4,0);

\end{tikzpicture}

\caption{A 3-multicover of length two (wiggly lines indicate possible edges)} \label{fig:2}
\end{figure}

For $C\subseteq V(G)$, a {\em $\rho$-multicover for $C$} in $G$
is a family $\mathcal{M}=(\mathcal{L}_i:i\in I)$, where $I$ is a set of integers,
such that
\begin{itemize}
\item for $1\le i\le m$, $\mathcal{L}_i$ is a $\rho$-cover for $C$;
\item for $1\le i<j\le m$, $V(\mathcal{L}_i)$ is disjoint from $V(\mathcal{L}_{j})$;
\item for all $i,j\in I$ with $i<j$, every vertex in $V(\mathcal{L}_{i})$ with a neighbour in $V(\mathcal{L}_{j})$
belongs to the base of $\mathcal{L}_{i}$.
\end{itemize}
We denote the union of the sets $V(\mathcal{L}_i)\;(i\in I)$
by $V(\mathcal{M})$. We call $|I|$ the {\em length} of the multicover, and $I$ is its {\em index set}.
The next two section are devoted to proving the following:
\begin{thm}\label{usetick2}
For all $\rho\ge 3$ and $\mu, \nu,\tau\ge 0$ there exist $m,c\ge 0$ with the following property.
Let $G$ be a $(\rho+2,\mu,\nu)$-restricted graph such that $\chi^{\rho-1}(G)\le \tau$.
If $C\subseteq V(G)$ with $\chi(C)>c$, then there is no $\rho$-multicover of $C$ in $G$ with length $m$.
\end{thm}

But first, let us assume the truth of \ref{usetick2}, and apply it
to prove a
result of great importance (for us), the following.
\begin{thm}\label{reducecontrol}
Let $\mu,\nu\ge 0$ and $\rho\ge 2$. Every $\rho$-controlled ideal of $(\rho+2,\mu,\nu)$-restricted graphs is $2$-controlled.
\end{thm}
\noindent{\bf Proof (assuming \ref{usetick2}).\ \ }The result is trivial for $\rho=2$, and we proceed by induction on $\rho$. 
Let $\rho\ge 3$, and let $\mathcal{C}$
be a $\rho$-controlled ideal of $(\rho+2,\mu,\nu)$-restricted graphs. Let $\phi$ be nondecreasing such that every graph in $\mathcal{C}$
is $(\rho,\phi)$-controlled. 

Let $\tau\ge 0$, and let $\mathcal{D}$ be the set of all graphs $H\in \mathcal{C}$ with $\chi^{\rho-1}(H)\le \tau$.
Let $m,c$ satisfy \ref{usetick2}.
Define $c_0 = c$, and inductively $c_t = \phi(c_{t-1} + \tau)$ for $t>0$.
We claim:
\\
\\
(1) {\em For $0\le t\le m$, if $H\in \mathcal{D}$ with $\chi(H)>c_t$ then 
there is a $\rho$-multicover in $H$ with length $t$ of some set $C\subseteq V(H)$ where $\chi(C)>c$.}
\\
\\
The claim is trivial if $t=0$, and we proceed by induction on $t$.
Let $H\in \mathcal{D}$ with $\chi(H)>c_t=\phi(c_{t-1} + \tau)$; then since $H$ is $(\rho,\phi)$-controlled, it follows that
$\chi(H)\le \phi(\chi^{\rho}(H))$, and so 
$\chi^{\rho}(H)>c_{t-1}+\tau$. Choose $x\in V(H)$ so that $\chi(N^{\rho}[x])>c_{t-1}+\tau$. Since 
$\chi(N^{\rho-1}[x])\le \tau$, it follows that $\chi(N^{\rho}(x))>c_{t-1}$. For each $i\ge 0$, let $L_i$
be the set of vertices in $H$ with distance exactly $i$ from $x$, and let $J = H[L_{\rho}]$.
Since $\chi(J)>c_{t-1}$, from the inductive hypothesis 
there is a $\rho$-multicover in $J$ with length $t-1$ of some set $C$ where $\chi(C)>c$, say
$(\mathcal{L}_i: 2\le i\le t)$. Define $\mathcal{L}_1=(L_0, L_1\ll L_{\rho-1})$; then
$(\mathcal{L}_i: 1\le i\le t)$ satisfies (1). (Note that every edge between $V(\mathcal{L}_1)$ and $V(\mathcal{L}_i)$ for $i>1$
is also between $V(\mathcal{L}_1)$ and $L_{\rho}$, and therefore has an end in $L_{\rho-1}$.) This proves (1).

\bigskip

From (1) and \ref{usetick2}, it follows that every member of $\mathcal{D}$ has chromatic number at most $c_m$.
At the start of the proof we made an arbitrary choice of $\tau$, and all the subsequent variables in (1)
(such as $\mathcal{D}, m$ and the sequence
$c_0, c_1,\ldots$) depend on $\tau$. In particular, $c_m$ is a function of $\tau$, say $\phi'(\tau)$.
Thus, if $H\in \mathcal{C}$, then $\chi(H)\le \phi'(\chi^{\rho-1}(H))$. 

We may assume that $\phi'$ is nondecreasing; and so every graph in $\mathcal{C}$ is $(\rho-1,\phi')$-controlled, and so
$\mathcal{C}$ is $(\rho-1)$-controlled, and hence $2$-controlled, from the inductive hypothesis. This proves
\ref{reducecontrol}.~\bbox

Next we will deduce \ref{3control}, but before that, here is a useful lemma.
\begin{thm}\label{fixclique}
Let $\rho\ge 2$, and let $\mathcal{C}$ be an ideal of graphs, such that for all $\nu\ge 0$,
the ideal $\mathcal{C}_{\nu}$ of graphs $G\in \mathcal{C}$
with $\omega(G)\le \nu$ is $\rho$-controlled. Then $\mathcal{C}$ is $\rho$-controlled.
\end{thm}
\Proof For each $\nu\ge 0$, let $\phi_{\nu}$ be a function such that each graph $G$ in
$\mathcal{C}_{\nu}$ is $(\rho, \phi_{\nu})$-controlled. For $c\ge 0$, let $\psi(c) = \max_{\nu\le c}\phi_{\nu}(c)$.
We claim that $\mathcal{C}$ is $(\rho,\psi)$-controlled. For let $G\in \mathcal{C}$, and let $H$ be an induced subgraph
of $G$ such that $\chi(H)>\psi(c)$, for some $c$. Let $\nu=\omega(H)$. If $\nu>c$, choose a clique $X$ of $H$
with $|X|>c$, and choose $v\in X$; then $X$ belongs to $N^{\rho}_H[v]$, and so $\chi^{\rho}(H)\ge |X|>c$ as required.
Thus we may assume that $\nu\le c$, and so $\chi(H)>\psi(c)\ge \phi_{\nu}(c)$. Since $G$ is $(\rho, \phi_{\nu})$-controlled,
it follows that $\chi^{\rho}(H)>c$ as required. This proves \ref{fixclique}.~\bbox

Now we prove \ref{3control}, which we restate.

\begin{thm}\label{3controlagain}
Let $\mu\ge 0$ and
$\rho\ge 2$, and let $\mathcal{C}$ be a $\rho$-controlled ideal of graphs. The ideal of all graphs in $\mathcal{C}$
that do not contain
any of $K_{\mu,\mu}^1\ll K_{\mu,\mu}^{\rho+2}$ as an induced subgraph is $2$-controlled.
\end{thm}
\noindent{\bf Proof (assuming \ref{usetick2}).\ \ }Let $\mathcal{D}$ be the ideal of all graphs in $\mathcal{C}$ that do not contain
any of $K_{\mu,\mu}^1\ll K_{\mu,\mu}^{\rho+2}$ as an induced subgraph. 
Let $\nu\ge 0$, and let $\mathcal{D}_{\nu}$ be the ideal of all graphs $G\in \mathcal{D}$ with $\omega(G)\le \nu$.
From \ref{impression}, there exists $n\ge 0$ such that no $G\in \mathcal{D}_{\nu}$ contains an impression of $K_{n,n}$ 
of order at most $\lambda+1$ as an induced subgraph; and consequently 
every graph in $\mathcal{D}_{\nu}$ is 
$(\rho+2,n,\nu)$-restricted. Therefore $\mathcal{D}_{\nu}$ is $2$-controlled by \ref{reducecontrol}, and from
\ref{fixclique} it follows that $\mathcal{D}$ is 2-controlled. This proves \ref{3controlagain}.~\bbox

\section{Extracting ticks from $\rho$-multicovers}

In this section and the next we prove \ref{usetick2}.
Let $\mathcal{M}=(\mathcal{L}_i:i\in I)$ and $\mathcal{M}'=(\mathcal{L}_i':i\in I')$ be $\rho$-multicovers in $G$
for $C$ and for $C'$, respectively, where $C'\subseteq C$.
If $I'\subseteq I$,
and $\mathcal{L}_i'$ is contained in $\mathcal{L}_i$ for each $i\in I'$,
we say that $\mathcal{M}'$ is {\em contained in} $\mathcal{M}$.

Let $\mathcal{M}=(\mathcal{L}_i:i\in I)$ be a $\rho$-multicover for $C$ in $G$.
Let $z\in V(G)\setminus (V(\mathcal{M})\cup C)$, and for each $i\in I$ let $S_i$ be an induced path of $G$ 
between $z$ and the apex $x_i$ say of $\mathcal{L}_i$, such that
\begin{itemize}
\item $z$ has no neighbours in $V(\mathcal{M})\cup C$;
\item for each $i\in I$, $V(S_i)\cap (V(\mathcal{M})\cup C)=\{x_i\}$; and
\item for each $i\in I$, every vertex in $V(\mathcal{M})\cup C$ with a neighbour in $V(S_i)$ belongs to 
$V(\mathcal{L}_i)$.
\end{itemize}
(We do not require the paths $S_i$ to be pairwise internally disjoint; they may intersect one another arbitrarily.)
We say that the family $(S_i:i\in I)$ is a {\em tick} of $G$ on $(\mathcal{M},C)$, and $z$ is its {\em head},
and its {\em order} is the maximum length of the paths $S_i$ for $i\in I$.
We will prove the following. 

\begin{thm}\label{gettick}
For all $\rho\ge 3$ and $\mu,\nu,\tau, m',c'\ge 0$ there exist $m,c\ge 0$ with the following property.
Let $G$ be a  $(1,\mu,\nu)$-restricted graph such that $\chi^{\rho-1}(G)\le \tau$. 
Let $C\subseteq V(G)$ with $\chi(C)>c$, and let
$\mathcal{M}=(\mathcal{L}_i:i\in I)$ be a $\rho$-multicover for $C$ with length $m$. Then 
there exist $C'\subseteq C$ with $\chi(C')> c'$, and a $\rho$-multicover
$\mathcal{M}'$ for $C'$ contained in $\mathcal{M}$ with length $m'$, indexed by $I'\subseteq I$,
and a tick $(S_i:i\in I')$ 
on $(\mathcal{M}',C')$ of order at most $\rho+3$, such that for each $i\in I'$, every vertex of $S_i$ belongs either to
$V(\mathcal{L}_i)$, or to $C$, or to $V(\mathcal{L}_{k})$ for some $k\in I\setminus I'$.
\end{thm}

Before we prove \ref{gettick}, let us see that it implies \ref{usetick2}, which we restate:

\begin{thm}\label{usetick2again}
For all $\rho\ge 3$ and $\mu, \nu,\tau\ge 0$ there exist $m,c\ge 0$ with the following property.
Let $G$ be a $(\rho+2,\mu,\nu)$-restricted graph such that $\chi^{\rho-1}(G)\le \tau$.
If $C\subseteq V(G)$ with $\chi(C)>c$, then there is no $\rho$-multicover of $C$ in $G$ with length $m$.
\end{thm}
\noindent{\bf Proof, assuming \ref{gettick}.\ \ }First, here is a sketch.
By starting with a $\rho$-multicover $\mathcal{M}$ with large enough length, for a set $C$ with chromatic number
large enough, and applying \ref{gettick} repeatedly, we obtain 
a sequence of multicovers, each contained in its predecessor, of successively smaller (but still large) lengths, and
a sequence of ticks all on the last multicover of the sequence $\mathcal{M}'$ say. The ticks are vertex-disjoint except 
for their vertices
in $V(\mathcal{M}')$. There may be edges between them, but if say $(S_i:i\in I)$ and $(T_i:i\in I)$ are two of these ticks, and some
vertex in $S_i$ is adjacent to some vertex in $T_j$, then $i=j$. Consequently we have obtained an impression of $K_{n,n}$ of order
at most $\rho+3$, with $n$ large, which is impossible if $G$ is $(\rho+2,\mu,\nu)$-restricted.

Now let us say it precisely. By \ref{impression}, there exists an integer $n\ge 0$ such that if $G$ is
$(\rho+2,\mu,\nu)$-restricted then there is no impression of $K_{n,n}$ in $G$
of order at most $\rho+3$. Define $m_n=n$ and $c_n=0$; and for $j=n-1,n-2\ll 0$ choose $m_j,c_j$
so that \ref{gettick} holds with $m',c', m,c$ replaced by $m_{j+1}, c_{j+1}, m_j, c_j$ respectively.

Let $m=m_0$ and $c=c_0$; we claim that $m,c$ satisfy the theorem. For let $G$ be $(\rho+2,\mu,\nu)$-restricted 
with $\chi^{\rho-1}(G)\le \tau$,
let $C_0\subseteq V(G)$ with $\chi(C_0)>c_0$, and suppose that $\mathcal{M}_0=(\mathcal{L}_{i0}:i\in I_0)$ is 
a $\rho$-multicover for $C$ with length $m_0$, indexed by $I_0$. 
Inductively, for $1\le j\le n$, we define $C_{j}$, $\mathcal{M}_{j}$, $I_j$ and $\mathcal{T}_j$
as follows. Since $G$ is $(\rho+2,\mu,\nu)$-restricted and hence $(1,\mu,\nu)$-restricted, and
$\mathcal{M}_{j-1}$ is a $\rho$-multicover for $C_{j-1}$ with length $m_{j-1}$, and 
$\chi(C_{j-1})>c_{j-1}$, we can apply \ref{gettick}. We deduce that 
there exist $C_j\subseteq C_{j-1}$ with $\chi(C_j)> c_j$, and a $\rho$-multicover
$\mathcal{M}_j=(\mathcal{L}_{ij}:i\in I_j)$ for $C_j$ contained in $\mathcal{M}_j$ with length $m_j$, and a tick
$\mathcal{T}_j = (S_{ij}:i\in I_j)$ on $(\mathcal{M}_j,C_j)$ of order at most $\rho+3$, such that for each $i\in I_j$,
every vertex of $S_i$ belongs either to
$V(\mathcal{L}_{i,j-1})$, or to $C_{j-1}$, or to $V(\mathcal{L}_{k,j-1})$ for some $k\in I_{j-1}\setminus I_j$.

For $1\le j\le n$ let $\mathcal{T}_j$ have head $z_j$,
and for $1\le i\le n$ let $\mathcal{L}_{in}$ have apex $x_i$. Thus for $i,j\in I_n$, $S_{ij}$ is a path joining $x_i$ and $z_j$,
and we claim that these paths
form an impression of $K_{n,n}$. To show this, we must show:
\\
\\
(1) {\em For all $i,j,i',j'\in I_n$, if $i\ne i'$ and $j\ne j'$ then $V(S_{ij})$ is disjoint from and anticomplete to $V(S_{i'j'})$.}
\\
\\
We may assume that $j<j'$, from the symmetry. Suppose that $v\in V(S_{ij})$ and $v'\in V(S_{i'j'})$ are either equal or adjacent.
Now $v'\in V(S_{i'j'})$ and so $v'$ belongs either to 
$V(\mathcal{L}_{i',j'-1})$, or to $C_{j'-1}$, or to $V(\mathcal{L}_{k,j'-1})$ for some $k\in I_{j'-1}\setminus I_{j'}$.
Hence $v'$ belongs either to
$V(\mathcal{L}_{i'j})$, or to $C_{j}$, or to $V(\mathcal{L}_{kj})$ for some $k\in I_{j}\setminus I_{n}$.
But $\mathcal{T}_j$ is a tick on $(\mathcal{M}_j,C_j)$, and hence
\begin{itemize}
\item $V(S_{ij})\cap (V(\mathcal{M}_j)\cup C_j)=\{x_i\}$, and so $v\ne v'$; and 
\item every vertex in $V(\mathcal{M}_j)\cup C_j$ with a neighbour in $V(S_{ij})$ belongs to
$V(\mathcal{L}_{ij})$.
\end{itemize}
It follows in particular that $v'\in V(\mathcal{L}_{ij})$; but we already showed that $v'$ 
belongs either to
$V(\mathcal{L}_{i'j})$, or to $C_{j}$, or to $V(\mathcal{L}_{kj})$ for some $k\in I_{j}\setminus I_{n}$, a contradiction.
This proves (1).

\bigskip

Since each $S_{ij}$ has length at most $\rho+3$, it follows that $G$ contains an impression of $K_{n,n}$ 
of order at most $\rho+3$, a contradiction. This proves \ref{usetick2again}.~\bbox

The proof of \ref{gettick} breaks into two cases, depending whether $\rho=3$ or not. In this section we handle the easier case $\rho\ge 4$,
and postpone $\rho=3$ until the next section.
When $\rho\ge 4$, a stronger statement holds, the following:

\begin{thm}\label{gettick3}
For all $\rho\ge 4$ and $\tau, m,c'\ge 0$ there exists $c\ge 0$ with the following property.
Let $G$ be a graph such that $\chi^{\rho-1}(G)\le \tau$.
Let $C\subseteq V(G)$ with $\chi(C)>c$, and let
$\mathcal{M}=(\mathcal{L}_i:i\in I)$ be a $\rho$-multicover for $C$, with $|I|= m$. Then
there exist $C'\subseteq C$ with $\chi(C')> c'$, and a $\rho$-multicover
$\mathcal{M}'$ for $C'$ contained in $\mathcal{M}$ with length $m$, and a tick
$(S_i:i\in I)$ on $(\mathcal{M}',C')$ with head $z\in C\setminus C'$, such that for each $i\in I$, $S_i$ has length $\rho$, 
and $V(S_i)\subseteq V(\mathcal{L}_i)\cup \{z\}$ 
(and so the paths $S_i\;(i\in I)$ are pairwise disjoint except for $z$).
\end{thm}
\Proof Let $c=c'+(m(\rho-1)+1)\tau$, and let $G,C$ and $\mathcal{M}=(\mathcal{L}_i:i\in I)$ be as in the theorem. Let $x_i$ be the apex of $\mathcal{L}_i$
for each $i\in I$, and let $X=\{x_i:i\in I\}$.
For each $i\in I$, let $C_i$ be the set of vertices in $C$ with distance at most $\rho-1$ from $x_i$ in $G$. Then by hypothesis,
$\chi(C_i)\le \tau$; let $D$ be the set of vertices in $C$ that do not belong to the union of the sets $C_i\;(i\in I)$.
It follows that $\chi(D)>c-m\tau$. Since $c\ge m\tau$, there exists $z\in D$; choose some such $z$.
For each $i\in I$ let $S_i$ be some $\mathcal{L}_i$-radius for $z$. 
\\
\\
(1) {\em For all distinct $i,i'\in I$, $x_{i'}$ has no neighbours in $V(S_i)$.}
\\
\\
Suppose that some $x_{i'}$ is adjacent to a vertex in $S_i$. Since
$S_i$ has length $\rho$, and the distance from $x_{i'}$ to $z$ is at least $\rho$ (because $z\notin C_{i'}$), 
it follows that $x_{i'}$ is adjacent to $x_i$ or to the neighbour of $x_i$ in $S_i$; but this contradicts that $\mathcal{M}$
is a multicover, since $\rho\ge 4$. This proves (1).

\bigskip

Let $S$ be the union of the sets $V(S_i)\;(i\in I)$.
Thus $|S|=m(\rho-1)+1$. Let $C'$ be the set of vertices in $C$ with distance at least $\rho$ in $G$ from every vertex in $S$. Since
$X\subseteq S$ it follows that $C'\subseteq D$, and $z\in D\setminus C'$, and $\chi(C')>c-(m(\rho-1)+1)\tau=c'$. 
For each $j\in I$, let $\mathcal{L}_j=(L_{0,j}\ll L_{\rho-1,j})$ say, and for $0\le i\le \rho-1$ let $L_{i,j}'$ be the set of vertices $v\in L_{i,j}$
such that some $\mathcal{L}_j$-radius contains both $v$ and a vertex in $C'$; and let $\mathcal{L}'_j = (L_{0,j}'\ll L_{\rho-1,j}')$. Then
$\mathcal{L}'_j$ is a $\rho$-cover for $C'$; let $\mathcal{M}' = (\mathcal{L}_j':j\in I)$, and then $\mathcal{M}'$
is a $\rho$-multicover for $C'$ contained in $\mathcal{M}$. We claim that it satisfies the theorem.
Certainly $z\in C\setminus C'$. 
\\
\\
(2) {\em $V(S_i)\cap V(\mathcal{M}') = \{x_i\}$ for each $i\in I$. }
\\
\\
For suppose that $u\in V(S_j)\cap V(\mathcal{M}')$, and choose $j'\in I$ so that
$u\in V(\mathcal{L}_{j'}')$. Since $V(S_j)\subseteq V(\mathcal{L}_j)$ and $V(\mathcal{L}_{j'}')\subseteq V(\mathcal{L}_{j'})$, it follows that 
$V(\mathcal{L}_j)$ is not disjoint from $V(\mathcal{L}_{j'})$, and so $j'=j$. Since $u\in  V(\mathcal{L}_{j}')$, there exists $i$ with $0\le i\le \rho-1$
such that $u\in L_{i,j}'$; and so the distance in $G$ between $u$ and some vertex in $C'$ is at most $\rho-i$. But from the definition
of $C'$, since $u\in S$ it follows that this distance is at least $\rho$, and so $i=0$, that is, $u=x_j$. This proves (2).
\\
\\
(3) {\em For each $j\in I$, if some $u\in V(S_j)$ is adjacent to some $v\in  V(\mathcal{M}')\cup C'$ then
$v\in V(\mathcal{L}_j')$.}
\\
\\
Assume that $u\in V(S_j)$ and $v\in  V(\mathcal{M}')\cup C'$ are adjacent. Since $u\in S$
and so has distance at least $\rho$ from every vertex in $C'$, it follows that $v\notin C'$, and so $v\in V(\mathcal{L}_{j'}')$ for some $j'\in I'$.
Choose $i$ so that $v\in L_{i,j'}'$; then the distance in $G$ between $v$ and some vertex in $C'$ is at most $\rho-i$, and so
the distance between $u$ and some vertex in $C'$ is at most $\rho+1-i$. Since this distance is at least $\rho$, it follows
that $i\le 1$, and so $v$ is equal to or adjacent to $x_{j'}$, and in either case $v$ does not belong to the base of $\mathcal{L}_{j'}$.
If $u$ belongs to the base of $\mathcal{L}_j$, then $u$ is adjacent
to $z$ (because only one vertex in $S_j$ belongs to the base of $\mathcal{L}_j$, namely the neighbour of $z$); and since $i\le 1$,
and therefore the distance between $u$ and $x_{j'}$ in $G$ is at most $2$, it follows that the distance between $z$ and $x_{j'}$ is at most $3$,
contrary to the definition of $D$ (since $\rho\ge 4$). Thus $u$ does not belong to the base of $\mathcal{L}_j$; and since
$\mathcal{M}$ is a multicover, it follows that $j=j'$. This proves (3).

\bigskip
From (1), (2) and (3) it follows that $(S_i:i\in I)$ is a tick on $(\mathcal{M}',C')$. This proves \ref{gettick3}.~\bbox

\section{Extracting ticks from $3$-multicovers}

In this section we prove \ref{gettick} when $\rho=3$.
We will need the following lemma, proved in~\cite{dsw}:
\begin{thm}\label{dswlemma}
Let $\mathcal{A}$ be a set of nonempty subsets of a finite set $V$, and let $k\ge 0$ be an integer. Then either:
\begin{itemize}
\item there exist $A_1,A_2\in \mathcal{A}$ with $A_1\cap A_2=\emptyset$;
\item there are $k$ distinct members $A_1\ll A_k\in \mathcal{A}$, and for all $i,j$ with $1\le i<j\le k$ an element $v_{ij}\in V$,
such that for all $h,i,j\in \{1\ll k\}$ with $i<j$, $v_{ij}\in A_h$ if and only if $h\in \{i,j\}$; or
\item there exists $X\subseteq V$ with $|X|\le 11(k+4)^5$ such that $X\cap A\ne \emptyset$ for all $A\in \mathcal{A}$.
\end{itemize}
\end{thm}

The idea of using \ref{dswlemma} in this context is due to Bousquet and Thomass\'e~\cite{stephan}.
We use it to prove the following.

\begin{thm}\label{farapart}
For all $\mu,\nu\ge 0$, there exists $m\ge 0$ with the following property. Let $G$ be $(1,\mu,\nu)$-restricted,
and let $X\subseteq V(G)$, such that every two vertices in $X$ have distance at most two in $G$. Then there exists $Y\subseteq V(G)$
with $|Y|\le m$ such that every vertex in $X\setminus Y$ has a neighbour in $Y$.
\end{thm}
\Proof Choose $k$ so that \ref{stephanproof} holds with $m$ replaced by $k$, and let $m=11(k+4)^5$.
We claim that $m$ satisfies the theorem; for let $G,X$ be as in the theorem.
For each $x\in X$, let $N[x]$ be the set of all vertices equal to or adjacent in $G$ to $x$, and let $\mathcal{A}$ be the set $\{N[x]:x\in X\}$. 
By hypothesis, no two members of $\mathcal{A}$ are disjoint. Suppose that $A_1\ll A_k\in \mathcal{A}$ are distinct, where $A_i=N[x_i]$
for $1\le i\le k$; then by \ref{stephanproof} and the choice of $k$, there exist $i,j$ with $1\le i<j\le k$ such that
$x_i, x_j$ are nonadjacent, and every vertex of $G$ adjacent to both $x_i, x_j$ has a third neighbour in $\{x_1\ll x_k\}$.
Consequently there is no vertex $v_{ij}$ in $V(G)$
such that for all $h\in \{1\ll k\}$ with $i<j$, $v_{ij}\in A_h$ if and only if $h\in \{i,j\}$. 

From \ref{dswlemma} we deduce that there exists $Y\subseteq V$ with $|Y|\le 11(k+4)^5=m$ such that $Y\cap A\ne \emptyset$ 
for all $A\in \mathcal{A}$. But then every vertex in $X$ either belongs to $Y$ or has a neighbour in $Y$. This proves \ref{farapart}.~\bbox

If $\mathcal{M}=(\mathcal{L}_i:i\in I)$ is a $3$-multicover of $C$, and $i,j\in I$ are distinct, and $z\in C$, let $P,Q$ be $\mathcal{L}_i$-
and $\mathcal{L}_j$-radii for $z$ respectively; 
then $P\cup Q$ is a path of $G$ (not necessarily induced), and we call such a path an {\em $(\mathcal{L}_i,\mathcal{L}_j)$-diameter}.
We need another lemma.
\begin{thm}\label{gettick2lemma}
For all $\mu, \nu,\tau, c'\ge 0$ and $m>0$ there exist $c\ge 0$ with the following property.
Let $G$ be a $(1,\mu,\nu)$-restricted graph such that $\chi^2(G)\le \tau$.
Let $C\subseteq V(G)$ with $\chi(C)>c$, and let
$\mathcal{M}=(\mathcal{L}_i:i\in I)$ be a $3$-multicover for $C$ with $|I|=m$. Let $x_i$ be the 
apex of $\mathcal{L}_i$ for $i\in I$. Let $k\in I$ be maximum. For each $g\in I\setminus \{k\}$,
there exist 
\begin{itemize}
\item a subset $I'\subseteq I\setminus \{k\}$ with $|I'|\ge m/2$ and with $\{i\in I: i\le g\}\subseteq I'$;
\item a subset $C'\subseteq C$ with $\chi(C')>c'$;
\item for each $i\in I'$, a $3$-cover $\mathcal{L}_{i}'$ for $C'$ contained in $\mathcal{L}_{i}$, such that for all distinct
$i,i'\in I'$, $x_i$ has no neighbour in $V(\mathcal{L}_{i'})$; and
\item an $(\mathcal{L}_g,\mathcal{L}_k)$-diameter $S$, such that $V(S)$ is anticomplete to $C'$, and
$V(S)$ is anticomplete to $V(\mathcal{L}'_{i})$ for each $i\in I'\setminus \{g\}$, and
$V(S)\cap V(\mathcal{L}'_{g})=\{x_g\}$, and $V(S)\subseteq V(\mathcal{L}_{g})\cup V(\mathcal{L}_{k})\cup C$.
\end{itemize}
\end{thm}
\Proof Choose $m_0$ so that \ref{farapart} holds with $m$ replaced by $m_0$. Let 
$$c=\max((m+m_0)\tau, (12+m)\tau+c'2^{m+1}).$$
We claim that $c$ satisfies the theorem. For let $G,C$, $\mathcal{M}=(\mathcal{L}_i:i\in I),k,g$
be as in the theorem, where $\mathcal{L}_{i} = (\{x_i\}, A_i, B_i)$ for each $i\in I$, say.
Since the set of vertices in $C$ with distance at most two from one of the vertices $x_i\;(i\in I)$ has chromatic number
at most $m\tau$, there exists $C_0\subseteq C$ with $\chi(C_0)> c-m\tau$ such that every vertex in $C_0$ has distance at 
least three from each $x_i$. Let $D$ be the set of vertices in $B_g$ with a neighbour in $C_0$.
\\
\\
(1) {\em There exist $y_1,y_2\in D$ with distance
at least three in $G$.}
\\
\\
For if not, then by \ref{farapart} applied with $X=D$, there exists $Y\subseteq V(G)$ with $|Y|\le m_0$ such that 
every vertex in $D\setminus Y$ has a neighbour in $Y$. Then every vertex in
$C_0$ has distance at most two from a vertex in $Y$, and so $\chi(C_0)\le |Y|\tau$; and since $\chi(C_0)>c-m\tau$, it follows
that $|Y|>c\tau^{-1}-m\ge m_0$, a contradiction. This proves (1). 

\bigskip

Choose $z_1,z_2\in C_0$
adjacent to $y_1,y_2$ respectively. Let $S_1$ be an $(\mathcal{L}_g,\mathcal{L}_k)$-diameter containing $y_1$ and $z_1$,
and choose $S_2$ for $y_2,z_2$ similarly. The union of $S_1$ and $S_2$ has at most $12$ vertices, and so the set
of vertices in $C_0$ with distance at most two from a vertex in $S_1\cup S_2$ has chromatic number at most $12\tau$.
Consequently there exists $C_1\subseteq C_0$ with $\chi(C_1)>c-m\tau-12\tau$ such that every vertex in $C_1$ has distance at least 
three from every vertex in $S_1\cup S_2$. For $1\le i\le g$, let $\mathcal{L}_{i}'$ be the levelling $(\{x_i\},A_i', B_i')$, where 
$B_i'$ is the set of vertices in $B_i$ with a neighbour in $C_1$, and $A_i'$ is the set of vertices in $A_i$ with a neighbour in $B_i'$.
Then $V(S_1\cup S_2)\cap V(\mathcal{L}_{g}')=\{x_g\}$, because every vertex in $C_1$ has distance at least three from $S_1\cup S_2$.
Also $V(S_1\cup S_2)$ is anticomplete to $V(\mathcal{L}_{i}')$ if $i<g$, since every vertex in $V(\mathcal{L}_{i})$ with a neighbour in
$S_1\cup S_2$ belongs to $B_i$ (from the definition of a $3$-multicover) and hence does not belong to $B_i'$ (because vertices in $B_i'$
have neighbours in $C_1$ and therefore have no neighbours in $S_1\cup S_2$). Also, for $j\in I$ with $j\ne g,k$, $x_j$ has 
no neighbour in $S_1\cup S_2$ (from the definition of a multicover, and since $z_1,z_2\in C_0$ and therefore have distance
at least three from $x_j$). Moreover, 
$$V(S_1\cup S_2)\subseteq V(\mathcal{L}_g)\cup V(\mathcal{L}_k)\cup C.$$

Now we shall choose one of $S_1,S_2$ to satisfy the other requirements of the theorem.
For each $j\in I\setminus \{k\}$ with $j>g$ and each $v\in C_1$, let $P_{jv}$ be an $\mathcal{L}_{j}$-radius for $v$. 
Fix $v\in C_1$ for the moment. 
Now $P_{jv}$ has length three;
let its vertices be $x_j\d  a_{jv}\d  b_{jv}\d v$ in order. We have seen that
$x_j$ has no neighbours in $S_1\cup S_2$.  Since $v\in C_1$ and therefore has distance at least three from every 
vertex in $S_1\cup S_2$, it follows that $v, b_{jv}$ have no neighbours in $S_1\cup S_2$; but $a_{jv}$ might have neighbours in $S_1\cup S_2$.
From the definition of a multicover, every neighbour of $a_{jv}$ in $S_1\cup S_2$ is one of $y_1,y_2$; and since $y_1,y_2$
have distance at least three in $G$, $a_{jv}$ is not adjacent to them both. Consequently $V(P_{jv})$ is anticomplete to
at least one of $S_1, S_2$. Choose $I_v\subseteq I\setminus \{k\}$ including $\{i\in I:i\le g\}$, with $|I_v|\ge m/2$,
such that for one of $S_1, S_2$ (say $S_v$),
each of the paths $P_{jv}\;(j\in I_v, j>g)$ is anticomplete to $S_v$. There are only $2^{m+1}$ possibilities for the pair
$(S_v,I_v)$; and so there exists $C'\subseteq C_1$ with $\chi(C')\ge \chi(C_1)2^{-m-1}>c'$, and one of $S_1, S_2$, say $S$, and a set $I'$,
such that $S_v=S$ and $I_v=I'$ for all $v\in C'$. For each $j\in I\setminus \{k\}$ with $j>g$, let $\mathcal{L}_{j}'$
be the levelling $(\{x_j\},A_j', B_j')$, where $A_j'=\{a_{jv}:v\in C'\}$ and $B_j' = \{b_{jv}:v\in C'\}$. 

We claim that for all distinct $i,i'\in I'$, $x_i$ has no neighbour in $V(\mathcal{L}_{i'}')$. Suppose it does; then $i>i'$ and
$x_i$ has a neighbour in $B_{i'}'$. But every vertex in $B_{i'}'$ has a neighbour in $C_1\subseteq C_0$, and the distance between 
$x_i$ and every vertex in $C_0$ is at least three, a contradiction. This proves the claim, and so proves \ref{gettick2lemma}.~\bbox

We deduce:

\begin{thm}\label{gettick2lemma2}
For all $\mu, \nu,\tau, c'\ge 0$, and $t>0$, and $m\ge t2^t$, there exist $c\ge 0$ with the following property.
Let $G$ be a $(1,\mu,\nu)$-restricted graph such that $\chi^2(G)\le \tau$.
Let $C\subseteq V(G)$ with $\chi(C)>c$, and let
$\mathcal{M}=(\mathcal{L}_i:i\in I)$ be a $3$-multicover for $C$ with $|I|=m$. Let $k\in I$ be maximum.
Then there exist
\begin{itemize}
\item a subset $I'\subseteq I\setminus \{k\}$ with $|I'|\ge m2^{-t}\ge t$; $I' = \{i_1\ll i_n\}$ say, where $i_1<i_2<\cdots<i_n$;
\item  a subset $C'\subseteq C$ with $\chi(C')>c'$;
\item for each $i\in I'$, a $3$-cover $\mathcal{L}_{i}'$ for $C'$, contained in $\mathcal{L}_{i}$;
\item for each $i\in \{i_1\ll i_t\}$, an $(\mathcal{L}_i,\mathcal{L}_k)$-diameter $S_i$, such that
$V(S_i)$ is anticomplete to $C'$,  and $V(S_i)$ is anticomplete to
$V(\mathcal{L}_{j}')$ for all $j\in I'\setminus \{i\}$, and $V(S_i)\cap V(\mathcal{L}_{i}')=\{x_{i}\}$, and
$V(S_i)\subseteq V(\mathcal{L}_{i})\cup V(\mathcal{L}_{k})\cup C$.
\end{itemize}
\end{thm}
\Proof We assume first that $t=1$. 
Choose $c$ so that \ref{gettick2lemma} is satisfied.
Choose $g\in I$, minimum; then the result follows from
\ref{gettick2lemma}. Thus the result holds if $t=1$.

We fix $\mu, \nu,\tau,m$, and proceed by induction on $t$ (assuming $m\ge t2^t$). Thus we assume that $t>1$ and the 
result holds with $t$ replaced by $t-1$.
Choose $c''$ so that \ref{gettick2lemma} is satisfied with $c$ replaced by $c''$ (and the given value of $m$).
Let $c$ have the value that satisfies the theorem with $t,c'$ replaced by $t-1,c''$; we claim that $c$ satisfies the theorem.

For let $G,C$ and $\mathcal{M}=(\mathcal{L}_i:i\in I),k$ be as in the theorem, where $|I|=m\ge t2^t$. 
From the inductive hypothesis, there exist
\begin{itemize}
\item a subset $I''\subseteq I\setminus \{k\}$ with $|I''|\ge m2^{1-t}$; $I'' = \{i_1\ll i_n\}$ say, where $i_1<i_2<\cdots<i_n$;
\item  a subset $C''\subseteq C$ with $\chi(C'')>c''$;
\item for each $i\in I''$, a $3$-cover $\mathcal{L}_{i}''$ for $C''$, contained in $\mathcal{L}_{i}$;
\item for each $i\in \{i_1\ll i_{t-1}\}$, an $(\mathcal{L}_i,\mathcal{L}_k)$-diameter $S_i$, such that
$V(S_i)$ is anticomplete to $C''$, and $V(S_i)$ is anticomplete to
$V(\mathcal{L}_{j}'')$ for all $j\in I''\setminus \{i\}$, and $V(S_i)\cap V(\mathcal{L}_{i}'')=\{x_{i}\}$, 
and $V(S_i)\subseteq V(\mathcal{L}_i)\cup V(\mathcal{L}_k)\cup C$.
\end{itemize} 
Let $\mathcal{L}_{k}''=\mathcal{L}_{k}$.
Thus $\mathcal{M}'' = (\mathcal{L}_i'':i\in I''\cup \{k\})$ is a $3$-multicover of $C''$, contained in $\mathcal{M}$. Also
$n\ge 2t$, since $n\ge m2^{1-t}$ and $m\ge t2^t$. From 
\ref{gettick2lemma} applied to $\mathcal{M}''$ taking $g=i_t$, we deduce that there exist
\begin{itemize}
\item a subset $I'\subseteq I''$ with $|I'|\ge (|I''|+1)/2\ge m2^{-t}$ and with $\{i_1\ll i_t\}\subseteq I'$;
\item a subset $C'\subseteq C''$ with $\chi(C')>c'$;
\item for each $i\in I'$, a $3$-cover $\mathcal{L}_{i}'$ for $C'$ contained in $\mathcal{L}_{i}''$;
\item an $(\mathcal{L}_{i_t}'',\mathcal{L}_k'')$-diameter $S_{i_t}$ (which is therefore also an $(\mathcal{L}_{i_t},\mathcal{L}_k)$-diameter), 
such that $V(S_{i_t})$ is anticomplete to $C'$, and
$V(S_{i_t})$ is anticomplete to $V(\mathcal{L}_{i}')$ for all $i\in I'\setminus \{i_t\}$, and 
$V(S_{i_t})\cap V(\mathcal{L}_{i_t}')=\{x_{i_t}\}$, and $V(S_{i_t})\subseteq V(\mathcal{L}_{i_t})\cup V(\mathcal{L}_{k})\cup C$.
\end{itemize}
But then $I', C'$, $\mathcal{L}_{i}'\;(i\in I')$, and the paths
$S_i\;(i\in \{i_1\ll i_t\})$ satisfy the theorem. This proves \ref{gettick2lemma2}.~\bbox

Now we prove the main result of this section, the case of \ref{gettick} for $3$-multicovers:
\begin{thm}\label{gettick2}
For all $\mu, \nu,\tau, m', c'\ge 0$ there exist $m,c\ge 0$ with the following property.
Let $G$ be a $(1,\mu,\nu)$-restricted graph such that $\chi^2(G)\le \tau$.
Let $C\subseteq V(G)$ with $\chi(C)>c$, and let
$\mathcal{M}=(\mathcal{L}_i:i\in I)$ be a $3$-multicover for $C$, with length $m$. Let $k\in I$ be maximum. Then
there exist $C'\subseteq C$ with $\chi(C')> c'$, and a $3$-multicover
$\mathcal{M}'$ for $C'$ contained in $\mathcal{M}$ with length $m'$, with index set some $I'\subseteq I\setminus \{k\}$, and a tick
$(S_i:i\in I')$ on $(\mathcal{M}',C')$ of order at most $6$, such that $V(S_i)\subseteq V(\mathcal{L}_i)\cup V(\mathcal{L}_k)\cup C$
for each $i\in I'$.
\end{thm}
\Proof Let $m=m'2^{m'}$ and let $c$ satisfy \ref{gettick2lemma2} with this choice of $m$, taking $t=m'$. 
We claim that $m,c$ satisfy the theorem. For let $G,C$, $\mathcal{M}=(\mathcal{L}_i:i\in I)$ and $k$
be as in the theorem. For each $i\in I$, let $\mathcal{L}_i = (\{x_i\}, A_i, B_i)$ say. 

By \ref{gettick2lemma2} applied to $\mathcal{M}$,
there exist
\begin{itemize}
\item a subset $I'\subseteq I\setminus \{k\}$ with $|I'|=|I|2^{-m'}=m'$ (we only take the first $m'$ elements of the 
set $I'$ claimed by \ref{gettick2lemma2}); 
\item  a subset $C'\subseteq C$ with $\chi(C')>c'$;
\item for each $i\in I'$, a $3$-cover $\mathcal{L}_{i}'$ for $C'$, contained in $\mathcal{L}_{i}$;
\item for each $i\in I'$, an $(\mathcal{L}_i,\mathcal{L}_k)$-diameter $S_i$, such that
$V(S_i)$ is anticomplete to $C'$, and $V(S_i)$ is anticomplete to
$V(\mathcal{L}_{j}')$ for all $j\in I'\setminus \{i\}$, and and $V(S_i)\cap V(\mathcal{L}_{i}')=\{x_{i}\}$, and
$V(S_i)\subseteq V(\mathcal{L}_i)\cup V(\mathcal{L}_k)\cup C$.
\end{itemize}
Let $\mathcal{M}'=(\mathcal{L}_{i}':i\in I')$. Then $\mathcal{M}'$ is a $3$-multicover of $C'$, and
$(S_i:i\in I')$ is a tick on $(\mathcal{M}',C')$ of order at most six, with head $x_k$. This proves \ref{gettick2}.~\bbox

Together \ref{gettick2} and \ref{gettick3} imply \ref{gettick}, so we have completed the proof of \ref{gettick}, and hence of 
\ref{usetick2}, \ref{reducecontrol} and \ref{3control}. Henceforth we need only consider 2-controlled ideals of graphs.

\section{Clique control}

Now we come to the second part of the paper, in which we handle 2-controlled graphs. 
We will follow the approach taken in 
\cite{longholes}; and in particular, it will be helpful to introduce a refinement of control, called ``clique-control''. 
If $X$ is a clique with $|X|=\xi$ we call $X$ a {\em $\xi$-clique.}
We denote by $N^1_G(X)$ the set of all vertices in $V(G)\setminus X$ that are complete to $X$;
and by $N^2_G(X)$ the set of all vertices in $V(G)\setminus X$ with a neighbour in $N^1(X)$ and with no neighbour in $X$.
When $X=\{v\}$ we write $N^i_G(v)$ for $N^i_G(X)$ ($i = 1,2$). (We omit the subscript $G$ when the graph is clear from context.)
We are assuming that in every induced subgraph $H$ of large $\chi$,
there is a vertex $v$ such that $N^2_H(v)$ also has large $\chi$; and perhaps the same is true for cliques larger than singletons.
For instance, it may or may not be true that in every induced subgraph $H$ of large $\chi$,
there is a $2$-clique $X$ such that $N^2_H(X)$ also has large $\chi$. If
this is false, we can find $H$ in the ideal with arbitrarily large $\chi$ such that $N^2_H(X)$ has bounded $\chi$ for all
$2$-cliques $X$, and we focus on these. If it is true, then we ask the same question for triples, and so on; we must
soon hit a clique-size for which the answer is ``false'', because none of our graphs have a clique larger than $\nu$. 
Let us say this more precisely.

Let $\phi:\mathbb{N}\rightarrow \mathbb{N}$ be a nondecreasing function, and
let $\xi\ge 1$ be an integer. We say a graph $G$ is {\em $(\xi,\phi)$-clique-controlled} if for every induced subgraph $H$ of $G$ and
every integer $n\ge 0$, if $\chi(H)>\phi(n)$ then there is a $\xi$-clique $X$ of $H$ such that $\chi(N^2_H(X))>n$.
Roughly, this means that in every induced subgraph $H$ of large chromatic number, there is a $\xi$-clique
$X$ with $N^2_H(X)$ of large chromatic number. We say an ideal of graphs $\mathcal{C}$ is {\em $\xi$-clique-controlled} if
there is a nondecreasing function $\phi$ such that every graph in $\mathcal{C}$ is $(\xi,\phi)$-clique-controlled.
An ideal $\mathcal{C}$ of graphs is {\em colourable} if there is an integer $k$ such that every graph in $\mathcal{C}$
has chromatic number at most $k$; and {\em non-colourable} if there is no such $k$.
\begin{thm}\label{maxclique}
Let $\nu\ge 1$ and $\tau_1\ge 0$, and let $\mathcal{C}$ be a non-colourable ideal of graphs such that 
\begin{itemize}
\item $\mathcal{C}$ is $2$-controlled;
\item $\omega(G)\le \nu$ for each $G\in \mathcal{C}$; and
\item  $\chi(G)\le \tau_1$ for every $G\in \mathcal{C}$ with $\omega(G)<\nu$.
\end{itemize}
Then there exists $\xi$
with $1\le \xi \le \nu$ such that $\mathcal{C}$ is $\xi$-clique-controlled; and there is a non-colourable subideal $\mathcal{C}'$
of $\mathcal{C}$ and  $\tau_2\ge 0$ such that 
$\chi(N^2_G(X))\le \tau_2$ for every $G\in \mathcal{C}'$ and for every $(\xi+1)$-clique $X$ of $G$.
\end{thm}
\Proof Suppose that $\mathcal{C}$ is $\nu$-clique-controlled, and choose a function $\phi$ such that every graph in $\mathcal{C}$
is $(\nu,\phi)$-clique-controlled. Let $c=\phi(0)$; then by hypothesis, there exists $G\in \mathcal{C}$ with $\chi(G)>c$.
From the definition of $(\nu,\phi)$-clique-controlled, there is a $\nu$-clique $X$ in $G$ with $\chi(N^2_G(X))>0$, which is 
impossible since $N^1(X)=\emptyset$ (because $\omega(G)\le \nu$).

This proves that $\mathcal{C}$ is not $\nu$-clique-controlled. We claim that $\mathcal{C}$ is $1$-clique-controlled.
Choose $\phi$ such that every graph in $\mathcal{C}$ is $(2,\phi)$-controlled, and let $\phi'(c) = \phi(c+\tau_1+1)$ for each $c\ge 0$.
We claim that every $G\in \mathcal{C}$ is $(1,\phi)$-clique-controlled. For let $c\ge 0$, and let
$H$ be an induced subgraph of $G\in \mathcal{C}$, with $\chi(H)>\phi'(c)$. Then $\chi(H)> \phi(c+\tau_1+1)$, and since
$G$ is $(2,\phi)$-controlled, it follows that
$\chi^2(H)>c+\tau_1+1$. Hence there is a vertex $v$ of $H$ such that $\chi(N^2_H[v])>c+\tau_1+1$. 
Now $\chi(N^1_H[v])\le \tau_1+1$, since
the subgraph of $H$ induced on $N^1_H(v)$ has clique number at most $\nu-1$. Consequently
$\chi(N^2_H(v))>c$. This proves that $\mathcal{C}$ is $1$-clique-controlled.

Choose $\xi$ maximum such that $\mathcal{C}$ is $\xi$-clique-controlled; then $1\le \xi< \nu$.
Suppose that for all $\kappa\ge 0$, there exists $m_{\kappa}$ such that for every $G\in \mathcal{C}$
with $\chi(G)>m_{\kappa}$, there is a $(\xi+1)$-clique $X$ of $G$ with $\chi(N^2_G(X))>\kappa$. Then every member of $\mathcal{C}$
is $(\xi+1,\phi')$-clique-controlled,
where we define $\phi'(\kappa) = m_{\kappa}$ for each $\kappa\ge 0$ (having arranged that $m_0\le m_1\le \ldots$). Consequently $\mathcal{C}$
is $(\xi+1)$-clique-controlled, a contradiction. 

Thus there exists $\kappa\ge 0$ 
such that for all $c$, there are graphs $G\in \mathcal{C}$ such that $\chi(G)>c$
and $\chi(N^2_G(X))\le \kappa$ for every $(\xi+1)$-clique $X$ of $G$. Let $\tau_2 = \kappa$, and let
$\mathcal{C}'$ be the subideal of all graphs $G\in \mathcal{C}$
such that $\chi(N^2_G(X))\le \tau_2$ for every $(\xi+1)$-clique $X$ of $G$. Then $\mathcal{C}'$ is non-colourable.
This proves \ref{maxclique}.~\bbox

The advantage of looking at an ideal of graphs that is $\xi$-clique-controlled is the following. Start with a graph in the ideal
with huge chromatic number. Consequently it contains a $\xi$-clique $X_1$ with $\chi(N^2(X_1))$ (not quite so) huge;
let $C_1=N^2(X_1)$. Since $\chi(C_1)$ is huge, there is a 
$\xi$-clique $X_2$ of $G_1=G[C_1]$ such that $\chi(N^2_{G_1}(X_2))$ is fairly huge; and so on. We generate a sequence
of ``$\xi$-clique-covers'' of some ultimate set $C$, of any desired length, and this gives us some structured thing 
to explore in the hope of finding the induced subgraph we want. We call this
a ``$\xi$-clique-multicover'' of $C$.

Formally:
let $G$ be a graph, and $X,N, C\subseteq V(G)$, such that
\begin{itemize}
\item $X,N, C$ are pairwise disjoint;
\item $X$ is a $\xi$-clique; 
\item $X$ is complete to $N$; 
\item $X$ is anticomplete to $C$; and
\item $N$ covers $C$.
\end{itemize}
We say that the pair
$\mathcal{L} = (X,N)$ is a {\em $\xi$-clique-cover} of $C$. We write
$X(\mathcal{L}) = X$, $N(\mathcal{L}) = N$, and $V(\mathcal{L})=X\cup N$.
Thus $(X,N)$ is a $1$-clique-cover of $C$ if and only if $(X,N)$ is a $2$-cover for $C$.

A {\em $\xi$-clique-multicover}
of $C$ of {\em length} $|I|$ is a family $(\mathcal{L}_i:i\in I)$ of $\xi$-clique-covers of $C$, where $I$ is a set of integers, such that:
\begin{itemize}
\item the sets $V(\mathcal{L}_i)\;(i \in I)$ are pairwise disjoint; and
\item for all $i,j\in I$ with $i<j$, $X(\mathcal{L}_i)$ is anticomplete to $V(\mathcal{L}_j)$.
\end{itemize}

\begin{figure}[h]
\centering
\tikzset{snake it/.style={decorate, decoration=snake}}
\begin{tikzpicture}[scale=.7,auto=left]
\node at (-4,1) {$X(\mathcal{L}_1)$};
\node at (4,1) {$X(\mathcal{L}_2)$};
\node at (-4,-2) {$N(\mathcal{L}_1)$};
\node at (4,-2) {$N(\mathcal{L}_2)$};
\node at (0,-4) {$C$};
\tikzstyle{every node}=[inner sep=1.5pt, fill=black,circle,draw]
\node (x1) at (-4,0) {};
\node (x2) at (-4.5,0.37) {};
\node (x3) at (-3.5, 0.37) {};
\node (y1) at (4,0) {};
\node (y2) at (4.5,0.37) {};
\node (y3) at (3.5, 0.37) {};

\draw (-4,-2) ellipse (2cm and .6cm);
\draw (4,-2) ellipse (2cm and .6cm);
\draw (0,-4) ellipse (3cm and 1cm);
\draw (-4,0) -- (-4,-1.6);
\draw (-4,0) -- (-3,-1.6);
\draw (-4,0) -- (-5,-1.6);
\draw (-4.5,0.37) -- (-4,-1.6);
\draw (-4.5,0.37) -- (-3,-1.6);
\draw (-4.5,0.37) -- (-5,-1.6);
\draw (-3.5, 0.37) -- (-4,-1.6);
\draw (-3.5, 0.37) -- (-3,-1.6);
\draw (-3.5, 0.37) -- (-5,-1.6);
\draw (-4,0) -- (-4.5,0.37);
\draw (-4,0) -- (-3.5, 0.37);
\draw (-4.5,0.37) -- (-3.5, 0.37);

\draw (4,0) -- (4,-1.6);
\draw (4,0) -- (3,-1.6);
\draw (4,0) -- (5,-1.6);
\draw (4.5,0.37) -- (4,-1.6);
\draw (4.5,0.37) -- (3,-1.6);
\draw (4.5,0.37) -- (5,-1.6);
\draw (3.5, 0.37) -- (4,-1.6);
\draw (3.5, 0.37) -- (3,-1.6);
\draw (3.5, 0.37) -- (5,-1.6);
\draw (4,0) -- (4.5,0.37);
\draw (4,0) -- (3.5, 0.37);
\draw (4.5,0.37) -- (3.5, 0.37);
\draw (-2.8, -2.4) -- (0.2,-3.1);
\draw (-4.9, -2.45) -- (-1.9,-3.35);
\draw (2.8, -2.4) -- (-0.2,-3.1);
\draw (4.9, -2.45) -- (1.9,-3.35);
\path[draw=black, snake it] (-2.7, -2.1) -- (2.7, -2.1);
\path[draw=black, snake it] (-3.0, -1.9) -- (4,0.2);

\end{tikzpicture}

\caption{A 3-clique-multicover of length two (wiggly lines indicate possible edges).} \label{fig:3}
\end{figure}
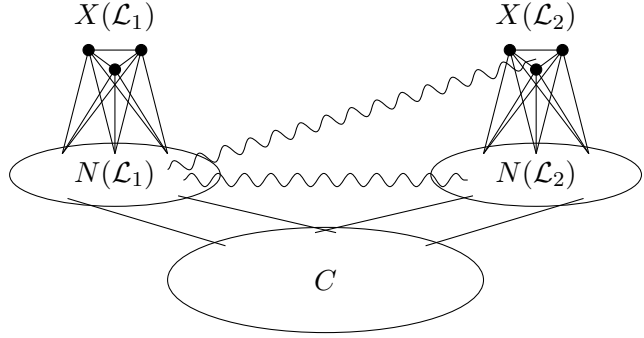

For $i,j\in I$ with $i<j$, we say that the pair $(\mathcal{L}_i,\mathcal{L}_j)$ is {\em independent (with respect to $C$)} if
there exists $x_j\in X(\mathcal{L}_j)$ such that no vertex in $N(\mathcal{L}_i)$ with a neighbour in $C$ is adjacent to $x_j$.
A $\xi$-clique-multicover $\mathcal{M} = (\mathcal{L}_i:i\in I)$ of $C$ is {\em independent} if 
all its pairs $(\mathcal{L}_i,\mathcal{L}_j)$ (where $j>i$) are independent with respect to $C$.
For brevity, let us say a graph $G$ is {\em $(\xi,\zeta,c)$-free}
if for each $C\subseteq V(G)$ with $\chi(C)>c$, there is no independent $\xi$-clique-multicover in $G$ of $C$ with length $\zeta$.

In \cite{longholes} we proved something like \ref{usetick2} for $\rho =2$, but it only applies to 
``strongly-independent'' 2-multicovers.
Let us say
a 2-multicover $\mathcal{M}=(\mathcal{L}_i:i\in I)$ is {\em strongly-independent} if for all $i,j\in I$ with $i<j$, the apex of $\mathcal{L}_j$
has no neighbour in the base of $\mathcal{L}_i$.
(Thus, any edge between $V(\mathcal{L}_i)$ and $V(\mathcal{L}_j)$
is between the two bases, so this is the same as independence as 1-clique-covers, 
except we are also forbidding vertices in $N(L_i)$ that have no neighbour in $C$ and are adjacent to the apex of $\mathcal{L}_j$.) 
A warning: in~\cite{longholes} we used the term ``multicover'' to mean what in this paper is called
a strongly-independent 2-multicover. The result of~\cite{longholes} that we need is the following, theorem 2.3 of that paper.

\begin{thm}\label{usetick1}
For all $n,\nu,\tau_1\ge 0$ there exist $m,d\ge 0$ with the following property.
Let $G$ be a graph with $\omega(G)\le \nu$, such that there is no impression of $K_{n,n}$ in $G$ of order two, and 
$\chi(H)\le \tau_1$
for every induced subgraph $H$ of $G$ with $\omega(H)<\nu$. If  $C\subseteq V(G)$ with $\chi(C)>d$,
then there is no strongly-independent 2-multicover
of $C$ in $G$ with length $m$.
\end{thm}

In view of \ref{impression}, we can strengthen this to:
\begin{thm}\label{betterusetick1}
For all $\mu,\nu,\tau_1\ge 0$ there exist $m,d\ge 0$ with the following property.
Let $G$ be $(1,\mu,\nu)$-restricted, and such that 
$\chi(H)\le \tau_1$ for every induced subgraph $H$ of $G$ with $\omega(H)<\nu$.
If  $C\subseteq V(G)$ with $\chi(C)>d$,
then there is no strongly-independent 2-multicover
of $C$ in $G$ with length $m$.
\end{thm}
\Proof Choose $n$ to satisfy \ref{impression} taking $\lambda=1$; and choose $m,d\ge 0$ to satisfy \ref{usetick1}.
Now let $G$ be as in the theorem; then $G$ is $(1,\mu,\nu)$-restricted, and so by \ref{impression}, there is
no impression of $K_{n,n}$ in $G$ of order at most $2$. The result follows from \ref{usetick1}. This proves
\ref{betterusetick1}.~\bbox

Because of \ref{betterusetick1}, for our pervasiveness problem, we win if we can find a strongly-independent 2-multicover
in $G$ of sufficient length and covering a set $C$ with large enough chromatic number; and so several theorems to come will
have as a hypothesis that there is no such 2-multicover. For brevity, let us say $G$ is {\em $(m,c)$-limited}
if for every subset $C\subseteq V(G)$ with $\chi(G)>c$, there is no strongly-independent 2-multicover
of $C$ of length $m$ in $G$.

The next result is closely related to theorem 3.1 of~\cite{longholes}.
\begin{thm}\label{getindpt}
For all $m\ge 0$ and $\xi\ge 1$, there exist $\zeta\ge 0$ such that for all $c\ge 0$,
every $(m,c)$-limited graph is $(\xi,\zeta,c)$-free.
\end{thm}
\Proof
Choose an integer $\zeta\ge 0$ such that
for every partition of the edges of $K_{\zeta}$ into $\xi$ classes, some $K_{m}$ subgraph has all its edges in the same class.
We claim that $\zeta$ satisfies the theorem. For let $G$ be a graph that is not
$(\xi,\zeta,c)$-free. Consequently for some $C\subseteq V(G)$ with $\chi(C)>c$, there is an independent $\xi$-clique-multicover
of $C$ with length $\zeta$, say $(\mathcal{L}_i:i\in I)$ where $|I|=\zeta$. For each $i\in I$, let $\mathcal{L}_i = (X_i, N_i)$,
and take an enumeration of $X_i$. Thus we may speak of the $p$th vertex of $X_i$ for $1\le p\le \xi$. For each $i$,
let $N_i'\subseteq N_i$ be the set of vertices in $N_i$ with a neighbour in $C$. For each pair
$i,j\in I$ with $i<j$, choose $p\in \{1\ll \xi\}$ such that the $p$th vertex of $X_j$ has no neighbours in $N_i'$
(this is possible since $(\mathcal{L}_i:i\in I)$ is independent); we call $p$ the {\em colour} of the pair $(i,j)$.
From the choice of $\zeta$, there exists $I'\subseteq I$
with $I'=m$ such that all pairs $(i,j)$ with $i,j\in I'$ and $i<j$ have the same colour, say $p$. For each $i\in I'$ let $x_i$
be the $p$th vertex of $X_i$; and let $\mathcal{L}'_i = (\{x_i\}, N_i')$. Then
$(\mathcal{L}_i':i\in I')$ is a strongly-independent 2-multicover
of $C$ in $G$ with length $m$; and so $G$ is not $(m,c)$-limited. This proves \ref{getindpt}.~\bbox

\section{Where are we going?}\label{sec:wherearewe}

It might be helpful at this stage if we try to sketch the difficulties that lie ahead and our route around them.
We have seen that we can assume we have a $\xi$-clique-multicover of huge length, covering some set $C$ with huge chromatic number.
Any subsequence is also a $\xi$-clique-multicover, and because of \ref{getindpt}, there is no long independent subsequence.
This is asking for us to apply Ramsey's theorem, and obtain a long sequence where each pair of terms are the ``opposite''
of independent, but what does that mean? Just ``not independent'' does not tell us anything worthwhile. Before we apply
Ramsey's theorem, it is better
to tidy up each pair of terms first, shrinking them as necessary, to make them either independent or ``very'' non-independent;
what can we arrange? 

If $(X_1, N_1)$ and $(X_2, N_2)$ are terms (in this order) of the $\xi$-clique-multicover of $C$, we would like to 
arrange 
that some vertex in $X_2$ has no neighbour in the set of vertices in $N_1$ that have neighbours in $C$; 
and it would be enough to arrange that no vertex in $N_1$ is complete
to $X_2$ (because then, since $X_2$ has bounded size, some vertex in $X_2$ would be nonadjacent to a big subset of $N_1$,
big enough to cover a large chromatic number part of $C$, and we could throw away the rest). So the problem is, vertices in $N_1$
that are complete to $X_2$. If the set of vertices in $N_1$ that are not complete to $X_2$ covers a big-$\chi$ part of $C$, we
could just take that, and delete the remainder of $N_1$; and if not then the vertices in $N_1$ that are complete to $X_2$
cover a big-$\chi$ part of $C$, so we could just take that. That would be one way to tidy up the pair; we would obtain
a pair that is either independent, or has the property that every vertex in $N_1$ is complete to $X_2$. We tidy up
every pair in this way, and then we apply Ramsey;
one outcome is a long sequence of $\xi$-clique-covers, pairwise independent, which is impossible; and the other is a long
sequence of $\xi$-clique-covers where the base of each is complete to the clique of every later term.
This unfortunately does not work; the second outcome is not rich enough to be useful. We have to tidy up the pairs more carefully.

When our sequence of $\xi$-clique-covers was created in the first place, we first chose one, say $(X_1, N_1)$, covering $C_1$;
then we chose $(X_2, N_2)$ covering $C_2$ in $G[C_1]$, and so on. In particular, every vertex of every later $X_j\cup N_j$ has a neighbour in every $N_i$.
So far we have used the fact that every vertex in the ultimate set $C$ has a neighbour in each $N_i$, 
and have been resigned to the fact that
vertices in $X_j\cup N_j$ might have neighbours in earlier $N_i$'s; but in fact they do have such neighbours, and these 
edges are useful and need to be carefully guarded, particularly in the case when we fail to get a long independent subsequence.
Here is a better way to tidy up the pairs, that is not so cavalier about the edges between $N_i$ and $N_j$.
(But it doesn't seem to work if we start with a sequence and try to tidy it; it only works if we grow the sequence term-by-term and tidy as we go.)

Again, start with $(X_1, N_1)$, covering $C_1$ say. For a clique $X\subseteq C_1$ (or a single vertex $X\in C_1$)
let us say the ``up-down-$\chi$'' of $X$ is the 
chromatic number of the set of vertices in $C_1$ that have a neighbour in $N_1$ that is complete to $X$.
Partition $C_1$ into two sets, one the union of all $\xi$-cliques with big up-down-$\chi$, and the other
its complement. One of them has big $\chi$, so we work inside that. 

Here there is a problem;  when we remove some of $C_1$, the up-down-$\chi$ of the $\xi$-cliques we keep might drop. So, perhaps we have
a subset of $C_1$ with big $\chi$, a union of $\xi$-cliques that all used to have big up-down-$\chi$. To make use of this
property, we need to keep track of the old $C_1$. As we grow more terms in the clique-multicover there will be more ``old'' sets 
that we need to keep track of, and we assemble them in a sequence called a ``world''.
Anyway, let us ignore the world for this sketch.

Choose a $\xi$-clique-cover
$(X_2, N_2)$ of $C_2$ say, all in $G[C_1]$, and let $Y$ be the set of vertices in $N_1$ complete to $X_2$.
The vertices in $C_2$ all have neighbours in $N_1$. If many (in the big-$\chi$ sense) have a neighbour
in $N_1\setminus Y$, we can tidy to make an independent pair of $\xi$-clique-covers by deleting the other part of $C_2$,
and we rejoice; so either that,
or by throwing away a small part of $C_2$, we can arrange that $C_2$ is anticomplete to $N_1\setminus Y$, and the $\xi$-clique
$X_2$ has big up-down-$\chi$ through $N_1$. Hence the vertices
in $N_2$ also belong to $\xi$-cliques that used to have big up-down-$\chi$, because of the way we partitioned $C_1$.
But each vertex $v$ in $N_2$ only
had small up-down-$\chi$ via $Y$, because any vertex that $v$ could reach in two steps via $Y$ belongs
to $N^2(X_2\cup \{v\})$, and the clique $X_2\cup \{v\}$ is too large to have second neighbours with big $\chi$. 
(This step is the primary reason why we are looking at $\xi$-clique-covers with $\xi$ maximum instead of 1-clique-covers.)
So $v$ had a neighbour in $N_1\setminus Y$, and therefore still has such a neighbour (we discarded part of $C_1$ but did not change
$N_1$). This is still the argument we used in~\cite{longholes}, but now comes a refinement;
$v$ has {\em many} neighbours in $N_1\setminus Y$, enough that it used to have big up-down-$\chi$ via these neighbours. This is a key
observation. 
The two possible
outcomes are, therefore, that either we obtain an independent pair, or we obtain a pair $(X_1, N_1), (X_2, N_2)$
where every vertex in $N_2$ belongs to a $\xi$-clique with big up-down-$\chi$ via $N_1\setminus Y$ (with notation as before) and 
some extra set $W_2$ 
(that was the old $C_1$ before we discarded some of it), and $C_2$
is anticomplete to $N_2\setminus Y$. We call this a ``$\beta$-skew'' pair ($\beta$ measures the size of the up-down-$\chi$,
but in this sketch we ignore $\beta$, and just call it a skew pair.)
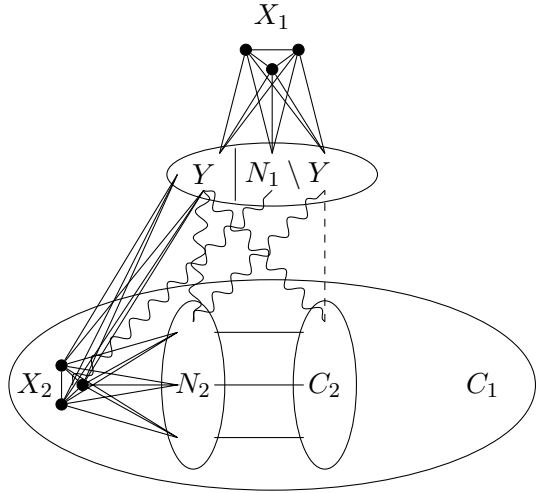
\begin{figure}[h]
\centering
\tikzset{snake it/.style={decorate, decoration=snake}}
\begin{tikzpicture}[scale=.7,auto=left]
\node at (0,1) {$X_1$};
\node at (-1.3,-2) {$Y$};
\node at (0.3,-2) {$N_1\setminus Y$};
\node at (4,-6) {$C_1$};
\node at (1,-6) {$C_2$};
\node at (-1.5,-6) {$N_2$};
\node at (-4.5,-6) {$X_2$};
\path (0,-1.6) coordinate (t1);
\path (1,-1.6) coordinate (t2);
\path (-1,-1.6) coordinate (t3);
\path (-1.8,-5) coordinate (s1);
\path (-1.8,-6) coordinate (s2);
\path (-1.8,-7) coordinate (s3);
\path (-1.8,-2.0) coordinate (r1);
\path (-1.3,-2.3) coordinate (r2);
\tikzstyle{every node}=[inner sep=1.5pt, fill=black,circle,draw]
\node (x1) at (0,0) {};
\node (x2) at (0.5,0.37) {};
\node (x3) at (-0.5, 0.37) {};

\node (y1) at (-3.6,-6) {};
\node (y2) at (-4,-6.37) {};
\node (y3) at (-4,-5.63) {};
\draw (0,-2) ellipse (2cm and .6cm);
\draw (0,-6) ellipse (5cm and 2cm);
\draw (-1.5,-6) ellipse (0.6cm and 1.6cm);
\draw (1,-6) ellipse (0.6cm and 1.6cm);
\foreach \i in {1,2,3}
{\foreach \j in {1,2}
\draw (y\i)--(r\j);}
\draw (x1) -- (t1);
\draw (x1) -- (t2);
\draw (x1) -- (t3);
\draw (x3) -- (t1);
\draw (x3) -- (t2);
\draw (x3) -- (t3);
\draw (x2) -- (t1);
\draw (x2) -- (t2);
\draw (x2) -- (t3);
\draw (x1) -- (x2);
\draw (x1) -- (x3);
\draw (x2) -- (x3);
\draw (-0.7,-1.5) -- (-0.7,-2.5);
\draw (y1) -- (y2);
\draw (y1) -- (y3);
\draw (y2) -- (y3);
\draw (y1)-- (s1)--(y2)--(s2)--(y3)--(s3)--(y1)--(s2);
\draw (y2) --(s3);
\draw (y3) -- (s1);

\foreach \i in {-5,-6,-7}
\draw (-1.1,\i) -- (0.6,\i);

\path[draw=black, snake it] (-1.3, -2.3) -- (1, -4.8);
\path[draw=black, snake it] (-1.3,-2.3) -- (-1.5,-4.8);
\path[draw=black, snake it] (-1.5,-4.8) -- (1,-2.3);
\path[draw=black, dashed] (1,-4.8) -- (1, -2.3);
\path[draw=black, snake it] (-3.8,-6) -- (0,-2.3);

\end{tikzpicture}
\caption{Birth of a skew pair (dashed $=$ anticomplete).} \label{fig:6}
\end{figure}

Now we go on to the birth of the third pair $(X_3, N_3)$, chosen within $G[C_2]$. We have to tidy up both the pairs 
$(X_1, N_1), (X_3, N_3)$ and $(X_2, N_2), (X_3, N_3)$, in the same way. One problem is, this might mess up what we already
did. For instance, perhaps we have arranged the pair $(X_1,N_1), (X_2,N_2)$ to be skew, and the pair $(X_1,N_1), (X_3,N_3)$
wants to be independent, and we therefore have to shrink $N_1$ to make this so. There is a danger that shrinking $N_1$ will mess up
the fact that $N_2$ is a union of $\xi$-cliques with big up-down-$\chi$ via $N_1\setminus Y$ (with notation as before). But we will be careful that
the vertices we remove from $N_1$ all have neighbours in $C_3$, and the vertices in $N_1\setminus Y$ do not.

So the third pair can be tidied, and so on; eventually we get a long sequence of $\xi$-clique covers of some set $C$, such that each pair is
either independent or skew. Now we apply Ramsey; and get a long subsequence such that all pairs are independent, or all pairs
are skew. The first is impossible, as always, so we have built a long sequence of $\xi$-clique-covers, all pairwise skew. 

This is an interesting object. We can show it contains any chandelier, and indeed any lamp, as an induced subgraph; it is much richer than
the thing we had before. One can greedily embed a tree into it; first embed the root at some vertex $v_k$ of some $N_k$
with $k$ large. Next we embed the neighbours of the root. There are vertices in each earlier $N_j$ that are adjacent to $v_k$; 
so choose one such vertex from $N_{k-1}$, one from $N_{k-2}$
and so on until we have enough.
We have to make these pairwise nonadjacent;
and this is where we use the key observation from above, that $v_k$ has many neighbours in $N_j$, enough
that it used to have big second neighbours via these neighbours, and we can argue that there is always one nonadjacent 
to all the vertices we have already chosen 
(except $v_k$).
Now start filling in the second neighbours of $v_k$ in the tree, and so on. To get a chandelier, arrange that each leaf of the 
tree is chosen from $N_1$; and then we can use a vertex from $X_1$ as the pivot. Lamps can be embedded the same way.

Unfortunately, this is not yet good enough: we don't want lamps, we want trees of lamps. How can we modify this to get a 
tree of lamps?
(Or tree of chandeliers, say, for this sketch -- though this method does not quite get every tree of chandeliers.) 
Notice that the pivot in the chandelier we just built 
could be chosen to be any vertex of $X_1$; so whenever we find a $\xi$-clique-cover $(X_1, N_1)$ of some set $C$ and we can extend
it to a long sequence of pairwise skew $\xi$-clique-covers, we can get a chandelier with pivot in $X_1$. And the definition
of ``big up-down-$\chi$'' ensures
that when we embed the chandelier, all the vertices we use belong to $\xi$-cliques $X$ such that there is a 
$\xi$-clique-cover $(X,N)$ of some ``semi-private'' big-$\chi$ set in which we can try to grow any desired pendant tree of lamps without
too much interruption from other vertices (again, this is a place where the world intrudes; and not true for the leaves
of the tree, embedded in $N_1$, which explains the curious composition rule for trees of lamps, and explains why we cannot get every
tree of chandeliers). 

So our problem is, we have
a $\xi$-clique-cover $(X,N)$ covering a set $C$ with big $\chi$, and we would be happy if we could prove that it can be extended
to a long sequence of pairwise skew $\xi$-clique-covers. Certainly it can be extended to a long sequence of 
$\xi$-clique-covers, and we can tidy them and then apply Ramsey; but the long skew subsequence we get 
might no longer include the first term. We have to do something so that we can get the long skew sequence without 
discarding the first term.

Can we always get a skew sequence of length two with specified first term? If we could, then look at the set they cover
in common, and do it again, tidying up all the pairs as we go; we would generate a long sequence of $\xi$-clique-covers, 
still including the given first term,
such that the first term and $i$th term are skew, for all~$i$. Then apply Ramsey to the sequence with first term removed,
get a long skew subsequence, and put the first term back, and we have won. So, the problem is just getting a skew
sequence of length two with a specific first term. 

Say a $\xi$-clique-cover, covering a set of large chromatic number, is ``bad'' if we cannot extend it (or a truncation of it)
to a skew sequence of length two, still covering
a set of large chromatic number. If we can move to a subideal, still with unbounded chromatic number, in which there are no bad
 $\xi$-clique-covers, do that. If not, then in some sense there are bad clique covers
everywhere; take a long sequence of them, and clean it up, and it turns into a long independent sequence, which is impossible.
This is the idea of the main proof of the next section.

\section{Skew pairs}

If $Z, W\subseteq V(G)$ are disjoint and $\beta\ge 0$ and $\xi>0$, we say that a clique $X\subseteq W$ is {\em $\beta$-earthed via $(Z,W)$}
if $\chi(M)>\beta$,
where $M$ is the set of all vertices in $W\setminus X$ that are anticomplete to $X$ and have a neighbour in $Z$ that is complete to $X$.
We say a vertex $v\in W$ is
{\em $(\beta,\xi)$-earthed via $(Z, W)$}
if there is a $\xi$-clique $X\subseteq W$ with $v\in X$, such that $X$ is $\beta$-earthed via $(Z,W)$.
(This is more-or-less the concept we called
``big up-down-$\chi$'' in section \ref{sec:wherearewe}.) We observe that
if $Z\subseteq Z'\subseteq V(G)\setminus W$, then every vertex that is $(\beta,\xi)$-earthed via $(Z,W)$ is also
$(\beta,\xi)$-earthed via $(Z',W)$.

Let $\mathcal{M} = (\mathcal{L}_i:i\in I)$ be a $\xi$-clique-multicover of $C$ in $G$.
A {\em world} for $\mathcal{M}, C$ is a family $\mathcal{W}=(W_i:i\in I)$ of subsets of $V(G)$ such that for all $i,j\in I$:
\begin{itemize}
\item if $i\le j$ then $W_i\supseteq W_j\supseteq C$;
\item if $i < j$ then $V(\mathcal{L}_i)\cap W_j=\emptyset$, and if $i\ge  j$ then $V(\mathcal{L}_i)\subseteq W_j$;
\item if $i <  j$ then $X(\mathcal{L}_i)$ is anticomplete to $W_j$.
\end{itemize}

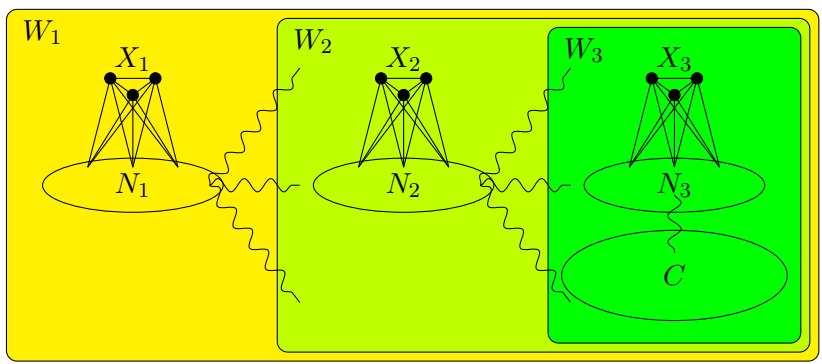
\begin{figure}[H]
\centering

\tikzset{snake it/.style={decorate, decoration=snake}}
\begin{tikzpicture}[scale=.6,auto=left]
\draw [rounded corners] (-8.8,1.9) rectangle (9.2,-5.9);
\draw [fill=lightgray, rounded corners] (-2.8,1.7) rectangle (9.0,-5.7);
\draw [fill=gray, rounded corners] (3.2,1.5) rectangle (8.8,-5.5);
\node at (6,-4) {$C$};
\node at (-8,1.4) {$W_1$};
\node at (-2,1.2) {$W_2$};
\node at (4,1.0) {$W_3$};
\node at (-6,-2) {$N_1$};
\node at (-0,-2) {$N_2$};
\node at (6,-2) {$N_3$};
\node at (-6,0.8) {$X_1$};
\node at (-0,0.8) {$X_2$};
\node at (6,0.8) {$X_3$};

\tikzstyle{every node}=[inner sep=1.5pt, fill=black,circle,draw]
\node (x1) at (-6,0) {};
\node (x2) at (-6.5,0.37) {};
\node (x3) at (-5.5, 0.37) {};
\node (y1) at (0,0) {};
\node (y2) at (0.5,0.37) {};
\node (y3) at (-0.5, 0.37) {};
\node (z1) at (6,0) {};
\node (z2) at (6.5,0.37) {};
\node (z3) at (5.5, 0.37) {};

\draw (-6,-2) ellipse (2cm and .6cm);
\draw (0,-2) ellipse (2cm and .6cm);
\draw (6,-2) ellipse (2cm and .6cm);
\draw (6,-4) ellipse (2.5cm and 1cm);
\draw (-6,0) -- (-6,-1.6);
\draw (-6,0) -- (-5,-1.6);
\draw (-6,0) -- (-7,-1.6);
\draw (-6.5,0.37) -- (-6,-1.6);
\draw (-6.5,0.37) -- (-5,-1.6);
\draw (-6.5,0.37) -- (-7,-1.6);
\draw (-5.5, 0.37) -- (-6,-1.6);
\draw (-5.5, 0.37) -- (-5,-1.6);
\draw (-5.5, 0.37) -- (-7,-1.6);
\draw (-6,0) -- (-6.5,0.37);
\draw (-6,0) -- (-5.5, 0.37);
\draw (-6.5,0.37) -- (-5.5, 0.37);

\draw (0,0) -- (0,-1.6);
\draw (0,0) -- (-1,-1.6);
\draw (0,0) -- (1,-1.6);
\draw (0.5,0.37) -- (0,-1.6);
\draw (0.5,0.37) -- (-1,-1.6);
\draw (0.5,0.37) -- (1,-1.6);
\draw (-0.5, 0.37) -- (0,-1.6);
\draw (-0.5, 0.37) -- (-1,-1.6);
\draw (-0.5, 0.37) -- (1,-1.6);
\draw (0,0) -- (0.5,0.37);
\draw (0,0) -- (-0.5, 0.37);
\draw (0.5,0.37) -- (-0.5, 0.37);

\draw (6,0) -- (6,-1.6);
\draw (6,0) -- (5,-1.6);
\draw (6,0) -- (7,-1.6);
\draw (6.5,0.37) -- (6,-1.6);
\draw (6.5,0.37) -- (5,-1.6);
\draw (6.5,0.37) -- (7,-1.6);
\draw (5.5, 0.37) -- (6,-1.6);
\draw (5.5, 0.37) -- (5,-1.6);
\draw (5.5, 0.37) -- (7,-1.6);
\draw (6,0) -- (6.5,0.37);
\draw (6,0) -- (5.5, 0.37);
\draw (6.5,0.37) -- (5.5, 0.37);

\path[draw=black, snake it] (-4.3, -2) -- (-2.3, 0.6);
\path[draw=black, snake it] (-4.3, -2) -- (-2.3,-2);
\path[draw=black, snake it] (-4.3, -2) -- (-2.3,-4.6);
\path[draw=black, snake it] (1.7, -2) -- (3.7, 0.6);
\path[draw=black, snake it] (1.7, -2) -- (3.7,-2);
\path[draw=black, snake it] (1.7, -2) -- (3.7,-4.6);
\path[draw=black, snake it] (6, -2.2) -- (6,-3.5);

\end{tikzpicture}

\caption{A world for a clique-multicover} \label{fig:5}
\end{figure}

Let $\mathcal{M} = (\mathcal{L}_i:i\in I)$ be a $\xi$-clique-multicover of $C$ in $G$, where $\mathcal{L}_i=(X_i, N_i)$ for
each $i\in I$, and let $\mathcal{W}=(W_i:i\in I)$ be a world for $\mathcal{M}, C$.
Let $i,j\in I$ with $i<j$, and let $Z$ be the set of vertices in $N_i$ that are not complete to $X_j$;
we say that the pair $(\mathcal{L}_i,\mathcal{L}_j)$ is
\begin{itemize}
\item {\em skew with respect to $\mathcal{M}, C,\mathcal{W}$} if $Z$ is
anticomplete to $C$ and to $W_k$ for all $k\in I$ with $k>j$;
\item {\em $\beta$-skew with respect to
$\mathcal{M},C,\mathcal{W}$} if it is skew with respect to
$\mathcal{M},C,\mathcal{W}$, and every vertex in $N_j$ is $(\beta,\xi)$-earthed via $(Z, W_{j})$.
\end{itemize}
We say that $\mathcal{M}$ is {\em skew with respect to $C,\mathcal{W}$} if all its pairs are skew
with respect to $\mathcal{M},C,\mathcal{W}$; and
similarly define {\em $\beta$-skew with respect to $C,\mathcal{W}$}
if all its pairs have the corresponding property.

Let $(X,N)$ be a $\xi$-clique-cover of $C$, and let $N'\subseteq N$. We call $(X,N')$ a {\em truncation} of $(X,N)$.
Let $\mathcal{M} = (\mathcal{L}_i:i\in I)$ be a $\xi$-clique-multicover of $C$, and for each $i\in I$ 
let $\mathcal{L}_i'$ be a truncation of $\mathcal{L}_i$. Then we say $(\mathcal{L}_i':i\in I)$ is a {\em truncation} of $\mathcal{M}$.

If $\mathcal{M}= (\mathcal{L}_i:i\in I)$ is a $\xi$-clique-multicover of $C$ in $G$, and $\mathcal{W}$ is a world for $\mathcal{M},C$,
a pair $(\mathcal{L}_i,\mathcal{L}_j)$ is
{\em $\beta$-tidy} with respect to $\mathcal{M}, C, \mathcal{W}$
if it is either independent with respect to $C$ or $\beta$-skew with respect to $\mathcal{M},C,\mathcal{W}$.
If every pair in $\mathcal{M}$ is $\beta$-tidy with respect to $\mathcal{M}, C, \mathcal{W}$, we say that $\mathcal{M}$ is
{\em $\beta$-tidy} with respect to $C,\mathcal{W}$.

It would be convenient if, given a $\xi$-clique-multicover $\mathcal{M}= (\mathcal{L}_i:i\in I)$ of $C$ in $G$, there is a $\beta$-tidy
truncation of $\mathcal{M}$ of the same length, covering some $C'\subseteq C$  where $\chi(C')$ is 
large (if we begin with $\chi(C)$ large enough). Unfortunately, this is false, even for multicovers of length two, and 
we need to work around this difficulty. It is true that, given a $\xi$-clique-multicover $\mathcal{M}$ of $C$, that is already
$\beta$-tidy, we can replace it by a truncation of the same length, and add
another term to the end, chosen with vertex set within $C$, and make a longer $\xi$-clique-multicover that is
still $\beta$-tidy; but this is not quite enough for what we need. We need to add a new last term in such a way that
the pair it makes with the first term is not just $\beta$-tidy but $\beta$-skew, and the following will help us 
to do that.

\begin{thm}\label{fixlast}
Let $\xi,t \ge 1$, and $\beta, \tau_2\ge 0$ and $c\ge \tau_2$.
Let $G$ be a graph such that
$\chi(N^2_G(X))\le \tau_2$ for every $(\xi+1)$-clique $X$ in $G$.
Let $\mathcal{M}= (\mathcal{L}_i:i\in I)$ be a $\xi$-clique-multicover of $C$ in $G$ with length nonzero and at most $t$, 
where $\chi(C)>(c+\beta)(\xi+1)^t$,
and let $\mathcal{W}=(W_i:i\in I)$ be a world for $\mathcal{M}, C$. Let $k\in I$ be maximum.
Let $(\mathcal{L}_i:i\in I\setminus \{k\})$ be $\beta$-tidy with respect to $W_k,(W_i:i\in I\setminus \{k\})$.
Suppose that for each $i\in I\setminus\{k\}$,
either
\begin{itemize}
\item the pair $(\mathcal{L}_i, \mathcal{L}_k)$ is $\beta$-tidy with respect to $\mathcal{M},C,\mathcal{W}$; or
\item every vertex in $X_k\cup N_k$ is $(\beta+\tau_2,\xi)$-earthed via $(N_i, W_k)$; or
\item no vertex in $X_k\cup N_k$ is $(\beta+\tau_2,\xi)$-earthed via $(N_i, W_k)$.
\end{itemize}
Then there exists $C'\subseteq C$ with $\chi(C')>c$, and a truncation $\mathcal{M}'=(\mathcal{L}_i':i\in I)$ of $\mathcal{M}$ covering $C'$,
such that
\begin{itemize}
\item $\mathcal{M}'$ is $\beta$-tidy with respect to $C',\mathcal{W}$; 
\item for all $i<j\in I$, if $(\mathcal{L}_i, \mathcal{L}_j)$
is independent with respect to $C$ then $(\mathcal{L}_i', \mathcal{L}_j')$
is independent with respect to $C'$; and
\item for all $i<j\in I$, if $(\mathcal{L}_i, \mathcal{L}_j)$
is $\beta$-skew with respect to $\mathcal{M},C,\mathcal{W}$, then $(\mathcal{L}_i', \mathcal{L}_j')$
is $\beta$-skew with respect to $\mathcal{M}',C',\mathcal{W}$.
\end{itemize}
\end{thm}
\Proof We are given that $t\ge |I|$; but for inductive purposes, let us weaken this hypothesis, and just assume that
$t$ is at least the number of $i\in I\setminus \{k\}$ such that
the pair $(\mathcal{L}_i, \mathcal{L}_k)$ is not $\beta$-tidy with respect to $\mathcal{M},C,\mathcal{W}$. We will prove the same conclusion.

For each $i\in I$ let $\mathcal{L}_i = (X_i, N_i)$.
We may assume that there exists $h\in I\setminus \{k\}$ such that
the pair $(\mathcal{L}_h, \mathcal{L}_k)$ is not
$\beta$-tidy with respect to $\mathcal{M},C,\mathcal{W}$, for if not then the result is true. 
Let $X_k=\{x_1\ll x_{\xi}\}$, for $1\le s\le \xi$
let $Y_s$ be the set of vertices in $N_h$
that are nonadjacent to $x_s$, and let $C_{s}$ be the set of vertices in $C$ that have a neighbour in $Y_s$. 
Let $C_{0}$ be the set of vertices in $C$ with no neighbour in $Y_1\cup\cdots\cup Y_{\xi}$.

One of $C_0\ll C_{\xi}$ has chromatic
number more than $(c+\beta)(\xi+1)^{t-1}$ say $C_s$. If $s>0$, define $C'=C_s$ and $N_h'=Y_s$; and otherwise define $C'=C_0$ and
$N_h'=N_h$.
In either case define $N_i'=N_i$ for $i\in I\setminus \{h\}$;
and for each $i\in I$, let $\mathcal{L}_i'=(X_i, N_i')$.
Let $\mathcal{M}'=(\mathcal{L}_i':i\in I)$. 
\\
\\
(1) {\em For $1\le i<j\le k$, if $(\mathcal{L}_i, \mathcal{L}_j)$ is 
independent with respect to $C$ then $(\mathcal{L}_i', \mathcal{L}_j')$
is independent with respect to $C'$, and if $(\mathcal{L}_i, \mathcal{L}_j)$
is $\beta$-skew with respect to $\mathcal{M},C,\mathcal{W}$, then $(\mathcal{L}_i', \mathcal{L}_j')$
is $\beta$-skew with respect to $\mathcal{M}',C',\mathcal{W}$.}
\\
\\
Suppose that $(\mathcal{L}_i, \mathcal{L}_j)$
is independent with respect to $C$. Then
there exists $x\in X_j$ such that no vertex in $N_i$ with a neighbour in $C$ is adjacent to $x$.
Consequently no vertex in $N_i'$ with a neighbour in $C'$ is adjacent to $x$, and so 
$(\mathcal{L}_i', \mathcal{L}_j')$
is independent with respect to $C'$.

Now suppose that $(\mathcal{L}_i, \mathcal{L}_j)$
is $\beta$-skew with respect to $\mathcal{M},C,\mathcal{W}$. 
Let $Z_{i,j}$ be the set of vertices in $N_i$ that are not complete to $X_j$;
then $Z_{i,j}$ is
anticomplete to $C$ and to $W_{j+1}\ll W_k$ 
and every vertex in $N_j$ is $(\beta,\xi)$-earthed via $(Z_{i,j}, W_{j})$. If $Z_{i,j}\subseteq N_i'$, then
$Z_{i,j}$ is the set of vertices in $N_i'$ that are not complete to $X_j$; and $Z_{i,j}$ is
anticomplete to $C'$ and to $W_{j+1}\ll W_k$; and every vertex in $N_j'$ is $(\beta,\xi)$-earthed via $(Z_{i,j}, W_{j})$, since
$N_j'\subseteq N_j$;
and so $(\mathcal{L}_i', \mathcal{L}_j')$
is $\beta$-skew with respect to $\mathcal{M}',C',\mathcal{W}$. Thus we may assume (for a contradiction) 
that $Z_{i,j}\not\subseteq N_i'$, and
consequently $i=h$ and $C'=C_s$ for some $s\in \{1\ll \xi\}$. 
Since $(\mathcal{L}_h, \mathcal{L}_k)$ is not
$\beta$-tidy with respect to $\mathcal{M},C,\mathcal{W}$, and  $(\mathcal{L}_i, \mathcal{L}_j)$
is $\beta$-skew with respect to $\mathcal{M},C,\mathcal{W}$, and $h=i$, it follows that $j\ne k$ and so $j<k$;
and therefore
$Z_{i,j}$ is anticomplete to $W_k$.  
Let $v\in Z_{i,j}\setminus N_i'$. Since $v\notin N_i'=Y_s$, it follows that $v$ is adjacent to $x_s$; but $x_s\in W_k$, and 
$Z_{i,j}$ is anticomplete to $W_k$, a contradiction.
This proves (1).

\bigskip

If $C'=C_s$ where $s>0$, then 
the pair $(\mathcal{L}_h', \mathcal{L}_k')$ is independent with respect to $C'$, and so from (1) and the inductive hypothesis
applied to $\mathcal{M}'$ and $C'$, the result follows. We may therefore assume that $C'=C_0$, and so $\mathcal{M}'=\mathcal{M}$.
We claim that the pair $(\mathcal{L}_h', \mathcal{L}_k')$ is $\beta$-skew with respect to $\mathcal{M},C',\mathcal{W}$. 
Let  $Z_{h,k}$ be the union of the sets $Y_1\ll Y_{\xi}$, that is, the set of vertices in $N_h$ with a nonneighbour in $X_k$.
We must check that:
\begin{itemize}
\item $Z_{h,k}$ is
anticomplete to $C'$; and
\item every vertex in $N_k$ is $(\beta,\xi)$-earthed via $(Z_{h,k}, W_{k})$.
\end{itemize}
The first claim follows from the definition of $C_0$, since $C'=C_0$. For the second, let $v\in N_k$. Now
every vertex in $C_0$ has a neighbour in $N_h$, and has no neighbour in $Z_{h,k}$; and so it
has a neighbour in $N_h\setminus Z_{h,k}$, and this neighbour is complete to $X_k$. 
But $\chi(C_0)>(c+\beta)(\xi+1)^{t-1}\ge \beta+\tau_2$, and since $C_0\subseteq W_k$, it follows that every vertex in $X_k$ is
$(\beta+\tau_2,\xi)$-earthed via $(N_h,W_k)$. Therefore, from the hypothesis, every vertex in $X_k\cup N_k$ is 
$(\beta+\tau_2,\xi)$-earthed via $(N_h, W_k)$, and in particular this is true for $v$. 
Let $X\subseteq W_k$ be a $\xi$-clique containing $v$ that is $\beta$-earthed via $(N_h, W_k)$, and let $M$ 
be the set of vertices in $W_k$
that are anticomplete to $X$ and have a neighbour in $N_h$ that is complete to $X$; thus $\chi(M)>\beta+\tau_2$.
Let $D$ be the set of $u\in M$
such that $u$ is adjacent to some vertex in $N_h$ that is complete to $X_k\cup \{v\}$. Then $\chi(D)\le \tau_2$, from
the hypothesis; and so $\chi(M\setminus D)>\beta$. But every vertex in $M\setminus D$ has a neighbour in $N_h$
that is complete to $X$ and not complete to $X_k$; and so this neighbour belongs to $Z_{h,k}$. This proves that $v$
is $(\beta,\xi)$-earthed via $(Z_{h,k}, W_{k})$, as claimed; and so proves that
the pair $(\mathcal{L}_h', \mathcal{L}_k')$ is $\beta$-skew with respect to $\mathcal{M},C',\mathcal{W}$. Consequently
the result follows from the inductive hypothesis, applied to $\mathcal{M}$ and $C'$. This proves \ref{fixlast}.~\bbox

For the next result, let us fix $\xi$; the functions we are about to describe depend on $\xi$, but it is cumbersome
to keep mentioning it, particularly since $\xi$ is constant throughout.

Let $\mathcal{L}=(X_1,N_1)$ be a $\xi$-clique-cover in $G$ of $C$. For $c\ge 0$, we say that 
$(\mathcal{L},C)$ is {\em $c$-skewable (in $G$)}
if 
there exist $N_1'\subseteq N_1$, and $C'\subseteq C$ with $\chi(C')>c$, and a $\xi$-clique-cover $(X_2,N_2)$ of $C'$
with $X_2,N_2\subseteq C$, such that
\begin{itemize}
\item  $((X_1,N_1'), (X_2,N_2))$ is a $\xi$-clique-multicover of $C'$ (of length two); 
\item
$Z$ is anticomplete to $C'$, and 
every vertex in $N_2$ is $(c,\xi)$-earthed via $(Z, C)$, where $Z$ is the set of vertices in $N_1'$ that are not complete to $X_2$.
\end{itemize}
(In other words, the $\xi$-clique-multicover $((X_1,N_1'), (X_2,N_2))$ of $C'$ is $c$-skew with respect to $C'$ and the world 
$(V(G), C)$.)

Let $\phi:\mathbb{N}\rightarrow \mathbb{N}$ be non-decreasing.
We say a graph $G$ is {\em $\phi$-skewable} if
for all $c\ge 0$, every $C\subseteq V(G)$ with $\chi(C)>\phi(c)$, and every $\xi$-clique-cover
$\mathcal{L}$ of $C$ in $G$, $(\mathcal{L},C)$ is $c$-skewable. An ideal $\mathcal{C}$ of graphs is {\em skewable} if there is a 
non-decreasing function $\phi$ such that every 
$G\in \mathcal{C}$ is $\phi$-skewable. As we said, all these definitions depend on $\xi$,
and to emphasize that we sometimes say ``$\phi$-skewable relative to $\xi$'', and similar expressions.

\begin{thm}\label{skewpair}
Let $\xi\ge 1$, and $\zeta,\tau_2, \tau_3\ge 0$. Let $\mathcal{C}$ be a non-colourable ideal of graphs such that 
for every $G\in \mathcal{C}$:
\begin{itemize}
\item $\chi(N^2_G(X))\le \tau_2$ for every $(\xi+1)$-clique $X$ in $G$; and
\item $G$ is $(\xi,\zeta,\tau_3)$-free.
\end{itemize}
Then there is a non-colourable subideal $\mathcal{C}'$ of $\mathcal{C}$ such that $\mathcal{C}'$ is skewable.
\end{thm}
\Proof
Let $\mathcal{C}'$ be a non-colourable subideal of $\mathcal{C}$. 
Let $K(\mathcal{C}')$ be the set of all integers $k\ge 0$ such that there exists $c\ge 0$ with the following property:
\begin{itemize}
\item 
For all $d\ge 0$, there exists $G\in \mathcal{C}'$, and $C\subseteq V(G)$ with $\chi(C)>d$, and
an independent $\xi$-clique-multicover $(\mathcal{L}_i:i\in I)$ of $C$ in $G$, where $|I|=k$,
such that for each $i\in I$, $(\mathcal{L}_i, C)$ is not $c$-skewable. 
\end{itemize}
Since $\mathcal{C}'$ is non-colourable,
$0\in K(\mathcal{C}')$; and since every graph $G\in \mathcal{C}$
is $(\xi,\zeta,\tau_3)$-free, $k<\zeta$ for all $k\in K(\mathcal{C}')$. 
Hence there is a largest number $k\in K(\mathcal{C}')$, and we call
$k$ the {\em rank} of $\mathcal{C}'$.
The rank of $\mathcal{C}'$ is zero if and only if $\mathcal{C}'$ is skewable.

Choose a non-colourable subideal $\mathcal{C}'$ of $\mathcal{C}$ with minimum rank; we claim that it satisfies the theorem, that is,
that its rank is zero. Suppose it has positive rank $k$ say. 
\\
\\
(1) {\em There exist $c_2\ge c_1\ge 0$ such that for every $c'\ge  0$, there exist a graph $G\in \mathcal{C}'$,
and subsets $C\subseteq D\subseteq V(G)$ with $\chi(C)>c'$, and a $\xi$-clique-multicover $(\mathcal{L}_i:1\le i\le k+1)$ of $C$ in $G$,
where $\mathcal{L}_i=(X_i, N_i)$ for $1\le i\le k+1$, 
such that:
\begin{itemize}
\item $(\mathcal{L}_i:1\le i\le k)$ is an independent $\xi$-clique-multicover of $D$;
\item $V(\mathcal{L}_{k+1})\subseteq D$;
\item for $1\le i\le k$, either every vertex in $V(\mathcal{L}_{k+1})$ is $(c_1+\tau_2,\xi)$-earthed via $(N_i,D)$, or none are;
\item for each $i\in \{1\ll k\}$, $(\mathcal{L}_i, D)$ is not $c_1$-skewable, and $(\mathcal{L}_{k+1}, C)$ 
is not $c_2$-skewable.
\end{itemize}}
\noindent Since $\mathcal{C}'$ has rank $k$, there exists $c_1\ge 0$ such that for each $d\ge 0$, we can 
choose $G_d\in \mathcal{C}'$, and $D\subseteq V(G_d)$
with $\chi(D)>2^kd$, such that there is an independent $\xi$-clique-multicover $(\mathcal{L}_i:1\le i\le k)$ of $D$ in $G_d$, 
and for $1\le i\le k$, $(\mathcal{L}_i, D)$ is not $c_1$-skewable in $G_d$. Let $\mathcal{L}_i = (X_i, N_i)$ for $1\le i\le k$.
For
$1\le i\le k$, there is a partition of $D$ into two parts,  where one of the parts consists of all vertices in $D$
that are $(c_1+\tau_2,\xi)$-earthed via $(N_i,D)$. Hence there is a partition of $D$ into $2^k$ parts, such that for
each part $B$, and for $1\le i\le k$, either every vertex in $B$ is $(c_1+\tau_2,\xi)$-earthed via $(N_i,D)$, or none are.
Since there are only $2^k$ parts, one of them, $D_d$ say, has chromatic number more than $d$. 

Let $\mathcal{C}''$ be the minimal subideal
of $\mathcal{C}$ that contains all the graphs 
$G_d[D_d]\:(d\ge 0)$ (that is, the ideal containing these graphs and all their induced subgraphs). Since $G_d[D_d]$
has chromatic number more than $d$, $\mathcal{C}''$ is non-colourable; and so it has rank at least $k$, from the choice of $\mathcal{C}'$.
In particular, it has rank at least one. Consequently there exists $c_2\ge 0$ such that,
for all $c'\ge 0$, there exists $d\ge 0$, and $C\subseteq D_d$ with $\chi(C)>c'$, and
a $\xi$-clique-cover $\mathcal{L}_{k+1}=(X_{k+1},N_{k+1})$ of $C$ in $G_d[D_d]$,  such that 
$(\mathcal{L}_{k+1}, C)$ is not $c_2$-skewable. We may assume that $c_2\ge c_1$, by replacing $c_2$ by $\max(c_1,c_2)$.
This proves (1) (with $G=G_d$).

\bigskip

Let $c_1,c_2$ be as in (1), let $d\ge \max(c_1, \tau_2)$, let $c'= (d+\beta)(\xi+1)^\zeta$, and let $G,C, D$ etc.\ be as in (1), 
where $\chi(C)>c'$.
Let $W_1=V(G)$, and for $2\le i\le k$ let $W_i$ be the set of vertices in $W_{i-1}$ that are anticomplete to $X_{i-1}$.
Let $W_{k+1}=D$. Then $\mathcal{W}=(W_1\ll W_{k+1})$ is a world for 
$\mathcal{M}=(\mathcal{L}_i\;:1\le i\le k+1)$, $C$.
By \ref{fixlast}
applied to $\mathcal{M}, \mathcal{W}$, and taking $\beta=c_1$, 
there exists $C'\subseteq C$ with $\chi(C')>d$, and a truncation $\mathcal{M}'=(\mathcal{L}_i':i\in I)$ of $M$ covering $C'$,
such that
\begin{itemize}
\item $\mathcal{M}'$ is $c_1$-tidy with respect to $C',\mathcal{W}$;
\item for $1\le i<j\le k$, 
$(\mathcal{L}_i', \mathcal{L}_j')$
is independent with respect to $C'$.
\end{itemize}
Let $\mathcal{L}_i'=(X_i, N_i')$ for $1\le i\le k+1$.
\\
\\
(2) {\em For $1\le i\le k$, $(\mathcal{L}_i', \mathcal{L}_{k+1}')$ is independent with respect to $C'$.}
\\
\\
Suppose not; then since $\mathcal{M}'$ is $c_1$-tidy with respect to $C',\mathcal{W}$, there exists $i\in \{1\ll k\}$
such that $(\mathcal{L}_i', \mathcal{L}_{k+1}')$ is $c_1$-skew with respect to 
$\mathcal{M}',C',\mathcal{W}$.
We claim that this shows that $(\mathcal{L}_i,D)$ is $c_1$-skewable. To show this, we must check:
\begin{itemize}
\item $N_i'\subseteq N_i$, and $C'\subseteq D$ with $\chi(C')>c_1$;
\item $X_{k+1},N_{k+1}'\subseteq D$;
\item $Z$ is anticomplete to $C'$,
and every vertex in $N_{k+1}'$ is $(c_1,\xi)$-earthed via $(Z, D)$, where $Z$ is the set of 
vertices in $N_i'$ that are not complete to $X_{k+1}$.
\end{itemize}
The first two are clear, since $d\ge c_1$. Since $(\mathcal{L}_i', \mathcal{L}_{k+1}')$ is $c_1$-skew with respect to 
$\mathcal{M}',C',\mathcal{W}$, it follows that $Z$ is
anticomplete to $C'$, and every vertex in $N_{k+1}$ is $(c_1,\xi)$-earthed via $(Z, W_{k+1})$. Since $W_{k+1}=D$, 
this shows the claim.
Consequently $(\mathcal{L}_i,D)$ is $c_1$-skewable, a contradiction. This proves (2).

\bigskip

From (2), we have shown that for all $d\ge  \max(c_1, \tau_2)$ (and hence for all $d\ge 0$) there exist 
$G\in \mathcal{C}'$, and $C'\subseteq V(G)$ with $\chi(C')>d$, and
an independent $\xi$-clique-multicover $(\mathcal{L}_i:1\le i\le k+1)$ of $C'$ in $G$, 
such that for $1\le i\le k+1$, $(\mathcal{L}_i, C')$ is not $c_2$-skewable. But this contradicts that $\mathcal{C}'$
has rank $k$.
This proves that $\mathcal{C}'$ has rank zero, and so satisfies the theorem; and hence proves \ref{skewpair}.~\bbox

\begin{thm}\label{getlongcable}
Let $\xi,t>0$ and $\tau_2\ge 0$, and let $\mathcal{C}$ be a skewable ideal of graphs (relative to $\xi$), such that
for each $G\in \mathcal{C}$, 
$\chi(N^2_G(X))\le \tau_2$ for every $(\xi+1)$-clique $X$ in $G$.
Then for all $\beta,c'\ge 0$ there exists $c\ge 0$ with the following property.
Let $G\in \mathcal{C}$, and let $\mathcal{L}$ be a $\xi$-clique-cover of $C$ in $G$, where $\chi(C)>c$.
Then there exist $C'\subseteq C$ with $\chi(C')>c'$, and 
a $\xi$-clique-multicover $\mathcal{M}=(\mathcal{L}_i:1\le i\le t)$ of $C'$, and a world $\mathcal{W}$
for $\mathcal{M}, C'$, such that:
\begin{itemize}
\item 
$\mathcal{L}_1$ is a truncation of $\mathcal{L}$;
\item $V(\mathcal{L}_i)\subseteq C$ for $2\le i\le t$;
\item $\mathcal{M}$ is $\beta$-tidy with respect to $C', \mathcal{W}$; and
\item for $2\le i\le t$, the pair $(\mathcal{L}_1,\mathcal{L}_i)$ is $\beta$-skew with respect to $\mathcal{M},C', \mathcal{W}$.
\end{itemize}
\end{thm}
\Proof The result is true when $t=1$, taking $c'=c$; so we assume that $t>1$ and the result holds for $t-1$. 
Let $\beta,c'\ge 0$. We may assume that $c'\ge \tau_2$.
Since $\mathcal{C}$ is skewable, there exists $c_0$ such that $(\mathcal{L},C)$ is $(c'+\beta)(\xi+1)^t$-skewable
for every $G\in \mathcal{C}$, every $C\subseteq V(G)$ with $\chi(C)>c_0$, and every $\xi$-clique-cover
$\mathcal{L}$ of $C$ in $G$.

Choose a value of $c$ such that the result holds with $t,c',c$ replaced by $t-1,c_02^t,c$ respectively.
We claim that $c$ satisfies the theorem. For let $G,C$ and $\mathcal{L}=(X,N)$
be as in the theorem, with
$\chi(C)>c$. From the choice of $c$, there exist $D_1\subseteq C$ with $\chi(D_1)>c_02^t$, and
a $\xi$-clique-multicover $\mathcal{M}_1=(\mathcal{L}_1'',\mathcal{L}_2\ll \mathcal{L}_{t-1})$ of $D_1$, and a world
$\mathcal{W}_1$ for $\mathcal{M}_1, D_1$, such that
\begin{itemize}
\item $\mathcal{L}_1''$ is a truncation of $\mathcal{L}$;
\item $V(\mathcal{L}_i)\subseteq C$ for $2\le i\le t-1$;
\item $\mathcal{M}_1$ is $\beta$-tidy with respect to $D_1, \mathcal{W}_1$; and
\item for $2\le i\le t-1$, the pair $(\mathcal{L}_1'',\mathcal{L}_i)$ is $\beta$-skew with respect to
$\mathcal{M}_1,D_1,\mathcal{W}_1$.
\end{itemize}
Choose $D_2\subseteq D_1$ with chromatic number at least $2^{-t}\chi(D_1)>c_0$, such that for $2\le i\le t-1$, either all vertices in $D_2$
are $(\beta+\tau_2,\xi)$-earthed via $(N(\mathcal{L}_i),D_1)$, or none are.

Let $\mathcal{W}_1=(W_1\ll W_{t-1})$, and define $W_t=D_1$ and $\mathcal{W} = (W_1\ll W_t)$.
Now $\mathcal{L}_1''$
is a $\xi$-clique-cover of $D_2$, and $\chi(D_2)>c_0$,
and so $(\mathcal{L}_1'',D_2)$ is $(c'+\beta)(\xi+1)^t$-skewable, by the choice of $c_0$. Hence
there exist $D_3\subseteq D_2$ with $\chi(D_3)>(c'+\beta)(\xi+1)^t$, and a truncation $\mathcal{L}_1$ of
$\mathcal{L}_1''$ covering $D_3$,
and a $\xi$-clique-cover $\mathcal{L}_t$ of $D_3$,
such that $V(\mathcal{L}_t)\subseteq D_2$, and the $\xi$-clique-multicover (of length two)
$(\mathcal{L}_1,\mathcal{L}_t)$ is $(c'+\beta)(\xi+1)^t$-skew, and hence $\beta$-skew, with respect to $D_3$ and the world $(V(G), D_2)$.
Now $N(\mathcal{L}_1)\subseteq N(\mathcal{L}_1'')$, and we may assume that every vertex in $N(\mathcal{L}_1'')\setminus N(\mathcal{L}_1)$
has a neighbour in $D_2$; because if some $v\in N(\mathcal{L}_1'')\setminus N(\mathcal{L}_1)$ has no neighbour in $D_2$, then
we can add it to $N(\mathcal{L}_1)$ preserving all the conditions.

Let
$$\mathcal{M}_2=(\mathcal{L}_1, \mathcal{L}_2, \mathcal{L}_3\ll \mathcal{L}_{t-1})$$
and
$$\mathcal{M}_3=(\mathcal{L}_1, \mathcal{L}_2, \mathcal{L}_3\ll \mathcal{L}_{t-1}, \mathcal{L}_t);$$
these are both $\xi$-clique-multicovers of $D_3$.
Also, $\mathcal{W}_1$ is a world for $\mathcal{M}_2, D_3$, and $\mathcal{W}$ is a world for
$\mathcal{M}_3,D_3$. 
\\
\\
(1) {\em Every pair of $\mathcal{M}_3$ is $\beta$-tidy with respect to $\mathcal{M}_3,D_3, \mathcal{W}$ except possibly the pairs
$(\mathcal{L}_i,\mathcal{L}_t)$ where $2\le i\le t-1$; and in particular, for $2\le i\le t$, the pair
$(\mathcal{L}_1, \mathcal{L}_i)$ is $\beta$-skew with respect to $\mathcal{M}_3,D_3, \mathcal{W}$.}
\\
\\
To see this, there are three kinds of pairs to consider:
\begin{itemize}
\item The pair $(\mathcal{L}_1, \mathcal{L}_i)$ where $2\le i\le t-1$: the pair
$(\mathcal{L}_1'',\mathcal{L}_i)$ is $\beta$-skew with respect to $\mathcal{M}_1,D_1, \mathcal{W}_1$, and therefore
$(\mathcal{L}_1, \mathcal{L}_i)$ is $\beta$-skew with respect to $\mathcal{M}_2, D_3, \mathcal{W}_1$, 
since every vertex in $N(\mathcal{L}_1)\setminus N(\mathcal{L}_1'')$
has a neighbour in $D_2$.
Since $W_t=D_1$, it is also $\beta$-skew with respect to $\mathcal{M}_3, D_3,\mathcal{W}$.
\item The pair $(\mathcal{L}_1, \mathcal{L}_t)$: this is $\beta$-skew with respect to $\mathcal{M}_3, D_3, \mathcal{W}$, since
as a $\xi$-clique-multicover, it is $\beta$-skew with
respect to $D_3$ and the world $(V(G), D_2)$.
\item The pair $(\mathcal{L}_i, \mathcal{L}_j)$ where $2\le i<j\le t-1$: this is $\beta$-tidy with respect to $\mathcal{M}_1, D_1,\mathcal{W}_1$,
and therefore with respect to $\mathcal{M}_2, D_3,\mathcal{W}_1$;
and hence also with respect to $\mathcal{M}_3, D_3,\mathcal{W}$ since $W_t= D_1$.
\end{itemize}
This proves (1).

\bigskip

By \ref{fixlast}
we deduce that there exist $D_4\subseteq D_3$ with $\chi(D_4)>c'$, 
and a truncation $\mathcal{L}_i'$ of $\mathcal{L}_i$ for $1\le i\le t$, such that $\mathcal{M}=(\mathcal{L}_1'\ll \mathcal{L}_t')$ is a 
$\beta$-tidy $\xi$-clique-multicover
of $D_4$, and $\mathcal{W}$ is a world for $\mathcal{M}, D_4$, and 
for $2\le i\le t$ each pair $(\mathcal{L}_1', \mathcal{L}_i')$ is $\beta$-skew with respect to 
$\mathcal{M},D_4,\mathcal{W}$.
This proves \ref{getlongcable}.~\bbox

By choosing $t$ large enough in \ref{getlongcable}, and applying Ramsey's theorem to the sequence $(\mathcal{L}_2\ll \mathcal{L}_t)$,
we deduce since $G$ is $(\xi,\zeta, \tau_3)$-free that the same result as \ref{getlongcable} is true with ``$\beta$-tidy''
replaced by ``$\beta$-skew''. This result is important enough that it deserves to be said explicitly:

\begin{thm}\label{getlongskew}
Let $\xi,t\ge 1$ and $\tau_2,\tau_3,\zeta\ge 0$, and let $\mathcal{C}$ be a skewable ideal of graphs, such that
for each $G\in \mathcal{C}$,
\begin{itemize}
\item $\chi(N^2_G(X))\le \tau_2$ for every $(\xi+1)$-clique $X$ in $G$; and
\item $G$ is $(\xi,\zeta,\tau_3)$-free.
\end{itemize}
Then for all $\beta,c'\ge 0$ there exists $c\ge 0$
with the following property. Let $G\in \mathcal{C}$, and
let $\mathcal{L}$ be a $\xi$-clique-cover of $C\subseteq V(G)$, where $\chi(C)>c$.
Then there exist $C'\subseteq C$ with $\chi(C')>c'$, and
a $\xi$-clique-multicover $\mathcal{M}=(\mathcal{L}_i:1\le i\le t)$ of $C'$, and a world $\mathcal{W}$
for $\mathcal{M}, C'$, such that:
\begin{itemize}
\item $\mathcal{L}_1$ is a truncation of $\mathcal{L}$;
\item $V(\mathcal{L}_i)\subseteq C$ for $2\le i\le t$;
\item $\mathcal{M}$ is $\beta$-skew with respect to $C', \mathcal{W}$.
\end{itemize}
\end{thm}
\Proof
Choose an integer $s\ge 0$ such that for every partition of the edges of $K_{s-1}$ into two classes, either some $K_{t-1}$ subgraph has
all its edges in the first class, or some $K_{\zeta}$ subgraph has all its edges in the second.
Let $c$ satisfy \ref{getlongcable} with $t$ replaced by $s$, and $c'$ replaced by $\max(c', \tau_3)$.
We claim that $t$ satisfies the theorem. Let
$G, \mathcal{L}$ and $C$ be as in the theorem. By \ref{getlongcable}
there exist $C'\subseteq C$ with $\chi(C')>\max(c', \tau_3)$, 
and
a $\xi$-clique-multicover $\mathcal{M}'=(\mathcal{L}_i:1\le i\le t)$ of $C'$, and a world $\mathcal{W}'$
for $\mathcal{M}', C'$, such that:
\begin{itemize}
\item
$\mathcal{L}_1$ is a truncation of $\mathcal{L}$;
\item $V(\mathcal{L}_i)\subseteq C$ for $2\le i\le t$;
\item $\mathcal{M}$ is $\beta$-tidy with respect to $C', \mathcal{W}'$; and
\item for $2\le i\le t$, the pair $(\mathcal{L}_1,\mathcal{L}_i)$ is $\beta$-skew with respect to $\mathcal{M}',C', \mathcal{W}'$.
\end{itemize}
For each pair $(i,j)$ with $2\le i<j\le s$, the pair $(\mathcal{L}_i, \mathcal{L}_j)$ is $\beta$-tidy with respect to
$\mathcal{M}', C'$, and so is either independent with respect to $C'$, or $\beta$-skew with respect to $\mathcal{M}', C',\mathcal{W}'$.
From the choice of $s$, either
\begin{itemize}
\item there exists $I\subseteq \{2\ll s\}$ with $|I|=t-1$ such that
$(\mathcal{L}_i, \mathcal{L}_j)$ is $\beta$-skew with respect to $\mathcal{M}, C',\mathcal{W}'$ for all $i<j$ with $i,j\in I$, or
\item there exists $J\subseteq \{2\ll s\}$ with $|J|=\zeta$ such that
$(\mathcal{L}_i, \mathcal{L}_j)$ is independent with respect to $C$, for all $i<j$ with $i,j\in J$.
\end{itemize}
The second is impossible, since $G$ is $(\xi,\zeta, \tau_3)$-free and $\chi(C')>\tau_3$, and so the first holds.
Let $\mathcal{W}' = (W_1 \ll W_t)$, and define $\mathcal{W} = (W_i:i\in \{1\}\cup I)$.
Then every pair of terms in $\mathcal{M}=(\mathcal{L}_i:i\in \{1\}\cup I)$ is $\beta$-skew
with respect to $\mathcal{M}, C',\mathcal{W}$,
and so $\mathcal{M}$ is $\beta$-skew with respect to $C',\mathcal{W}$.
This proves \ref{getlongskew}.~\bbox

\section{Finding a tree of lamps}

Now we come to reap the benefit of all the complications of \ref{getlongskew}: we show that any graph satisfying the conditions
of \ref{getlongskew} contains any given tree of lamps as an induced subgraph,
if the number $t$ and the chromatic number are large enough.

First, we need two lemmas.
Let $\mathcal{M}=(\mathcal{L}_1,\mathcal{L}_2\ll \mathcal{L}_t)$ be a $\xi$-clique-multicover of $C\subseteq V(G)$,
that is $\beta$-skew with respect to $C, \mathcal{W}$.
For $1\le i\le t$, let $\mathcal{L}_i=(X_i, N_i)$, and let $\mathcal{W}=(W_1\ll W_t)$. Define 
$W_{t+1}=C$ (thus, $C\cup W_{j+1}\cup\cdots\cup W_t=W_{j+1}$ for all $j\in \{1\ll t\}$).
For $1\le i<j\le t$, let $Z_{i,j}$ be the set of vertices in $N_i$
that have a neighbour in $W_j$ and are anticomplete to $W_{j+1}$.
We call the family of sets $Z_{i,j}\; (1\le i<j\le t)$ the {\em standard refinement} of $\mathcal{M}, C$.
\begin{thm}\label{breakdown}
In the notation just given:
\begin{itemize}
\item the sets $Z_{i,i+1}\ll Z_{i,t}$ are pairwise disjoint subsets of $N_i$;
\item $X_j$ is complete to $Z_{i,k}$ for $1\le i\le j<k\le t$, and to every vertex in $N_i$ with a neighbour in $C$, 
for $1\le i\le j$;   
\item $X_j$ is anticomplete to $Z_{i,k}$ for all $i,j,k\in \{1\ll t\}$ with $i<k$ if $j<i$ or $k<j$; and 
\item every vertex in $N_j$ is $(\beta,\xi)$-earthed via $(Z_{i,j}, W_j)$ for $1\le i<j\le t$.
\end{itemize}
\end{thm}
\Proof
The first statement is clear from the definition. 
Let $1\le i<j\le t$, and let $Z$ be the set of all vertices in $N_i$ anticomplete to $W_{j+1}$. 
Thus $Z= Z_{i,i+1}\cup\cdots\cup Z_{i,j}\cup U_i$, where $U_i$ is the set of vertices in $N_i$ anticomplete to
$W_{i+1}$. 
From the definition of ``$\beta$-skew'', every vertex in $N_i\setminus Z$ is complete to $X_j$, so the second
statement follows if $i<j$; and if $i=j$ then it follows since $X_i$ is complete to $N_i$. 
Now $X_j$ is anticomplete to $Z_{i,k}$ if $j<i$ from the definition of a $\xi$-clique-multicover;
and $X_j$ is anticomplete to $Z_{i,k}$ if $k<j$, since $Z_{i,k}$ is anticomplete to $W_{k+1}\supseteq X_j$, 
so the third statement follows. From the definition of ``$\beta$-skew'', every vertex in $N_j$ is $(\beta,\xi)$-earthed via
$(Z,W_j)$, and since $Z_{i,j}$ includes the set of all vertices in $Z$ that have a neighbour in $N_j$, the fourth statement follows.
This proves \ref{breakdown}.~\bbox

\begin{thm}\label{findtreelemma}
Let $\xi\ge 1$ and $\tau_1,\tau_2, \beta\ge 0$.
Let $G$ be such that
\begin{itemize}
\item $\chi(N^1_G(v))\le \tau_1$ for every $v\in V(G)$; and
\item $\chi(N^2_G(X))\le \tau_2$ for every $(\xi+1)$-clique $X$ in $G$;
\end{itemize}
Let $\mathcal{W}=(W_1\ll W_t)$, define $W_{t+1}=C\subseteq V(G)$, 
let $\mathcal{M}=(\mathcal{L}_1,\mathcal{L}_2\ll \mathcal{L}_t)$ be a $\xi$-clique-multicover of $C$
that is $\beta$-skew with respect to $C,\mathcal{W}$, 
and let $Z_{i,j}\; (1\le i<j\le t)$ be its
standard refinement.
Let $1\le i<j\le t$, and let
$$r\in \left(\bigcup_{1\le h< i} X_h\cup (N_h\setminus Z_{h,i})\right)\cup \left(\bigcup_{i\le h< j}N_h\right)\cup W_{j+1}.$$
Let $A$ be the set of vertices in $V(G)$
that are equal or adjacent to $r$, or
have a neighbour in $Z_{i,j}$ adjacent to $r$.
Then $\chi(A)\le  \tau_2+(\xi+1)(\tau_1+1)$.
\end{thm}
\Proof
If $r$ has no neighbour in $Z_{i,j}$
then every vertex in $A$ is equal to or adjacent to $r$ and hence $\chi(A)\le \tau_1+1$ and the result holds.
So we may assume that $r$ has a neighbour in $Z_{i,j}$, and so $r\notin W_{j+1}$; choose $h\in \{1\ll j-1\}$ with $r\in X_h\cup N_h$.
\\
\\
(1) {\em One of $X_h, X_i$ is complete to $Z_{i,j}\cup \{r\}$.}
\\
\\
If $i\le h<j$, then $r\in N_h$ by hypothesis, and then $X_h$ is complete
to $r$ and to $Z_{i,j}$ by \ref{breakdown}; so we assume that $h<i$. Then since $r$ has a neighbour in $Z_{i,j}$, 
it follows that $r\in N_h$. If $r$ is complete to $X_i$ then the claim holds, so we assume not. 
Consequently
\ref{breakdown} implies that $r$ has no neighbour in $C$; and therefore
$r\in Z_{h,k}$ for some $k$. Again, since $r$ is not complete to $X_i$, \ref{breakdown} implies that $k\le i$. 
Since $r$ has a neighbour in $N_i$, it follows that $k=i$, contrary to the hypothesis. This proves (1).

\bigskip

Let $X$ be a $\xi$-clique that is complete to $Z_{i,j}\cup \{r\}$.
Since $N^2_G(X\cup \{r\})\le \tau_2$ (because $X\cup \{r\}$ is a $(\xi+1)$-clique),
and $X$ is complete to $Z_{i,j}$,
it follows that the set of vertices in $A$ that are adjacent to a neighbour of $r$ in $Z_{i,j}$ and anticomplete to $X\cup \{r\}$
has chromatic number at most $\tau_2$. But the chromatic number of the set of vertices in $A$ that belong to or have a neighbour
in $X\cup \{r\}$ is at most $(\xi+1)(\tau_1+1)$; and so $\chi(A)\le \tau_2+(\xi+1)(\tau_1+1)$.
This proves \ref{findtreelemma}.~\bbox

Here at last is a definition of a tree of lamps. (See figure 2.) 
Start with a tree $T$, and select a vertex of $T$ called the {\em root};
then every vertex different from the root has a unique {\em parent}, its neighbour on the path towards the root.
Take a map $w$ from $V(T)$ into the set of positive integers, such that
\begin{itemize}
\item for all $u,v\in V(T)$, if $v$ is the parent of $u$ then $w(v)>w(u)$
(and consequently the $w$-value of the root is strictly larger than all the other values);
\item there is a vertex $v$ with $w(v)=1$ (necessarily, either $v$ is the root and $|V(T)|=1$, or $v$ is a leaf of $T$);
\item for all vertices $u,v$ with $u\ne v$, if $w(u)=w(v)$ then $w(u)=1$.
\end{itemize}
We call such a function $w$ a {\em height function} for $T$.
Let $w(V(T))$ denote the set $\{w(v):v\in V(T)\}$.

Now choose a set $J$ of integers, each at least 1 and at most the $w$-value of the root, with $J\cap w(V(T))=\{1\}$.
For each $j\in J$, take a new vertex $x_j$; and make $x_j$ adjacent to $v$ for every edge $uv$ of $T$ such that $w(v)\le j$ and $w(u)>j$.
(If $|V(T)|=1$, make $x_1$ adjacent to the root.)
A graph constructed this way is called a {\em lamp}, and $x_1$ is its {\em plug}.
Thus every chandelier is a lamp, but many lamps are not chandeliers.

Analogously to trees of chandeliers, we can make trees of lamps, by taking a new lamp, and 
attaching trees of lamps already constructed to this new lamp by their plugs.
However, we are not permitted to attach anything to neighbours of the plug of the new lamp. Let us say this more precisely.
A {\em spotlight} is a one-vertex graph, with plug its vertex.
No tree of lamps has negative height; and the spotlight is the only 
tree of lamps of height zero.
Inductively for $r>0$, having defined trees of lamps of height $\le r-1$ and their plugs, we proceed as follows.
Let $L$ be a lamp with plug $\ell$.
For each  $v\in V(L)$, let $Q_v$ be a tree of lamps of height at most $r-1$, 
such that all the graphs $L$ and $Q_v\; (v\in V(L))$
are pairwise anticomplete, and such that if $v$ is equal to or adjacent to $\ell$, then $Q_v$ is a spotlight. 
Now identify $v$ with the plug of $Q_v$, for each $v\in V(L)$. (More precisely, add new edges
joining $v$ to every neighbour of the plug of $Q_v$, and then delete the plug of $Q_v$, for each $v\in V$.) Let the result be $Q$.
Any such graph $Q$, with plug $\ell$, is said to be a tree
of lamps of height $\le r$ (and so is the spotlight). 

We mentioned earlier that we think that not every tree of chandeliers is a tree of lamps; the reason for this (if true) is the
more restrictive composition rule. In fact, there is a third class: we have
\begin{itemize}
\item trees of lamps (call this $\mathcal{A}$)
\item connected induced subgraphs of trees of lamps ($\mathcal{B}$)
\item trees of chandeliers ($\mathcal{C}$).
\end{itemize}
Evidently $\mathcal{A}\subseteq \mathcal{B}$, but we are not sure whether equality holds, or whether $\mathcal{C}$
is a subclass of either of the other two, although we expect the answer is ``no'' in each case.

We used earlier the fact that for every tree of lanterns $H$, there is a tree of lamps $Q$ such that some subdivision
of $H$ is an induced subgraph of $Q$. We leave it to the reader to verify this. (When growing a tree of lanterns, there is
no need to attach new lanterns to the pivot of what we have already built, because a graph formed by
two lanterns with their pivots identified
is an induced subgraph of one bigger lantern with the same pivot. So, grow it adding one lantern at a time, and identifying 
the pivot of the new lantern with a non-pivot vertex of what we have already built. Now change this; for each new lantern
that we want to attach, first subdivide all the edges incident with its pivot and attach that instead. What we 
construct is a tree of lamps that is a subdivision of our original tree of lanterns.)

We will show the following.

\begin{thm}\label{gettree}
Let $\xi>0$ and $\tau_1,\tau_2,\tau_3, \zeta\ge 0$, let $\phi:\mathbb{N}\rightarrow \mathbb{N}$ be a non-decreasing function, and
let $Q$ be a tree of lamps. Then there exists $c\ge 0$ with the following property.
Let $G$ be a graph such that:
\begin{itemize}
\item $G$ is $\phi$-skewable, relative to $\xi$;
\item $\chi(N^1_G(v))\le \tau_1$ for every $v\in V(G)$;
\item $\chi(N^2_G(X))\le \tau_2$ for every $(\xi+1)$-clique $X$ in $G$; and
\item $G$ is $(\xi,\zeta, \tau_3)$-free.
\end{itemize}
Let $\mathcal{L}_0$ be a $\xi$-clique-cover of $C\subseteq V(G)$, where $\chi(C)>c$, 
and let $a\in X(\mathcal{L}_0)$.
Then there is an isomorphism from $Q$ to an induced subgraph of $G$, mapping the plug of $Q$ to $a$ and mapping
all other vertices of $Q$ into $N(\mathcal{L}_0)\cup C$.
\end{thm}
\Proof We proceed by induction on $|V(Q)|$.
Certainly it is true if $|V(Q)| = 1$, so we assume that $|V(Q)|>1$ and the result holds for all smaller trees of lamps. Since, up to
isomorphism, there are only finitely many smaller trees of lamps, we can choose $c_0\ge 0$ such that
the theorem is true with $c$ replaced by $c_0$ for every tree of lamps with at most $|V(Q)|-1$ vertices.

There is a lamp $L$ with plug $\ell$ say, and
trees of lamps $Q_v\;(v\in V(L))$ such that $Q$ is obtained from $L$ and the graphs
$Q_v\;(v\in V(L))$ as in the definition.

There is a tree $T$, a height function $w$, a set $J$ of integers, and vertices $x_j\;(j\in J)$ in $L$,
as in the definition of a lamp. Choose $w$ and $J$ such that $w(v)$ is congruent to $1$ modulo $3$ for all $v$, and 
every member of $J$ is also congruent to 1 modulo 3.
Let $q_0$ be the root of $T$, and let $t=w(q_0)$. 

Let $\mathcal{C}$ be the ideal of all graphs that satisfy the four bullets of the theorem. Thus $\mathcal{C}$ is skewable.
Let $\beta = c_0 + |V(Q)|(\tau_2+(\xi+1)(\tau_1+1))$, and choose $c$ such that \ref{getlongskew} holds for $\mathcal{C}$, taking
$c'=0$.
We claim that $c$ satisfies the theorem.

Let $G, \mathcal{L}_0$ and $C$ be as in the theorem. By \ref{getlongskew}, applied to the graph $G[V(\mathcal{L}_0)\cup C]$
and to $\mathcal{C}$,
there exist $C'\subseteq C$ 
with $\chi(C')>0$, and 
a $\xi$-clique-multicover $\mathcal{M}=(\mathcal{L}_i:1\le i\le t)$ of $C'$, and a world $\mathcal{W}$
for $\mathcal{M}, C'$, such that:
\begin{itemize}
\item $\mathcal{L}_1$ is a truncation of $\mathcal{L}_0$;
\item $V(\mathcal{L}_i)\subseteq C$ for $2\le i\le t$;
\item $\mathcal{M}$ is $\beta$-skew with respect to $C', \mathcal{W}$; and
\item every term of $\mathcal{W}$ is a subset of $V(\mathcal{L}_0)\cup C$.
\end{itemize}
For $1\le i\le t$ let $\mathcal{L}_i=(X_i, N_i)$, and let $Z_{i,j}\; (1\le i<j\le t)$ be the
standard refinement of $\mathcal{M}, C'$.

Now we begin to construct the isomorphism $\eta$ from $Q$ to an induced subgraph of $G$. We recall that $q_0$ is the root of $T$;
choose some vertex in $N_{t}$, and call it $\eta(q_0)$. At a general stage of the process, we will have defined $\eta(p)$ only for
the vertices $p$ in a subset $\dom(\eta)$ of $V(Q)$.
We will ensure that $\eta$ is injective, and for all $u,v\in \dom(\eta)$,
$u,v$ are adjacent in $Q$ if and only if $\eta(u), \eta(v)$ are adjacent in $G$. If $|V(T)|=1$, then $|J|=1$, and
(since no pendant lamp
can be attached at the plug or at one of its neighbours) it follows that $|V(Q)|\le 2$ and the claim is trivial; so we may
assume that $|V(T)|\ge 2$.

First we extend $\dom(\eta)$ to equal $V(T)$, in such a way that $\eta(p)\in N_{w(p)}$ for each $p\in V(T)$, 
by repeating the following process.
\begin{itemize}
\item Choose an integer $n$ maximum such that $w(v)=n$ for some $v\in V(T)\setminus \dom(\eta)$. (When $\dom(\eta)=V(T)$, stop).
\item Let $u$ be the neighbour of $v$ in $\dom(\eta)$ (necessarily unique). Note that $w(v)<w(u)$.
\item Choose a vertex $y\in Z_{w(v),w(u)}$ adjacent to $\eta(u)$
and nonadjacent to
all the vertices $\eta(p)\;(p\in \dom(\eta)\setminus \{u\})$. To see that this is possible, let $p\in \dom(\eta)\setminus \{u\}$. 
Since $w(u)>w(v)\ge 1$, and therefore 
$w(p)\ne w(u)$, it follows from \ref{findtreelemma},  and from the fact that $\eta(p)\in N(w(p))$, that the set of vertices in $V(G)$ that have a neighbour in 
$Z_{w(v),w(u)}$ adjacent to $\eta(p)$ has chromatic number at most $\tau_2+(\xi+1)(\tau_1+1)$. Consequently the set of
vertices in $W_{w(u)}$ that have a neighbour in
$Z_{w(v),w(u)}$ with a neighbour in $\{\eta(p):p\in \dom(\eta)\setminus \{u\}\}$ has chromatic number at most 
$|V(Q)|(\tau_2+(\xi+1)(\tau_1+1))$. 
Since $\eta(u)$ is $(\beta,\xi)$-earthed via $(Z_{w(v), w(u)}, W_{w(u)})$ by \ref{breakdown}, and $\beta\ge |V(Q)|(\tau_2+(\xi+1)(\tau_1+1))$, there
is at least one vertex $x\in W_{w(u)}$ that has a neighbour $y\in Z_{w(v), w(u)}$ adjacent to $\eta(u)$,
and has no neighbour in $Z_{w(v), w(u)}$ that is adjacent to any of $\eta(p)\;(p\in \dom(\eta)\setminus \{u\})$. In particular,
$y$ is nonadjacent to all of $\eta(p)\;(p\in \dom(\eta)\setminus \{u\})$. This shows the existence of the vertex $y$ as claimed. 

\item Define $\eta(v)=y$, and add $v$ to $\dom(\eta)$.
\end{itemize}
Note that for all $i,j$ with $1\le i<j\le t$, if some vertex of $T$ is mapped into $Z_{i,j}$ by $\eta$, then both $i,j$ are equal
to 1 modulo 3.

Next we add all the vertices $x_j\;(j\in J)$ to $\dom(\eta)$, defining $\eta(x_j)$ to be some vertex in $X_j$ for each $j\in J$, and
in particular choosing $\eta(x_1)=a$. We claim that 
$\eta$ still defines an isomorphism from $\dom(\eta)$ into $G$. To see this, let  $j\in J$ and $v\in V(T)$. 
We must check that $\eta(x_j), \eta(v)$ are adjacent if and only if $v$ has a parent $u$ in $T$ and 
$w(u)>w(x_j)\ge w(v)$.
Let $\eta(v)\in Z_{i,k}$ say. If $i>j$ then $\eta(x_j), \eta(v)$ are nonadjacent since $X_j$ is anticomplete to 
$N_i$; so we may assume that $i\le j$.  Consequently, if $v$ has no parent, then $i=1$ and $|V(T)|=1$, a contradiction;
so $v$ has a parent $u$. From the construction, $\eta(u)\in N_k$. Now $Z_{i,k}$ is 
anticomplete to $X_j$ if $k<j$, from \ref{breakdown}, so we may assume that $j\le k$; and so $j<k$ since $k\ne 1$. Thus
$i\le j<k$; and so $\eta(x_j), \eta(v)$ are adjacent since $X_j$ is complete to $Z_{i,k}$ by \ref{breakdown}. This proves that
we can add all the vertices $x_j\;(j\in J)$ to $\dom(\eta)$ so that $\eta$ still defines an isomorphism.
At this stage, then, $\dom(\eta)=V(L)$.

Now we turn to adding the ``pendant'' trees of lamps $Q_v\;(v\in V(L))$. The plug of each $Q_v$, namely $v$, already belongs to
$\dom(\eta)$, and we must add the other vertices of $Q_v$; and we shall do so mapping $V(Q_v)\setminus \{v\}$ into 
$W_{w(v)-1}$. We do them in order:
for $n=t, t-3, t-6\ll 1$ in turn, if there is a vertex $v\in \dom(\eta)$ with $w(v)=n$, we shall extend $\dom(\eta)$ to include 
$V(Q_v)\setminus \{v\}$. If $n=1$, then since all the $Q_v$ are spotlights when $w(v)=1$, the process stops.
At the start of a general step of the process, $n\ge 2$ and $n=1$ modulo 3. 
Let $R=\{\eta(v):v\in \dom(\eta)\}$; then $|R|\le |V(Q)|$,
and every $r\in R$ belongs either to $W_{n+2}$, or to some $X_i\cup N_i$ where $i\le  n$ and $i=1$ modulo 3. Moreover,
if $R\cap Z_{h,i}\ne \emptyset$ where $h\le n+1$, then both $h,i$ equal 1 modulo 3.

If there is no $v\in V(L)$ with $w(v)=n$, go on to the next value of $n$. So now, there is such a vertex $v$, unique since
$n>1$, and $\eta(v)\in X_n\cup N_n$. Either $v\in V(T)$ or $v=x_n$; the arguments in the 
two cases are almost identical,
but slightly different (this is why we need two values of $m$ in (1)).
\\
\\
(1) {\em For each $r\in R\setminus \{\eta(v)\}$, and for $m=n,n+1$, the set of vertices in $V(G)$ that either are equal or adjacent to $r$, or have a neighbour in $Z_{n-1,m}$ 
adjacent to $r$, has chromatic number at most $\tau_2+(\xi+1)(\tau_1+1)$.}
\\
\\
Let $r\in R\setminus \{\eta(v)\}$. Then $r$ belongs either to $W_{n+2}$, or to some $X_i\cup N_i$ where $i<n$ and $i=1$ modulo 3. 
Moreover,
if $R\cap Z_{h,i}\ne \emptyset$ where $h\le n+1$, then both $h,i$ equal 1 modulo 3. Since $W_{n+2}\subseteq W_{m+1}$,
and $n-1$ does not equal 1 modulo 3, it follows that
$$r\in \left(\bigcup_{1\le h< n-1} X_h\cup (N_h\setminus Z_{h,n-1})\right)\cup \left(\bigcup_{n-1\le h< m}N_h\right)\cup W_{m+1}.$$
Hence the claim follows from \ref{findtreelemma}. This proves (1). 

\bigskip
Now there are two cases, depending whether $v\in V(T)$ or $v=x_n$.
\begin{itemize}
\item Assume that $v\in V(T)$. 
Since $\eta(v)$ is $(\beta,\xi)$-earthed via $(Z_{n-1, n}, W_{n})$, by \ref{breakdown},
there is a $\xi$-clique $X\subseteq W_n$ with $\eta(v)\in X$ such that $X$ is $\beta$-earthed via $(Z_{n-1, n}, W_{n})$.
Let $M'$ be the set of vertices in $W_n$ that are anticomplete to $X$ and have a neighbour in $Z_{n-1, n}$ that is complete to $X$; thus
$\chi(M')> \beta$.
Let $Z$ be the set of vertices in $Z_{n-1,n}$ with no neighbour in $R\setminus \{\eta(v)\}$, and let $W$ be
the set of vertices in $W_n$ with no neighbour in $R\setminus \{\eta(v)\}$.
By (1), the set of vertices in $V(G)$ that either belong to $W_n\setminus W$ 
or have a neighbour in $Z_{n-1,n}\setminus Z$ has chromatic number at most
$|Q|(\tau_2+(\xi+1)(\tau_1+1))$. Consequently there exists $M\subseteq M'\cap W$ with 
$$\chi(M)>\beta-|Q|(\tau_2+(\xi+1)(\tau_1+1))=c_0$$
such that $M$ is anticomplete to $Z_{n-1,n}\setminus Z$.
Thus $(X,Z)$ is a $\xi$-clique-cover of $M$, and $\eta(v)\in X$;
and from the inductive hypothesis, there is an isomorphism from $Q_v$ to
an induced subgraph of $G[Z\cup W\cup \{\eta(v)\}]$, mapping the plug of $Q_v$ to $\eta(v)$. This provides the desired extension of
$\eta$ and $\dom(\eta)$ to include $V(Q_v)$. Then go to the next value of $n$.

\item Assume that $v=x_n$, and so $n<t$ and there are vertices in $N_{n+1}$; choose one. Since it is
$(\beta,\xi)$-earthed via $(Z_{n-1,n+1}, W_{n+1})$, by \ref{breakdown}, it follows that
the set $M'$ of vertices in $W_{n+1}$ that have a neighbour in
$Z_{n-1, n+1}$ has chromatic number more than $\beta$. 

Let $Z$ be the set of vertices in $Z_{n-1,n+1}$ with no neighbour in $R\setminus \{\eta(v)\}$, and let $W$ be
be the set of vertices in $W_{n+1}$ with no neighbour in $R\setminus \{\eta(v)\}$.
By (1), the set of vertices in $V(G)$ that either belong to $W_{n+1}\setminus W$
or have a neighbour in $Z_{n-1,n+1}\setminus Z$ has chromatic number at most
$|Q|(\tau_2+(\xi+1)(\tau_1+1))$; and since $\chi(M')>\beta$, it follows that there exists $M\subseteq M'\cap W$ with $\chi(M)>c_0$,
such that $M$ is anticomplete to $Z_{n-1,n+1}\setminus Z$.
Hence $(X_n, Z)$ is a $\xi$-clique-cover of $M$ (because $X_n$ is complete to $Z_{n-1,n+1}$ and anticomplete to $W_{n+1}$).
From the inductive hypothesis, there is an isomorphism from $Q_v$ to
an induced subgraph of $G[Z\cup W\cup \{\eta(v)\}]$, mapping the plug of $Q_v$ to $\eta(v)$. This provides the desired extension of
$\eta$ and $\dom(\eta)$ to include $V(Q_v)$. Then go to the next value of $n$.
\end{itemize}
This completes the construction of the isomorphism, and so completes the proof of \ref{gettree}.~\bbox

\section{Putting the pieces together}

From \ref{gettree}, we deduce:

\begin{thm}\label{findchand5}
Let $\nu,m,\tau_3\ge 0$, and let $\mathcal{C}$ be a non-colourable ideal of graphs
such that
\begin{itemize}
\item $\omega(G)\le \nu$ for all $G\in \mathcal{C}$;
\item $\mathcal{C}$ is $2$-controlled; and
\item all graphs in $\mathcal{C}$ are $(m,\tau_3)$-limited.
\end{itemize}
Then $\mathcal{C}$ contains every tree of lamps.
\end{thm}
\Proof
We proceed by induction on $\nu$. We may assume that $\nu\ge 1$ and the result holds for $\nu-1$. Let $\mathcal{D}$ be the ideal of
all $H\in \mathcal{C}$ with $\omega(H)<\nu$. Thus
by the inductive hypothesis, we may assume that there exists $\tau_1$ such that
all graphs in $\mathcal{D}$ have chromatic number at most $\tau_1$. In particular, for all $G\in \mathcal{C}$,
$\chi(N^1_G(v))\le \tau_1$ for every vertex $v\in V(G)$.

By \ref{maxclique}, there exists $\xi>0$ such that $\mathcal{C}$ is $\xi$-clique-controlled, and
there is a non-colourable subideal $\mathcal{C}'$
of $\mathcal{C}$ and  $\tau_2\ge 0$ such that
$\chi(N^2_G(X))\le \tau_2$ for every $G\in \mathcal{C}'$ and for every $(\xi+1)$-clique $X$ of $G$.

By \ref{getindpt} applied to $\mathcal{C}'$, there exists
$\zeta\ge 0$ such that every graph in $\mathcal{C}'$ is $(\xi,\zeta, \tau_3)$-free. 
By \ref{skewpair} applied to $\mathcal{C}'$, there is a non-colourable subideal $\mathcal{C}''$ of $\mathcal{C}'$ 
such that $\mathcal{C}''$ is skewable. Let $\phi:\mathbb{N}\rightarrow \mathbb{N}$ be a nondecreasing 
function such that every graph in $\mathcal{C}''$
is $\phi$-skewable relative to $\xi$.
Let $Q$ be a tree of lamps, and let $c$ satisfy \ref{gettree}.
Since $\mathcal{C}''$ is $\xi$-clique-controlled, there exists $c'\ge 0$ such that for all $G\in \mathcal{C}''$ with $\chi(G)>c'$,
there is a $\xi$-clique $X_1$ of $G$ with $\chi(N^2_G(X_1))>c$.

Since $\mathcal{C}''$ is non-colourable, there exists $G\in \mathcal{C}''$ with $\chi(G)>c'$. Consequently
there is a $\xi$-clique $X_1$ of $G$ with $\chi(N^2_G(X_1))>c$.
By \ref{gettree}, $G$ contains $Q$ as an induced
subgraph. This proves \ref{findchand5}.~\bbox

Because of \ref{betterusetick1}, we have the corollary:
\begin{thm}\label{findchand6}
Let $\mu,\nu\ge 0$, and let $\mathcal{C}$ be a non-colourable ideal of graphs
such that
\begin{itemize}
\item $\mathcal{C}$ is $2$-controlled;
\item all graphs in $\mathcal{C}$ are $(1,\mu,\nu)$-restricted.
\end{itemize}
Then $\mathcal{C}$ contains every tree of lamps.
\end{thm}
\Proof We proceed by induction on $\nu$; so, as in \ref{findchand5}, 
we may assume that there exists $\tau_1$ such that
all graphs in $\mathcal{D}$ have chromatic number at most $\tau_1$.
Choose $m,d$ as in \ref{betterusetick1}; then since every graph in $\mathcal{C}$ is
$(1,\mu,\nu)$-restricted, they are all $(m,d)$-limited by \ref{betterusetick1}, and the result follows
from \ref{findchand5}.~\bbox

We see that \ref{2control} is an immediate consequence of \ref{findchand6}.
Let us prove \ref{rhopervasive}, which we restate:
\begin{thm}\label{rhopervasiveagain}
For all $\rho\ge 2$, every forest of lanterns is pervasive in every $\rho$-controlled ideal.
\end{thm}
\Proof Let $\mathcal{C}$ be a $\rho$-controlled ideal, let $T$ be a forest of lanterns, and let $\nu,\ell \ge 0$.
We must show that there exists $c$ such that for every graph $G\in \mathcal{C}$ with $\omega(G)\le \nu$ and $\chi(G)>c$,
there is an induced subgraph of $G$ isomorphic to an ($\ge \ell$)-subdivision of $T$.
Let $T_1$ be the $\ell$-subdivision of $T$; then $T_1$ is also a forest of lanterns.
Choose a tree of lamps $Q$ such that some subdivision of $T_1$ is an induced subgraph of $Q$, and choose
$\mu\ge 0$ such that some subdivision of $T_1$ is an induced subgraph of $K_{\mu,\mu}^1$ (and hence
every proper subdivision of $K_{\mu,\mu}$ contains some $(\ge \ell)$-subdivision of $T$ as an induced subgraph).
Let $\mathcal{D}$ be the ideal of graphs $G\in \mathcal{C}$ with clique number at most $\nu$ such that
no induced subgraph of $G$ is an $(\ge \ell)$-subdivision of $T$. It follows that every graph in $\mathcal{D}$ is
$(\rho+2,\mu,\nu)$-restricted, and hence $\mathcal{D}$ is $2$-controlled by \ref{reducecontrol}. By \ref{findchand6} applied to
$\mathcal{D}$ and $Q$, the members of $\mathcal{D}$ have bounded chromatic number. This proves \ref{rhopervasiveagain}.~\bbox

\section{String graphs}

A {\em curve} means a subset of the plane which is homeomorphic to the interval $[0,1]$. Given a finite set $C$ of curves in
the plane, its {\em intersection graph} is the graph with vertex set $C$ in which distinct $S,T\in C$ are adjacent
if $S\cap T\ne \emptyset$; and the intersection graphs of sets of curves are called {\em string graphs}.
Every string graph can be realized by a set of piecewise linear curves, and in this paper, a {\em string} means a piecewise
linear curve. In this section we prove that the ideal of string graphs is $3$-controlled, and consequently the theorems of this
paper can be applied to the ideal. The proof that they are 3-controlled is a modification and simplification of an argument
of McGuinness~\cite{mcguinness}, who showed that a similar statement holds for a triangle-free subideal of string graphs satisfying 
another condition that we omit. 

Let $(v_1\ll v_n)$ be a sequence of distinct vertices of a graph $G$. 
We say that $(v_1\ll v_n)$ has the {\em cross property} if for all
$h,i,j,k$ with $1\le h<i<j<k\le n$, if $P, Q$ are paths of $G$ between $v_h, v_j$ and between $v_i, v_k$ respectively, then 
$V(P)$ is not anticomplete to $V(Q)$.
We need the following.

\begin{thm}\label{discthm}
Let $\Delta$ be a closed disc in the plane, and let $C$ be a finite set of strings all within $\Delta$. 
Let $C_1$ be the set of members 
of $C$ with nonempty intersection with the boundary of $\Delta$. Then $C_1$ can be ordered as $\{v_1\ll v_n\}$
such that $(v_1\ll v_n)$ has the cross property in the string graph of $C$.
\end{thm}
\Proof Let $G$ be the string graph of $C$.
Choose a point $d\in bd(\Delta)$ such that every member of $C_1$ contains a point of $bd(\Delta)\setminus \{d\}$, 
and for each $x\in C_1$
choose a point $f(x) \in x\cap (bd(\Delta)\setminus \{d\})$. Number $C_1$ so that the points $f(x)\:(x\in C_1)$ are in 
clockwise order, starting from $d$ and breaking ties arbitrarily. Let the numbering of $C_1$ be $\{v_1\ll v_n\}$. 
If $1\le h<i<j<k\le n$, and
$P$ is a path of $G$ between $v_h$ and $v_j$, then the union of the strings in $V(P)$ is an arcwise connected subset of $\Delta$,
containing $f(v_h)$ and $f(v_j)$; and therefore includes a string $s$ with ends $f(v_h)$ and $f(v_j)$ (not necessarily in $C$)
with $s\subseteq \Delta$. Similarly if $Q$ is between $v_i, v_k$, there is a string $t$ between $f(v_i)$ and $f(v_k)$. The strings
$s,t$ intersect, and so one of the strings in $V(P)$ has nonempty intersection with one of the strings in $V(Q)$. 
This proves \ref{discthm}.~\bbox

A {\em homomorphism} from a graph $H$ to a graph $G$ is a map $\eta:V(H)\rightarrow V(G)$,
such that for all adjacent $u,v\in V(H)$,
$\eta(u), \eta(v)$ are distinct and adjacent in $G$.

\begin{thm}\label{crossprop}
Let $G$ be a non-null string graph. Then there is a graph $H$ and $W=\{v_1\ll v_n\}\subseteq V(H)$, such that
\begin{itemize}
\item $(v_1\ll v_n)$ has the cross property in $H$;
\item every vertex in $V(H)\setminus W$ has a neighbour in $W$;
\item there is a homomorphism from $H$ to $G$; and
\item $\chi(H\setminus W)\ge \chi(G)/2$.
\end{itemize}
\end{thm}
\Proof We may assume that $\chi(G)\ge 3$ for otherwise the result is trivial.
Choose a component $D$ of $G$ with maximum chromatic number, and let $z\in D$. For $i\ge 0$ let $L_i$ be the set of vertices of $D$
with distance $i$ from $z$. Choose $k$ such that $\chi(L_k)\ge \chi(G)/2$. Thus $k\ne 0$, and if $k=1$ then let $H$ be
the subgraph induced on $L_0\cup L_1$, and let $n=1$ and $v_1 = z$, and the theorem holds. So we may assume that $k\ge 2$.
Let $D'$ be a component of $G[L_k]$ with maximum chromatic number. The union of the set of strings in $D'$ is a closed arcwise
connected subset of the plane, say $S_1$; and also the union of the strings in $L_0\cup \cdots\cup  L_{k-2}$ is nonnull, 
closed and arcwise connected, say
$S_2$; and $S_1\cap S_2 = \emptyset$. Consequently there is a closed disc $\Delta$ in the plane disjoint from $S_2$ and with $S_1$
in its interior. Moreover, we can choose $\Delta$ such that for each string 
in $L_{k-1}$, its intersection with $\Delta$ is the disjoint union of a finite set of strings. Let $W$ be the set of all strings $s$
such that $s$ is a component of the intersection with $\Delta$
of a string in $L_{k-1}$, and let $H$ be the intersection graph
of the set of strings $W\cup D'$.
For each $s\in W$, we claim that $s\cap bd(\Delta)\ne \emptyset$. For there exists $t\in L_{k-1}$ such that $s$ is a 
component of $t\cap \Delta$; then since $t$ is adjacent in $G$ to a vertex in $L_{k-2}$,
and consequently $t\cap S_2\ne \emptyset$, it follows that every component of $t\cap \Delta$ has nonempty intersection with
$bd(\Delta)$, and in particular, $s\cap bd(\Delta)\ne \emptyset$ as claimed.
The map $\eta:V(H)\rightarrow V(G)$ mapping each string in $V(H)$ to the string in $V(G)$ of which
it is a component, is a homomorphism. Moreover, let $r\in V(H)\setminus W = D'$; we claim that $r$ is adjacent in $H$ to a vertex in
$W$. For let $t\in L_{k-1}$ be adjacent to $r$ in $G$; then $r\cap t\ne\emptyset$, and since $r\subseteq S_1$, 
it follows that $r\cap s\ne \emptyset$ for some $s\in W$. Consequently $r$ is adjacent in
$H$ to a vertex in $W$. 
The result follows from \ref{discthm}. This proves \ref{crossprop}.~\bbox

Finally we need: 
\begin{thm}\label{usecross}
Let $H$ be a graph, let $W\subseteq V(H)$, and let $W=\{v_1\ll v_n\}$ where
$(v_1\ll v_n)$ has the cross property in $H$. Assume also that every vertex in $V(H)\setminus W$ 
has a neighbour in $W$. Then 
$$\chi^3(H)\ge \chi(H\setminus W)/20.$$
\end{thm}
\Proof Let $\kappa = \chi^3(H)$, and suppose that $\chi(H\setminus W)>20\kappa$. 
We may assume that $H$ is connected (by choosing a component of $H$ with maximum chromatic number, and working
inside that). 
For each $i\ge 0$, let $L_i$ be the set of vertices of $H$ with distance exactly $i$ from $v_1$. Choose
$k$ such that $\chi(L_k\setminus W)\ge \chi(H\setminus W)/2$.
Thus $\chi(L_k\setminus W)> 10\kappa$.
Since every vertex in $L_k\setminus W$ has a neighbour in $W$, there are disjoint subsets $X_1\ll X_n$ of $L_k\setminus W$
with union $L_k\setminus W$, such that every vertex in $X_i$
is adjacent to $v_i$ for $1\le i\le n$. Consequently $\chi(X_i)\le \kappa$ for $1\le i\le n$.
\\
\\
(1) {\em There exist $a,b,c,d$ with $1\le a<b<c<d\le n$, such that there is a path of length three between 
$v_a, v_d$, and both its internal vertices belong to $L_k\setminus W$, 
and the subgraph of $H$ induced on $\bigcup_{b\le i\le c} X_i$
has chromatic number more than $4\kappa$.}
\\
\\
For $0\le h\le j\le n$, let $Y(h,j)=\bigcup_{h<i\le j} X_i$. Let $i_0=0$.
Inductively, having defined $i_{j-1}$, choose $i_{j}$ with $i_{j-1}\le i_{j}\le n$ minimal such that 
$\chi(Y(i_{j-1}, i_{j}))> 4\kappa$, if such a choice is possible; and otherwise let $i_{j} = n$ and stop.
Let this process stop with $j=t$ and $i_t=n$ say. For $1\le j< t$, the minimality of $i_j$ implies that 
$\chi(Y(i_{j-1},i_j))\le 5\kappa$, since $\chi(X_{i_j})\le \kappa$. Also $\chi(Y(i_{t-1},i_t))\le 4\kappa$ since the sequence stopped.
Since each of $Y(i_0, i_1), Y(i_1,i_2)\ll Y(i_{t-1}, i_t)$ has chromatic number at most  $5\kappa$, and $\chi(L_k\setminus W)> 10\kappa$,
there exist $h,k$ with $1\le h\le k\le t$ and $h+2\le k$ such that there is an edge between $Y_{i_{h-1},i_{h}}$ and 
$Y_{i_{k-1},i_k}$. Choose $j$ with $h<j<k$; then, taking $b=i_{j-1}+1$ and $c= i_j$, and choosing $a\le i_{j-1}$ 
and $d>i_j$ such that there is an edge between $X_a$ and $X_d$, this proves (1).

\bigskip
Choose $a,b,c,d$ as in (1), and let $Q$ be a path between $v_a,v_d$ of length three.
\\
\\
(2) {\em For each $v\in \bigcup_{b\le i\le c} X_i$, there is a vertex $q$ of $Q$ such that the distance 
between $v,q$ is at most three.}
\\
\\
Since $v\in L_k$, there is a path $P$ between $v_1,v$ of length $k$. Let its vertices be $p_0\d p_1\c p_k$ in order, where
$p_0=v_1$ and $p_k = v$. Choose $e$ with $b\le e\le c$ such that $v$ is adjacent to $v_e$. Then there is a path of $H$
between $v_e,v_1$ with interior included in $V(P)$. 
By the cross property, there is a vertex $q\in V(Q)$ that either belongs to $V(P)\cup \{v_e\}$ 
or has a neighbour in $V(P)\cup \{v_e\}$. Now since the interior vertices of $Q$ belong to $L_k$, it follows that
for $0\le i\le k-3$, $p_i\notin V(Q)$ and has no neighbour in $V(Q)$. So $q$ equals or is adjacent to one of
$p_{k-2}, p_{k-1}, p_k=v, v_e$. In each case the distance between $v,q$ is at most three. This proves (2).

\bigskip
Since the subgraph of $H$ induced on $\bigcup_{b\le i\le c} X_i$
has chromatic number more than $4\kappa$, (2) implies that for one of the four vertices of $Q$, say $q$,
$\chi(N^3[q])>\kappa$, a contradiction. Thus $\chi(H\setminus W)\le 20\kappa$. This proves \ref{usecross}.~\bbox

From \ref{crossprop} and \ref{usecross}, we deduce:

\begin{thm}\label{3string}
For every string graph $G$, $\chi(G)\le 40 \chi^3(G)$.
\end{thm}
\Proof Let $G$ be a string graph, and choose $H$ and $W$ as in \ref{crossprop}. Thus $\chi(H\setminus W)\ge \chi(G)/2$.
By \ref{usecross}, $\chi^3(H)\ge \chi(H\setminus W)/20$, and so $\chi^3(H)\ge \chi(G)/40$. But $\chi^3(G)\ge \chi^3(H)$
since there is a homomorphism from $H$ to $G$. This proves \ref{3string}.~\bbox

In particular, the ideal of string graphs is $3$-controlled. Since no string graph has an induced subgraph which is a proper
subdivision of $K_{3,3}$, \ref{reducecontrol} and \ref{fixclique} imply a result mentioned in section 1, which we restate:

\begin{thm}\label{shortstringagain}
The ideal of string graphs is $2$-controlled.
\end{thm}

Consequently the theorems of this paper apply to string graphs, and in particular, \ref{findchand6} implies a result
mentioned in section 1,
which we restate:
\begin{thm}\label{stringcolouragain}
Let $\nu\ge 0$, and let $H$ be a tree of lamps.
Then there exists $c$ such that every string graph with clique number at most $\nu$ and chromatic number greater than $c$
contains $H$ as an induced subgraph.
\end{thm}

Finally, here is a nice question, raised by Bartosz Walczak (private communication). We proved in \ref{3string} that
$\chi(G)\le 40 \chi^3(G)$ for every string graph $G$, which implies this class is 3-controlled; but we also proved it is 2-controlled.
Is there an analogous result that says $\chi(G)\le K \chi^2(G)$ for every string graph $G$, where $K$ is some constant? 
We think the proofs of this paper
give bounds that are linear if $\omega(G)$ is bounded, but what if $\omega(G)$ is not bounded?

\section{Acknowledgement}
The authors are very grateful for (and impressed by) the excellent and extraordinarily thorough referee report.
They would also like to thank Louis Esperet for discussions on the results of~\cite{chandeliers}, and 
Sean McGuinness for his advice on string graphs.

\end{document}